\documentclass[11pt,reqno]{article}
\usepackage{amssymb, amsmath, amsthm, amsfonts, amscd, epsfig}
\usepackage[english]{babel}
\usepackage{graphicx, epsfig}
\usepackage{color}

\newtheorem{theorem}{Theorem}[section]
\newtheorem{cor}[theorem]{Corollary}
\newtheorem{lemma}[theorem]{Lemma}
\newtheorem{prop}[theorem]{Proposition}

\newcommand{\nm}{\noalign{\smallskip}}
\newcommand{\ds}{\displaystyle}
\newcommand{\pf}{\noindent {\sl Proof}. \ }
\newcommand{\p}{\partial}

\newcommand{\eqnref}[1]{(\ref {#1})}

\newcommand{\Cbb}{\mathbb{C}}

\newcommand{\Ibb}{\mathbb{I}}

\newcommand{\Kbb}{\mathbb{K}}
\newcommand{\Rbb}{\mathbb{R}}
\newcommand{\Sbb}{\mathbb{S}}

\newcommand{\la}{\langle}
\newcommand{\ra}{\rangle}

\newcommand{\Acal}{\mathcal{A}}

\newcommand{\Ecal}{\mathcal{E}}

\newcommand{\Ical}{\mathcal{I}}
\newcommand{\Jcal}{\mathcal{J}}
\newcommand{\Lcal}{\mathcal{L}}
\newcommand{\Kcal}{\mathcal{K}}
\newcommand{\Dcal}{\mathcal{D}}

\newcommand{\Scal}{\mathcal{S}}

\def\Bb{{\bf b}}
\def\Bc{{\bf c}}

\def\Be{{\bf e}}
\def\Bf{{\bf f}}
\def\Bg{{\bf g}}
\def\Bh{{\bf h}}

\def\Bn{{\bf n}}

\def\Bp{{\bf p}}
\def\Bq{{\bf q}}
\def\Br{{\bf r}}
\def\Bs{{\bf s}}
\def\Bt{{\bf t}}
\def\Bu{{\bf u}}
\def\Bv{{\bf v}}
\def\Bw{{\bf w}}
\def\Bx{{\bf x}}
\def\By{{\bf y}}
\def\Bz{{\bf z}}

\def\BH{{\bf H}}
\def\BI{{\bf I}}

\newcommand{\Ga}{\alpha}
\newcommand{\Gb}{\beta}
\newcommand{\Gd}{\delta}
\newcommand{\Ge}{\epsilon}

\newcommand{\Gvf}{\varphi}
\newcommand{\Gg}{\gamma}

\newcommand{\Gk}{\kappa}

\newcommand{\Gl}{\lambda}
\newcommand{\Gn}{\eta}
\newcommand{\Gm}{\mu}
\newcommand{\Gv}{\nu}

\newcommand{\Gt}{\theta}

\newcommand{\Gs}{\sigma}

\newcommand{\Go}{\omega}

\newcommand{\Gy}{\psi}
\newcommand{\Gz}{\zeta}
\newcommand{\GD}{\Delta}

\newcommand{\GG}{\Gamma}

\newcommand{\GO}{\Omega}

\newcommand{\GY}{\Psi}

\newcommand{\BGG}{{\bf \GG}}

\newcommand{\BGvf}{\mbox{\boldmath $\Gvf$}}
\newcommand{\Bpsi}{\mbox{\boldmath $\Gy$}}
\newcommand{\BPsi}{\mbox{\boldmath $\GY$}}
\newcommand{\BGs}{\mbox{\boldmath $\Gs$}}

\newcommand{\beq}{\begin{equation}}
\newcommand{\eeq}{\end{equation}}

\def\ol{\overline}
\newcommand{\hatna}{\widehat{\nabla}}
\numberwithin{equation}{section}
\numberwithin{figure}{section}

\begin{document}

\title{Quantitative characterization of stress concentration in the presence of closely spaced hard inclusions in two-dimensional linear elasticity\thanks{This work is supported by NRF 2016R1A2B4011304 and 2017R1A4A1014735}}

\author{Hyeonbae Kang\thanks{Department of Mathematics, Inha University, Incheon
    22212, S. Korea (hbkang@inha.ac.kr)} \and Sanghyeon Yu\thanks{Seminar for Applied Mathematics, Department of Mathematics, ETH Z\"urich, R\"amistrasse 101, CH-8092 Z\"urich, Switzerland (sanghyeon.yu@sam.math.ethz.ch)}}

\maketitle

\begin{abstract}
In the region between close-to-touching hard inclusions, the stress may be arbitrarily large as the inclusions get closer. The stress is represented by the gradient of a solution to the Lam\'e system of linear elasticity. We consider the problem of characterizing the gradient blow-up of the solution in the narrow region between two inclusions and estimating its magnitude. We introduce singular functions which are constructed in terms of nuclei of strain and hence are solutions of the Lam\'{e} system, and then show that the singular behavior of the gradient in the narrow region can be precisely captured by singular functions. As a consequence of the characterization, we are able to regain the existing upper bound on the blow-up rate of the gradient, namely, $\epsilon^{-1/2}$ where $\epsilon$ is the distance between two inclusions. We then show that it is in fact an optimal bound by showing that there are cases where $\epsilon^{-1/2}$ is also a lower bound. This work is the first to completely reveal the singular nature of the gradient blow-up in the context of the Lam\'{e} system with hard inclusions. The singular functions introduced in this paper play essential roles to overcome the difficulties in applying the methods of previous works. Main tools of this paper are the layer potential techniques and the variational principle. The variational principle can be applied because the singular functions of this paper are solutions of the Lam\'{e} system.
\end{abstract}



\tableofcontents

%%%%%%%%%%%%%%%%%%%%%%%%%%%%%%%%%%%%%%%%%%%%%%
\section{Introduction}
%%%%%%%%%%%%%%%%%%%%%%%%%%%%%%%%%%%%%%%%%%%%%%

When two inclusions are close to touching, the physical field such as the stress or the electric field may be arbitrarily large in the narrow region between the inclusions. It is quite important to understand the field concentration precisely.
Stress concentration may occur in fiber-reinforced composites where elastic inclusions are densely packed \cite{bab}.  The electric field can be greatly enhanced in the conducting inclusions case. It can be utilized to achieve subwavelength imaging and sensitive spectroscopy \cite{YA-SIREV}.

In response to such importance there has been much progress in understanding the field concentration in the last decade or so. In the context of electrostatics (or anti-plane elasticity), the field is the gradient of a solution to the Laplace equation and the precise estimates of the gradient were obtained. It is discovered that when the conductivity of the inclusions is $\infty$, the blow-up rate of the gradient is $\Ge^{-1/2}$ in two dimensions \cite{AKL-MA-05,Yun-SIAP-07}, where $\Ge$ is the distance between two inclusions, and it is $|\Ge \ln \Ge|^{-1}$ in three dimensions \cite{BLY-ARMA-09}. There is a long list of literature in this direction of research, e.g., \cite{AKLLL-JMPA-07, AKLLZ-JDE-09, BLY-CPDE-10, Gorb-MMS-16, GN-MMS-12, Lekner-PRSA-12, LLBY-QAM-14, LY-CPDE-09, LY-JDE-11, Yun-JMAA-09, Yun-JDE-16}.
While these works are related to the estimate of the blow-up rate of the gradient, there is other direction of research to characterize the singular behavior of the gradient \cite{ACKLY-ARMA-13, KLY-MA-15, KLY-JMPA-13, KLY-SIAP-14, LY-JMAA-15}. An explicit function, which is called a singular function, is introduced and the singular behavior of the gradient is completely characterized by this singular function. Since the singular function is closely related to this work, we include a brief discussion on it at the beginning of subsection \ref{subsec:singular}. All the work mentioned above are related to the homogenous equation and inclusions with smooth boundaries. Recently there have been important extensions to the inhomogeneous equation \cite{DL-arXiv} and inclusions with corners (the bow-tie shape) \cite{KYun}.

In this paper, we consider a similar problem in linear elasticity, i.e., the Lam\'e system.  
We assume two hard inclusions, which have infinite shear modulus, are presented with a small separation distance $\Ge$. The stress is represented in terms of the gradient of a solution to the Lam\'e system. 
We are interested in the singular behavior of the stress (or the gradient) when the distance $\Ge$ goes to zero.

Even though much progress has been made for the Laplace equation of the anti-plane elasticity as mentioned above, not much is known about the gradient blow-up in the context of the full elasticity, e.g., the Lam\'e system. Recently, a significant progress has been made by Bao {\it et al} \cite{BLL-ARMA-15, BLL-arXiv}: it is proved in \cite{BLL-ARMA-15} that $\Ge^{-1/2}$ is an upper bound on the blow-up rate of the gradient in the two-dimensional Lam\'e system. We emphasize that there is significant difficulty in applying the methods for scalar equations to systems of equations. For instance, the maximum principle does not hold for the system. In \cite{BLL-ARMA-15} they come up with an ingenious iteration technique to overcome this difficulty and obtain the upper bound on the blow-up rate. However, it was still not known if it is also a lower bound. 

The purpose of this paper is to construct singular functions for the two-dimensional Lam\'e system, like the one for electrostatics, and to characterize the singular behavior of the gradient using singular functions. In fact, we construct singular functions as elaborated linear combinations of nuclei of strain, and show that they capture the singular behavior of the gradient precisely. Nuclei of strain are the columns of the Kelvin matrix of the fundamental solution to the Lam\'e system and their variants. As a consequence of such characterization, we are able to reobtain the result of \cite{BLL-ARMA-15} with a different proof, which states that  $\Ge^{-1/2}$ is an upper bound on the blow-up rate of the gradient. More importantly, the characterization enables us to show that the rate $\Ge^{-1/2}$ is actually optimal, optimal in the sense that there are cases where $\Ge^{-1/2}$ is a lower bound on the blow-up rate. 
To the best of our knowledge, this work is the first to completely reveal the singular nature of the gradient blow-up  in the context of the Lam\'e system with hard inclusions. The singular functions introduced in this paper play essential roles to overcome the difficulties in applying the methods of previous works.

We emphasize that the nuclei of strain and singular functions are solutions of the Lam\'e system. This has a significant implication. We heavily use the variational principle for proving the characterization of the stress concentration in section \ref{sec:bvp}, which is only possible since singular functions are solutions of the Lam\'e system. This makes the method of this paper significantly different from that of \cite{BLL-ARMA-15}. We include a brief comparison of two methods at the end of subsection \ref{subsec:bvp}.

It is worth mentioning that the singular functions constructed in this paper are applied to an important problem other than analysis of the gradient blow-up. In fact, quantitative analysis of the gradient is closely related to the computation of the effective property of densely packed 
composites. In \cite{BK-ARMA-01}, Beryland et al. provided the first rigorous justification of the asymptotic formula for the effective conductivity, which was found by Keller \cite{Keller-JAP-63}. However, the corresponding formulas of Flaherty-Keller \cite{FK-CPAM-73} for the effective elastic properties have not been rigorously proved to the best of our knowledge. Using singular functions of this paper we are able to prove the formulas in a mathematically rigorous way. We emphasize that this is possible only because singular functions are solutions of the Lam\'e system. We report this result in a separate paper \cite{KY}.

Accurate numerical computation of the gradient in the presence of closely spaced hard inclusions is a well-known challenging problem in computational mathematics and sciences. When computing the gradient, a serious difficulty arises since a fine mesh is required to capture the gradient blow-up in the narrow region. The precise characterization of the gradient blow-up can be utilized for designing an efficient numerical scheme to compute the gradient. This was done for the conductivity case in \cite{KLY-JMPA-13}. The result of  this paper may open up a way to do computation for the isotropic elasticity.

It is worth mentioning that, for the Lam\'e system where two inclusions are circular holes, the gradient blow-up is recently characterized by a singular function in \cite{LY-arXiv}. Moreover, the optimal blow-up rate of the gradient is obtained. The holes are characterized by the vanishing traction condition on the boundary, and the blow-up rate is the same as the hard inclusion case, namely, $\Ge^{-1/2}$. We emphasize that unlike the anti-plane elasticity, the hole case is not the dual problem of the hard inclusion case, and a different method is required to handle the hole case.

This paper consists of six sections including introduction and appendices. In section \ref{sec:prob}, we formulate the problem to be considered, derive some preliminary results which will be used in later sections, and describe geometry of two inclusions. In section \ref{sec:singular}, singular functions are constructed in terms of nuclei of strain and their properties are derived for later use. Section \ref{sec:bvp} and \ref{sec:free} respectively deal with the problem of characterizing the stress concentration in a bounded domain and in the free space. In section \ref{sec:symmetric_case} we consider the case when inclusions are symmetric, in particular, when inclusions are disks of the same radius, and show that $\Ge^{-1/2}$ is a lower bound on the blow-up rate of the gradient when Lam\'e constants satisfy a certain constraint. Since each section is rather long and its subject can be viewed as independent, we include an introduction in each section. Appendices are to prove some results used in the text, especially existence and uniqueness of the solution to the exterior problem of the Lam\'e system and the layer potential representation of the solution to the boundary value problem and the free space problem.

Throughout this paper, we use the expression $A \lesssim B$ to imply that there is a constant $C$ independent of $\Ge$ such that $A \le CB$. The expression $A \approx B$ implies that both $A \lesssim B$ and $B \lesssim A$ hold.

%%%%%%%%%%%%%%%%%%%%%%%%%%%%%%%%%%%%%%%%%%%%%%
\section{Problem formulation and preliminaries}\label{sec:prob}
%%%%%%%%%%%%%%%%%%%%%%%%%%%%%%%%%%%%%%%%%%%%%%

In this section we formulate the problem of characterizing the stress concentration. The main tools in dealing with the problem are the layer potential technique and the variational principle. We introduce them in this section. We then consider the existence and uniqueness question of the exterior problem for the Lam\'e system with arbitrary Dirichlet data. The final subsection is to describe the geometry of two inclusions in a precise manner.

%%%%%%%%%%%%%%%%%%%%%%%%%%%%%%%%%%%%%%%%%%%%%%%%%%%%%%%%%%%%%%%%%%
\subsection{Lam\'{e} system with hard inclusions: a problem formulation}
%%%%%%%%%%%%%%%%%%%%%%%%%%%%%%%%%%%%%%%%%%%%%%%%%%%%%%%%%%%%%%%%%%

We consider two disjoint elastic inclusions $D_1$ and $D_2$ which are embedded in $\Rbb^2$ occupied by an elastic material. We assume that $D_1$ and $D_2$ are simply connected bounded domains with $C^4$-smooth boundaries. We emphasize that the results of this paper are valid even if boundaries are $C^{3,\Ga}$ for some $\Ga>0$. But we assume that they are $C^4$ for convenience. Advantage of assuming $C^4$ is made clear in subsection \ref{subsection:geo_two_incl}. We also assume some convexity of the boundaries which is precisely described in the same subsection.

Let $(\Gl,\Gm)$ be the pair of Lam\'{e} constants of $D^e:= \Rbb^2 \setminus \ol{D_1 \cup D_2}$ which satisfies the strong ellipticity conditions
$\Gm>0$ and $\Gl+\Gm>0$. Then the elasticity tensor is given by $\Cbb=(C_{ijkl})$ with
$$
C_{ijkl} = \Gl \Gd_{ij} \Gd_{kl} + \Gm (\Gd_{ik} \Gd_{jl} + \Gd_{il} \Gd_{jk}), \quad i,j,k,l=1,2,
$$
where $\Gd_{ij}$ denotes Kronecker's delta. The Lam\'e operator $\Lcal_{\Gl,\Gm}$ of the linear isotropic elasticity is defined by
\beq
\Lcal_{\Gl,\Gm} \Bu:= \nabla \cdot \Cbb \hatna\Bu = \Gm \GD \Bu+(\Gl+\Gm)\nabla \nabla \cdot \Bu,
\eeq
where $\hatna$ denotes the symmetric gradient, namely,
$$
\hatna \Bu = \frac{1}{2} \left( \nabla \Bu +\nabla \Bu^T\right) \quad \text{($T$ for transpose)}.
$$
The corresponding conormal derivative $\p_\Gv \Bu$ on $\p D_j$ is defined as
\beq
\p_\Gv \Bu = (\Cbb \hatna\Bu) \Bn,
\eeq
where $\Bn$ is the outward unit normal vector to $\p D_j$ ($j=1,2$).

Given a displacement field $\Bu=(u_1,u_2)^T$, $\hatna \Bu$ is the strain tensor while the stress tensor $\BGs = (\Gs_{ij})_{i,j=1}^2$ is defined to be
\beq\label{def_stress}
\BGs := \Cbb \hatna\Bu = \Gl \mbox{tr} ( \hatna \Bu)\BI +2\Gm \hatna \Bu,
\eeq
namely,
\begin{align}
\Gs_{11} &= (\Gl+2\Gm) \p_1 u_1 + \Gl\p_2 u_2,
\nonumber
\\
\Gs_{22} &= \Gl \p_1 u_1 + (\Gl+2\Gm) \p_2 u_2, \label{stress_cartesian}
\\
\Gs_{12} &= \Gs_{21} = \Gm (\p_2 u_1 + \p_1 u_2). \nonumber
\end{align}
Here and throughout this paper, tr stands for the trace and $\p_j$ denotes the partial derivative with respect to the $x_j$-variable for $j=1,2$.

Let $\BPsi$ be the collection of all functions $\Bpsi$ such that $\hatna \Bpsi = 0$ in $\Rbb^2$, {\it i.e.}, the three-dimensional vector space spanned by
the displacement fields of the rigid motions $\{\Psi_j\}_{j=1}^3$ defined as follows:
\beq\label{Psidef}
\Psi_1(\Bx)=\begin{bmatrix} 1\\0 \end{bmatrix}, \quad
\Psi_2(\Bx)=\begin{bmatrix} 0\\1 \end{bmatrix}, \quad
\Psi_3(\Bx)=\begin{bmatrix} -y\\ x \end{bmatrix}.
\eeq
Throughout this paper we denote the point $\Bx$ in $\Rbb^2$ by either $(x_1,x_2)^T$ or $(x,y)^T$ at its convenience.

We assume $D_1$ and $D_2$ are hard inclusions. This assumption is inscribed on the boundary conditions on $\p D_j$ in the following problem: Let $\GO$ be a bounded domain in $\Rbb^2$ containing $D_1$ and $D_2$ such that $\mbox{dist}(\p \GO, D_1 \cup D_2) \geq C$ for some constant $C>0$.
Let us denote
$$
\widetilde{\GO}=\GO\setminus\overline{D_1\cup D_2}.
$$
For a given Dirichlet data $\Bg$ we consider the following problem:
\beq\label{elas_eqn_bdd}
 \ \left \{
 \begin{array} {ll}
\ds \Lcal_{\Gl,\Gm} \Bu= 0 \quad &\mbox{ in } \widetilde{\GO},\\[2mm]
\ds \Bu=\sum_{j=1}^3 c_{ij} \Psi_j(\Bx) \quad &\mbox{ on } \p D_i, \quad i=1,2 ,\\[2mm]
\ds \Bu = \mathbf{g} \quad& \mbox{ on } \p \GO,
 \end{array}
 \right.
 \eeq
where the constants $c_{ij}$ are determined by the conditions
\beq\label{int_zero}
\int_{\p D_i} \p_\Gv \Bu |_+ \cdot \Psi_j \,d\Gs =0, \quad i=1,2, \, j=1,2,3.
\eeq
Here and afterwards, the subscript $+$ denotes the limit from outside $\p D_j$.

Let
\beq
\Ge:= \mbox{dist}(D_1, D_2).
\eeq
The gradient $\nabla \Bu$ of the solution $\Bu$ to \eqnref{elas_eqn_bdd} may become arbitrarily large as two inclusions get closer, namely, as $\Ge \to 0$. The main purpose of this paper is to characterize the blow-up of $\nabla \Bu$. Roughly speaking, we show that $\Bu$ can be decomposed as
\beq\label{BuBsBb}
\Bu=\Bs+ \Bb,
\eeq
where $\nabla \Bs$ has the main singularity of $\nabla \Bu$ while $\nabla \Bb$ is regular or less singular. So the singular behavior of $\nabla \Bu$ is characterized by that of $\nabla \Bs$. We will find $\Bs$ in an explicit form. The characterization of the gradient blow-up enables us to show that the optimal blow-up rate of $\nabla \Bu$ in terms of $\Ge$ is $\Ge^{-1/2}$. It is proved in \cite{BLL-ARMA-15} that $\Ge^{-1/2}$ is an upper bound on the blow-up rate of $\nabla \Bu$ as mentioned before.

The problem in the presence of hard inclusions may be considered as the limiting problem of a high contrast elasticity problem when the shear modulus of the inclusions degenerates to infinity \cite{BLL-ARMA-15}. When the shear modulus is bounded away from zero and infinity, it is known that the gradient is bounded regardless of the distance between inclusions \cite{LN}.

We also consider the free space problem: For a given function $\BH$ satisfying $\Lcal_{\Gl,\Gm} \BH= 0 \mbox{ in } \Rbb^2$, the displacement field $\Bu$ satisfies
\beq\label{elas_eqn_free}
 \ \left \{
 \begin{array} {ll}
\ds \Lcal_{\Gl,\Gm} \Bu= 0 \quad &\mbox{in } D^e,\\[2mm]
\ds
 \Bu=\sum_{j=1}^3 d_{ij} \Psi_j \quad &\mbox{on } {\p D_i}, \quad i=1,2 ,
\\[2mm]
\ds \Bu(\Bx)-\BH(\Bx) = O(|\Bx|^{-1}) \quad& \mbox{as } |\Bx| \rightarrow \infty,
 \end{array}
 \right.
 \eeq
where the constants $d_{ij}$  are determined by the condition \eqnref{int_zero}. We will obtain the decomposition of the form \eqnref{BuBsBb} and estimates of $\nabla \Bu$ for this problem as well.

%%%%%%%%%%%%%%%%%%%%%%%%%%%%%%%%%%%%%%%%%%%%%%%%%%%%%%%%%%
\subsection{Layer potentials for 2D Lam\'{e} system}
%%%%%%%%%%%%%%%%%%%%%%%%%%%%%%%%%%%%%%%%%%%%%%%%%%%%%%%%%%

The Kelvin matrix of fundamental solutions $\BGG = \left( \GG_{ij} \right)_{i, j = 1}^2$ to the Lam\'{e} operator $\Lcal_{\Gl, \Gm}$ is given by
\beq\label{Kelvin}
  \GG_{ij}(\Bx) =
    \ds \Ga_1 \Gd_{ij} \ln{|\Bx|} - \Ga_2 \displaystyle \frac{x_i x_j}{|\Bx| ^2}
\eeq
where
\beq
  \Ga_1 = \frac{1}{4\pi} \left( \frac{1}{\Gm} + \frac{1}{\Gl + 2 \Gm} \right) \quad\mbox{and}\quad
  \Ga_2 = \frac{1}{4\pi} \left( \frac{1}{\Gm} - \frac{1}{\Gl + 2 \Gm} \right).
\eeq
In short, $\BGG$ can be expressed as
\beq\label{Kelvintensor}
\BGG(\Bx) = \Ga_1 \ln|\Bx-\By| \BI - \Ga_2  \Bx \otimes \nabla(\ln |\Bx|),
\eeq
where $\BI$ is the identity matrix.

For a given bounded domain $D$ with $C^2$ boundary, the single and double layer potentials on $\p D$ associated with the pair of Lam\'{e} parameters $(\Gl,\Gm)$ are defined by
\begin{align}
& \Scal_{\p D} [\BGvf] (\Bx) := \int_{\p D} \BGG  (\Bx-\By) \BGvf(\By) d \Gs(\By), \quad \Bx \in \Rbb^2,\\
& \Dcal_{\p D} [\BGvf] (\Bx) := \int_{\p D} \p_\Gv \BGG  (\Bx-\By) \BGvf(\By) d \Gs(\By), \quad \Bx \in \Rbb^2 \setminus \p D,
\end{align}
where the conormal derivative $\p_\Gv \BGG (\Bx-\By)$ is defined by
$$
\p_\Gv \BGG  (\Bx-\By) \Bb = \p_\Gv (\BGG  (\Bx-\By)\Bb)
$$
for any constant vector $\Bb$.

Let $H^{1/2}(\p D)$ be the usual $L^2$-Sobolev space of order $1/2$ on $\p D$ and $H^{-1/2}(\p D)$ be its dual space. With functions $\Psi_j$ in \eqnref{Psidef} we define
\beq
H^{-1/2}_\Psi (\p D) := \{ \Bf\in H^{-1/2}(\p D)^2 : \int_{\p D} \Bf\cdot \Psi_j =0, \,j=1,2,3 \}.
\eeq

The following propositions for representations of the solutions to \eqnref{elas_eqn_bdd} and \eqnref{elas_eqn_free} can be proved in a standard way (see, for example, \cite{AK-book-07}). We include brief proofs in Appendix.

\begin{prop}\label{thm_u_layer_bdd_case}
Let $\Bu$ be the solution to \eqnref{elas_eqn_bdd} and let $\Bf:= \p_\Gv \Bu|_-$ on $\p\GO$. Define
\beq\label{eqn_def_H_Omega}
\BH_{\GO}(\Bx) =-\Scal_{\p\GO}[\Bf](\Bx)+ \Dcal_{\p\GO}[\Bg](\Bx), \quad \Bx \in \GO.
\eeq
Then there is a unique pair $(\BGvf_1,\BGvf_2)\in H^{-1/2}_\Psi (\p D_1) \times H^{-1/2}_\Psi(\p D_2)$ such that
\beq\label{rep_bdd}
\Bu(\Bx) = \BH_\GO(\Bx) + \Scal_{\p D_1}[\BGvf_1](\Bx) + \Scal_{\p D_2}[\BGvf_2](\Bx), \quad \Bx \in \GO.
\eeq
In fact, $\BGvf_j$ is given by $\BGvf_j = \p_\Gv \Bu|_+$ on $\p D_j$ for $j=1,2$.
\end{prop}

\begin{prop}\label{thm_u_layer}
Let $\Bu$ be the solution to \eqnref{elas_eqn_free}. Then there is a unique pair $(\BGvf_1,\BGvf_2)\in H^{-1/2}_\Psi (\p D_1) \times H^{-1/2}_\Psi(\p D_2)$ such that
\beq\label{singlerep}
\Bu(\Bx) = \BH(\Bx) + \Scal_{\p D_1}[\BGvf_1](\Bx) + \Scal_{\p D_2}[\BGvf_2](\Bx), \quad \Bx \in D^e.
\eeq
In fact, $\BGvf_j$ is given by $\BGvf_j = \p_\Gv \Bu|_+$ on $\p D_j$ for $j=1,2$.
\end{prop}
Note that $\int_{\p D_j} \BGvf_j =0$, which holds because $\BGvf_j$ belongs to $H^{-1/2}_\Psi (\p D_j)$. So, we have $\Scal_{\p D_j}[\BGvf_j](\Bx) = O(|\Bx|^{-1})$ as $|\Bx| \to \infty$. Thus $\Bu$ given by \eqnref{singlerep} satisfies the last condition in \eqnref{elas_eqn_free}.

Note that since the domains $D_1$ and $D_2$ are assumed to have $C^{4}$ boundaries, the solutions to \eqnref{elas_eqn_bdd} and \eqnref{elas_eqn_free} are $C^{3,\Ga}$ in $\GO \setminus (D_1 \cup D_2)$ including $\p D_1 \cup \p D_2$ for any $0<\Ga <1$.

We now prove an analogue of the addition formula for $\BGG (\Bx-\By)$. Let $\{ \Be_1, \Be_2 \}$ be the standard basis for $\Rbb^2$. For $n\in \mathbb{Z}$ let
\beq
P_n(\Bx) = r^{|n|} e^{in \Gt},
\eeq
where  $(r, \Gt)$ denotes the polar coordinates of $\Bx$. Let
\begin{align}
\Bv_n^{(i)}(\Bx) &= \Ga_1 P_n(\Bx) \Be_i -\Ga_2 x_i \nabla P_n(\Bx), \quad i=1,2, \\
\Bw_n(\Bx) &= \Ga_2 \nabla P_n (\Bx).
\end{align}
Since $P_n$ is harmonic in $\Rbb^2$, one can easily see that $\Bw_n$ is a solution to the Lam\'{e} system in $\Rbb^2$. To show that $\Bv_n^{(i)}$ is a solution to the Lam\'{e} system in $\Rbb^2$, we prove a more general fact:
\begin{lemma}\label{harLame}
If $h$ is a harmonic function, then a vector-valued function $\Bv$ of the form
\beq
\Bv(\Bx) = \Ga_1 h(\Bx)\Be_j-\Ga_2 x_j \nabla h(\Bx)
\eeq
for $j=1, 2$, is a solution of the Lam\'e system, namely, $\Lcal_{\Gl, \Gm} \Bv=0$.
\end{lemma}
\pf
We only prove the case when $j=1$.
Let us write $\Bv=(v_1, v_2)^T$. Simple computations show that
$$
\GD v_1 = -2\Ga_2 \p_1^2 h,
$$
and
$$
\GD v_2 = -2\Ga_2 \p_1 \p_2 h.
$$
We also have
$$
\nabla\cdot \Bv = \Ga_1 \p_1 h -\Ga_2 (x_1 \GD h + \p_1 h)
=(\Ga_1-\Ga_2) \p_1 h.
$$
Therefore we obtain
\begin{align*}
[\mu \GD \Bv + (\Gl+\mu) \nabla(\nabla\cdot \Bv)]\cdot\Be_k &=
-2 \mu \Ga_2 \p_1\p_k h + (\Gl+\mu)(\Ga_1-\Ga_2) \p_1\p_k h
\\
&=
\Big(-\frac{1}{2\pi} \frac{\Gl+\mu}{\Gl+2\mu}+\frac{1}{2\pi} \frac{\Gl+\mu}{\Gl+2\mu}\Big)\p_1 \p_k h
=0
\end{align*}
for $k=1,2$. This completes the proof.
\qed

We obtain the following proposition.
\begin{prop}\label{thm:funda_sol_series}
The fundamental solution $\BGG$ admits the following series expansion: for $|\Bx| > |\By|$ and for any constant vector $\Bb$ in $\Rbb^2$
\begin{align}
\BGG(\Bx-\By) \Bb &= -\sum_{n\neq 0}\frac{1}{2|n|}\frac{e^{-in\Gt}}{r^{|n|}} \sum_{i=1}^2(\Bv_n^{(i)}(\By) \cdot \Bb)\Be_i
\nonumber \\
& \quad + \sum_{n\neq 0}\frac{1}{2|n|}\frac{\Bx e^{-in\Gt }}{r^{|n|}} \big(\Bw_n(\By)\cdot \Bb\big)
 +\Ga_1 \ln|\Bx| \Bb, \label{addition}
\end{align}
where $\Bx=(r, \Gt)$ in the polar coordinates.
Moreover, the series converges absolutely and uniformly in $\Bx$ and $\By$ provided that there are numbers $r_1$ and $r_2$ such that $|\By| \le r_1 < r_2 \le |\Bx|$.
\end{prop}
\pf
By \eqnref{Kelvintensor}, we have
$$
\BGG(\Bx-\By)\Bb = \Ga_1 \ln|\Bx-\By|\Bb - \Ga_2 (\nabla_\By(\ln |\Bx-\By|)  \cdot \Bb) (\Bx-\By)(-1)
$$
for any constant vector $\Bb$. The addition formula for $\ln|\Bx-\By|$ reads
$$
\ln|\Bx-\By| = \ln|\Bx| - \sum_{n\neq 0} \frac{1}{2|n|}\frac{e^{-in\Gt}}{r^{|n|}} P_n(\By).
$$
By substituting this formula to the one above, we obtain \eqnref{addition}.
\qed

%%%%%%%%%%%%%%%%%%%%%%%%%%%%%%%%%%%%%%%%%%%%%%%%%%%%%%%%%%%%%%%%%%%%%%%%%%%
\subsection{The exterior problem and the variational principle}
%%%%%%%%%%%%%%%%%%%%%%%%%%%%%%%%%%%%%%%%%%%%%%%%%%%%%%%%%%%%%%%%%%%%%%%%%%%

In this subsection we consider the following exterior Dirichlet problem for the Lam\'e system:
\beq\label{eqn_ext_diri}
\left \{
\begin{array} {ll}
\ds \Lcal_{\Gl,\Gm} \mathbf{v}= 0 \quad &\mbox{ \rm in } D^e,\\[2mm]
\ds \mathbf{v}= \mathbf{g} \quad &\mbox{ \rm on } \p D^e=\p D_1 \cup \p D_2,
 \end{array}
 \right.
\eeq
for $\Bg \in H^{1/2}(\p D^e)^2:= H^{1/2}(\p D_1)^2 \times H^{1/2}(\p D_2)^2$. We seek a solution in the function space $\Acal^*$ defined as follows: Let $\Acal$ be the collection of all $\Bv \in H^1_\text{loc}(D^e)$ such that there exists a $2 \times 2$ symmetric matrix $B$ such that
\beq\label{2002}
\Bv(\Bx)= \sum_{j=1}^2 \p_j \BGG(\Bx) B \Be_j + O(|\Bx|^{-2}) \quad\mbox{as } |\Bx| \to \infty,
\eeq
where $\{\Be_1, \Be_2 \}$ is the standard basis of $\Rbb^2$. We emphasize that $\Bv(\Bx)=O(|\Bx|^{-1})$ as $|\Bx| \to \infty$. We then define
\beq
\Acal^* := \left\{ \Bu = \Bv + \sum_{j=1}^3 b_j \Psi_j ~|~ \Bv \in \Acal, \ \ b_j: \text{constant} \right\}.
\eeq

A proof of the following theorem is given in Appendix.

\begin{theorem}\label{thm_ext_diri}
For any $\Bg \in H^{1/2}(\p D^e)^2$, \eqnref{eqn_ext_diri} admits a unique solution in $\Acal^*$.
\end{theorem}
This theorem in a different form is proved in \cite{Constanda} when $D^e$ is the compliment of a simply connected domain. Here, $\p D^e$ has two components, namely, $\p D^e = \p D_1 \cup \p D_2$. Moreover, the proof of this paper is completely different from that of \cite{Constanda}. It is worth mentioning that the term $\sum_{j=1}^3 b_j \Psi_j$ plays the role of the solution corresponding to the component of $\Bg$ spanned by $\Psi_j$, $j=1,2,3$.

The condition \eqnref{2002} is somewhat unfamiliar. To motivate it we prove the following lemma. This lemma will be used in the proof of Theorem \ref{thm_ext_diri}.

\begin{lemma}\label{lem:Acal}
\begin{itemize}
\item[(i)] If $\BGvf=(\BGvf_1, \BGvf_2) \in H^{-1/2}(\p D_1)^2 \times H^{-1/2}(\p D_2)^2$ and satisfies
\beq\label{vanishingcon}
\int_{\p D_1} \BGvf_1 \cdot \Psi_k + \int_{\p D_2} \BGvf_2 \cdot \Psi_k =0, \quad k=1,2,3,
\eeq
then $\Bv$, defined by
\beq\label{Bvtemp}
\Bv(\Bx)= \Scal_{\p D_1} [\BGvf_1](\Bx) + \Scal_{\p D_2} [\BGvf_2](\Bx), \quad \Bx \in D^e,
\eeq
belongs to $\Acal$.
\item[(ii)] If $\Bpsi=(\Bpsi_1, \Bpsi_2)$ belongs to $H^{1/2}(\p D_1)^2 \times H^{1/2}(\p D_2)^2$, then $\Bw$, defined by
\beq
\Bw(\Bx)= \Dcal_{\p D_1} [\Bpsi_1](\Bx) + \Dcal_{\p D_2} [\Bpsi_2](\Bx), \quad \Bx \in D^e,
\eeq
belongs to $\Acal$.
\end{itemize}
\end{lemma}
\pf
If $\By \in \p D^e$ and $|\Bx| \to \infty$, then by the Taylor expansion we have
\beq\label{Taylor}
\BGG(\Bx-\By) = \BGG(\Bx) + \sum_{j=1}^2 \p_j \BGG(\Bx) y_j + O(|\Bx|^{-2}).
\eeq
So $\Bv$ defined by \eqnref{Bvtemp} takes the form
\beq\label{2000}
\Bv(\Bx)= \BGG(\Bx) \int_{\p D^e} \BGvf + \sum_{j=1}^2 \p_j \BGG(\Bx) \int_{\p D^e} y_j \BGvf + O(|\Bx|^{-2}).
\eeq
Here and throughout this paper we use $\int_{\p D^e} \BGvf$ to denote $\int_{\p D_1} \BGvf_1 + \int_{\p D_2} \BGvf_2$ for ease of notation. The assumption \eqnref{vanishingcon} for $k=1,2$ implies that the first term in the right-hand side of \eqnref{2000} vanishes. Define the matrix $B:= (b_{ij})_{i,j=1,2}$ by
$$
\begin{bmatrix} b_{11} \\ b_{21} \end{bmatrix} := \int_{\p D^e} y_1 \BGvf \quad\mbox{and}\quad
\begin{bmatrix} b_{12} \\ b_{22} \end{bmatrix} := \int_{\p D^e} y_2 \BGvf.
$$
Then, we may rewrite \eqnref{2000} as
$$
\Bv(\Bx)= \sum_{j=1}^2 \p_j \BGG(\Bx) B \Be_j + O(|\Bx|^{-2}) \quad\mbox{as } |\Bx| \to \infty.
$$
Note that the assumption \eqnref{vanishingcon} for $k=3$ implies $b_{12}=b_{21}$, namely, $B$ is symmetric.

To prove (ii), let $\Bu_j$ be the solution to $\Lcal_{\Gl, \Gm} \Bu_j =0$ in $D_j$ and $\Bu_j= \Bpsi_j$ on $\p D_j$. Then $\p_\Gv \Bu_j \in H^{-1/2}_\Psi (\p D_j)$ and Green's formula for the Lam\'e system shows that the following holds:
$$
\Dcal_{\p D_1} [\Bpsi_1](\Bx)= \Scal_{\p D_1} [\p_\Gv \Bu_1](\Bx), \quad \Dcal_{\p D_2} [\Bpsi_2](\Bx)= \Scal_{\p D_2} [\p_\Gv \Bu_2](\Bx), \quad \Bx \in D^e.
$$
So, we have
$$
\Bw(\Bx)= \Scal_{\p D_1} [\p_\Gv \Bu_1](\Bx) + \Scal_{\p D_2} [\p_\Gv \Bu_2](\Bx), \quad \Bx \in D^e.
$$
So, (ii) follows from (i).
\qed

The most important property of the function of the form $\sum_{j=1}^2 \p_j \BGG(\Bx) B \Be_j$ lies in the following fact.
\begin{lemma}\label{Bfunction}
Let $\Bv(\Bx)= \sum_{j=1}^2 \p_j \BGG(\Bx) B \Be_j$ for some symmetric matrix $B$. Then the following holds for any simple closed Lipschitz curve $C$ such that $0 \notin C$:
\beq\label{crux}
\int_{C} \p_\Gv \Bv \cdot \Psi_k=0, \quad k=1,2,3.
\eeq
\end{lemma}
\pf
Since the cases of $k=1,2$ are easier to prove, we only consider the case of $k=3$. Let $U$ be the bounded domain enclosed by $C$. If $0 \notin U$, then by Green's formula for the Lam\'e system, we have
$$
\int_{C} \p_\Gv \Bv \cdot \Psi_3 = \int_U \Cbb \hatna \Bv: \hatna \Psi_3=0.
$$
Suppose that $0 \in U$. Then choose $B_r$, the disk of radius $r$ centered at $0$, so that $\overline{B_r} \subset U$. Then, we see that
$$
\int_{C} \p_\Gv \Bv \cdot \Psi_3 = \int_{\p B_r} \p_\Gv \Bv \cdot \Psi_3.
$$
Straightforward but tedious computations show that on $\p B_r$
$$
\p_\Gv \Bv \cdot \Psi_3 = \frac{1}{2\pi r} (b_{21}-b_{12})+   \frac{\Gl+\mu}{\Gl+2\mu}(b_{11} -\frac{b_{12}}{2}-\frac{b_{21}}{2} - b_{22})\frac{\sin 2\theta}{2\pi r},
$$
where $(r,\theta)$ is the polar coordinates.
So we obtain
$$
\int_{\p B_R} \p_\Gv \Bv \cdot \Psi_3  = b_{21}-b_{12}.
$$
Since $b_{12}=b_{21}$, \eqnref{crux} follows. \qed

The following lemma shows that Green's formula holds for $\Bu, \Bv\in \Acal^*$. It is worth mentioning that the $-$ sign appears on the right-hand side of \eqnref{231} below since the normal vector on $\p D^e$ is directed outward.

\begin{lemma}\label{cor_betti}
If $\Bu, \Bv\in \Acal^*$ and $\Lcal_{\Gl,\Gm}\Bu=0$ in $D^e$, then
\beq\label{231}
\int_{D^e}\Cbb \hatna \Bu: \hatna \Bv  =-\int_{\p D^e} \p_\Gv \Bu|_+ \cdot \Bv,
\eeq
where the left-hand side is understood to be
\beq
\int_{D^e}\Cbb \hatna \Bu: \hatna \Bv  = \lim_{R \to \infty} \int_{B_R \setminus (D_1 \cup D_2)}\Cbb \hatna \Bu: \hatna \Bv.
\eeq
\end{lemma}

\pf
We have
$$
\int_{B_R \setminus (D_1 \cup D_2)}\Cbb \hatna \Bu: \hatna \Bv
= -\int_{\p D^e} \p_\Gv \Bu|_+ \cdot \Bv + \int_{\p B_R} \p_\Gv \Bu|_+ \cdot \Bv.
$$
So, it suffices to prove that
$$
\lim_{R \to \infty} \int_{\p B_R} \p_\Gv \Bu|_+ \cdot \Bv =0.
$$

Let $\Bu= \Bu_1+\Bu_2+\Bu_3$ where $\Bu_1$ is of the form $\sum_{j=1}^2 \p_j \BGG(\Bx) B \Be_j$, $\Bu_2(\Bx)=O(|\Bx|^{-2})$, and $\Bu_3$ is of the form $\sum_{k=1}^3 a_k \Psi_k$. We also let $\Bv= \Bv_1+\Bv_2$ where $\Bv_1(\Bx)=O(|\Bx|^{-1})$ and $\Bv_2$ is of the form $\sum_{k=1}^3 b_k \Psi_k$. Since $\p_\Gv \Bu_3=0$ on $\p B_R$ for any $R$, we have
$$
\int_{\p B_R} \p_\Gv \Bu|_+ \cdot \Bv = \int_{\p B_R} \p_\Gv (\Bu_1+\Bu_2)|_+ \cdot (\Bv_1+\Bv_2).
$$
We see
$$
\lim_{R \to \infty} \left[ \int_{\p B_R} \p_\Gv (\Bu_1+\Bu_2)|_+ \cdot \Bv_1 + \int_{\p B_R} \p_\Gv \Bu_2|_+ \cdot \Bv_2 \right]=0
$$
by considering the decay at $\infty$ of the functions involved. We also have from \eqnref{crux}
$$
\int_{\p B_R} \p_\Gv \Bu_1 |_+ \cdot \Bv_2 =0.
$$
This completes the proof. \qed

The following variational principle for the exterior Dirichlet problem plays a crucial role in what follows.
\begin{lemma}\label{lem_var_principle}
Define
\beq\label{eqn_def_Ecal}
\Ecal_{D^e}[\Bw]:= \int_{D^e} \Cbb \hatna \Bw :\hatna\Bw.
\eeq
Let $\Bu$ be the solution in $\Acal^*$ to \eqnref{eqn_ext_diri} with $\Bg\in H^{1/2}(\p D^e)^2$. Then the following variational principle holds:
\beq\label{variation}
\Ecal_{D^e}[\Bu] =\min_{\Bw\in W_\Bg} \Ecal_{D^e}[\Bw],
\eeq
 where
$$
W_{\Bg} = \big\{ \Bw\in \Acal^* : \Bw|_{\p D^e} = \Bg\big\}.
$$
 \end{lemma}

\pf
Let $\Bw \in W_{\Bg}$. By Lemma \ref{cor_betti}, we have
\begin{align*}
\int_{D^e}\Cbb \hatna\Bu:\hatna\Bu  = -\int_{\p D^e} \p_\Gv \Bu |_+ \cdot \Bg = \int_{D^e}\Cbb \hatna\Bu:\hatna\Bw.
\end{align*}
By the Cauchy-Schwartz inequality, we have
\begin{align*}
\int_{D^e} \Cbb\hatna\Bu:\hatna\Bu
= \int_{D^e} \Cbb\hatna\Bu:\hatna\Bw
\le \frac{1}{2}(\int_{D^e} \Cbb\hatna\Bu:\hatna\Bu
+ \int_{D^e} \Cbb\hatna\Bw:\hatna\Bw).
\end{align*}
Thus \eqnref{variation} holds. \qed

%%%%%%%%%%%%%%%%%%%%%%%%%%%%%%%%%%%%%%%%%%%%%%%%%%%%%%%%%%%%%%%%%%%%%%%%%%%%%%%%%%%
\subsection{An estimate for the free space problem}\label{subsec:nabla_k_u_H_bdd}
%%%%%%%%%%%%%%%%%%%%%%%%%%%%%%%%%%%%%%%%%%%%%%%%%%%%%%%%%%%%%%%%%%%%%%%%%%%%%%%%%%%

The purpose of this subsection is to prove the following proposition which will be used in section \ref{sec:free}.

\begin{prop}\label{freeest}
Let $\Bu$ be the solution to \eqnref{elas_eqn_free} for a given $\BH$. Then for any disk $B$ centered at $0$ containing $\overline{D_1 \cup D_2}$ and for $k=0,1,2, \ldots$, there is a constant $C_k$ independent of $\Ge$ (and $\BH$) such that
\beq\label{eqn_nablak_u_H_estim}
\| \nabla^k (\Bu-\BH)\|_{L^\infty(\Rbb^2\setminus B)}\le C_k \|\BH \|_{H^{1}(B)}.
\eeq
\end{prop}

The main emphasis of \eqnref{eqn_nablak_u_H_estim} is that the estimate holds independently of $\Ge$, the distance between $D_1$ and $D_2$. It shows that even if $\Bu$ depends on $\Ge$, the dependence is negligible far away from the inclusions.

To prove Proposition \ref{freeest}, we begin with the following lemma.
\begin{lemma}\label{lem:Ecal_u_H_estim}
Let $\Bu$ be the solution to \eqnref{elas_eqn_free}. There is a constant $C$ independent of $\Ge$ and $\BH$ such that
\beq
\Ecal_{D^e}[\Bu-\BH] \le C \|\BH\|^2_{H^{1}(B)},
\eeq
where $\Ecal_{D^e}$ is defined in \eqnref{eqn_def_Ecal} and $B$ is a disk centered at $0$ containing $\ol{D_1 \cup D_2}$.
\end{lemma}
\pf
We first observe from \eqnref{int_zero} and the second condition in \eqnref{elas_eqn_free} that
$$
\int_{\p D_i} \p_\Gv \Bu |_+ \cdot \Bu=0, \quad i=1,2.
$$
Since $\Lcal_{\Gl,\Gm} \BH=0$ in $\Rbb^2$, we also have
$$
\int_{\p D_i} \p_\Gv \BH |_+ \cdot \Bu=0, \quad i=1,2.
$$
So we have
\begin{align*}
\Ecal_{D^e}[\Bu-\BH] &= \int_{D^e}\Cbb \hatna(\Bu-\BH):\hatna(\Bu-\BH) \\
 &= -\int_{\p D^e} \p_\Gv (\Bu-\BH) |_+ \cdot (\Bu-\BH)
=\int_{\p D^e} \p_\Gv (\Bu-\BH) |_+ \cdot \BH.
\end{align*}

Let $R$ be the radius of $B$ and let $r$ be such that $r < R$ and $\overline{D_1 \cup D_2} \subset B_{r}$. Let $\chi$ be a smooth radial function such that $\chi(\Bx)=1$ if $|\Bx| \le r$ and $\chi(\Bx)=0$ if $|\Bx| \ge R$. Let $\Bw:= -\chi \BH$. Then we have
\begin{align*}
\int_{D^e}\Cbb \hatna(\Bu-\BH):\hatna\Bw
= -\int_{\p D^e} \p_\Gv (\Bu-\BH) |_+ \cdot \Bw
=\int_{\p D^e} \p_\Gv (\Bu-\BH) |_+ \cdot \BH.
\end{align*}
It then follows that
\begin{align*}
\Ecal_{D^e}[\Bu-\BH] = \int_{D^e} \Cbb\hatna(\Bu-\BH):\hatna\Bw \leq \frac{1}{2} \left( \Ecal_{D^e}[\Bu-\BH] + \Ecal_{D^e}[\Bw] \right).
\end{align*}
So we have
\begin{align*}
\Ecal_{D^e}[\Bu-\BH] &\leq
\Ecal_{D^e}[\Bw] \le C \| \BH\|_{H^1(B)}^2.
\end{align*}
The proof is completed.
\qed

\medskip
\noindent{\sl Proof of Proposition \ref{freeest}}.
By Proposition \ref{thm_u_layer}, the solution $\Bu$ is represented as
\beq
\Bu = \BH + \Scal_{\p D_1}[\BGvf_1]+\Scal_{\p D_2}[\BGvf_2]
\eeq
with $\BGvf_j = \p_\Gv \Bu|_+$ on $\p D_j$, $j=1,2$. Proposition \ref{thm:funda_sol_series} yields
\beq\label{u_H_far_expand}
(\mathbf{\Bu-\BH})(\Bx)
= \sum_{n\neq 0}\frac{1}{2|n|}\frac{e^{-in\Gt}}{r^{|n|}} \left( -M_n^{(1)} \Be_1 - M_n^{(2)} \Be_2 +   M_n^{(3)}\Bx \right), \quad r=|\Bx| \geq R,
\eeq
where $R$ is the radius of $B$ and $M_n^{(i)}$ is given by
\begin{align}
M_n^{(i)} &= \int_{\p D_1} \Bv_n^{(i)} \cdot \BGvf_1 \, d\sigma +
\int_{\p D_2} \Bv_n^{(i)} \cdot \BGvf_2 \, d\sigma, \quad i=1,2,
\\
M_n^{(3)} &= \int_{\p D_1} \Bw_n \cdot \BGvf_1 \, d\sigma +
\int_{\p D_2} \Bw_n \cdot \BGvf_2 \, d\sigma.
\end{align}
Observe that the dependence of $\Bu-\BH$ on $\Ge$ is contained only in the coefficients $M_n^{(i)}$, $i=1,2,3$.

Let $r$ be such that $r < R$ and $B_{r}$ contains $\overline{D_1 \cup D_2}$. We now show that
there is a constant $C$ independent of $\Ge$ and $n$ such that
\beq\label{M3_estim}
|M_n^{(i)}| \le C |n| r^{|n|+2} \|\BH \|_{H^1(B_r)}
\eeq
for all $n \neq 0$ and $i=1,2,3$.

For simplicity, we consider only $i=1$. The other cases can be proved in the exactly same way.
Let $\Bv$ be the solution to \eqnref{elas_eqn_free} when $\BH=\Bv^{(1)}_n$. Since $\BGvf_j = \p_\Gv \Bu|_+$ on $\p D_j$ and $\Bv|_{\p D_j} \in \Psi$, we have using \eqnref{int_zero} that
\begin{align*}
M_n^{(1)} &= \int_{\p D^e} \Bv_n^{(1)} \cdot \p_\Gv \Bu|_+  \\
&= -\int_{\p D^e} (\Bv-\Bv_n^{(1)}) \cdot \p_\Gv \Bu|_+  \\
&= -\int_{\p D^e} (\Bv-\Bv_n^{(1)}) \cdot \p_\Gv (\Bu-\BH)|_+
+\int_{\p D^e} \Bv_n^{(1)} \cdot \p_\Gv \BH|_+ \\
&= \int_{D^e} \Cbb \hatna (\Bv-\Bv_n^{(1)}): \hatna (\Bu-\BH)
+ \int_{\p D^e} \Bv_n^{(1)} \cdot \p_\Gv \BH|_+ .
\end{align*}
So, by applying Lemma \ref{lem:Ecal_u_H_estim} to $\Ecal_{D^e}[\Bv-\Bv_n^{(1)}]$ and $\Ecal_{D^e}[\Bu-\BH]$ on $B_{r}$, we obtain
\begin{align*}
|M_n^{(1)}| &\le \Ecal_{D^e}\big[\Bv-\Bv_n^{(1)}\big]^{1/2}
\Ecal_{D^e}\big[\Bu-\BH\big]^{1/2}
+ \sum_{i=1}^2 \|\Bv_n^{(1)} \|_{H^{1/2}(\p D_i)} \| \p_\Gv \BH \|_{H^{-1/2}(\p D_i)} \\
& \le C \|\Bv_n^{(1)} \|_{H^1(B_{r})} \|\BH \|_{H^1(B_{r})} + \sum_{i=1}^2 \|\Bv_n^{(1)} \|_{H^{1}(D_i)} \| \BH \|_{H^{1}(D_i)} \\
& \le C \|\Bv_n^{(1)} \|_{H^1(B_{r})} \|\BH \|_{H^1(B_{r})}.
\end{align*}
Here and throughout this paper, the constant $C$ appearing in the course of estimations may differ at each occurrence. Since $\Bv_n^{(1)}$ is a homogeneous polynomial of order $n$, there is a constant $C$ independent of $n$ such that
$$
\|\Bv_n^{(1)} \|_{H^1(B_{r})} \le C |n| r^{|n|+2},
$$
assuming that $r >1$.

It follows from \eqnref{u_H_far_expand} and \eqnref{M3_estim} that
\begin{align*}
\|\Bu-\BH\|_{L^\infty(\Rbb^2\setminus B_R)} & \le C \sum_{n\neq 0} \frac{1}{2|n|} \frac{1}{R^{|n|}}(|n| r^{|n|+2})\|\BH \|_{H^1(B_r)} \\
& \le C \sum_{n\neq 0} \left( \frac{r}{R} \right)^{|n|} \|\BH \|_{H^1(B)} \le C \|\BH \|_{H^1(B)}.
\end{align*}
This proves \eqnref{eqn_nablak_u_H_estim} for $k=0$.

If $k>0$, we differentiate \eqnref{u_H_far_expand} to obtain \eqnref{eqn_nablak_u_H_estim}. This completes the proof.
\qed

%%%%%%%%%%%%%%%%%%%%%%%%%%%%%%%%%%%%%%%%%%%%%%%%%%%%%%%%%%%%%%%%%%%%%%%%%%%
\subsection{Geometry of two inclusions}\label{subsection:geo_two_incl}
\begin{figure*}%[h!]
\begin{center}
\epsfig{figure=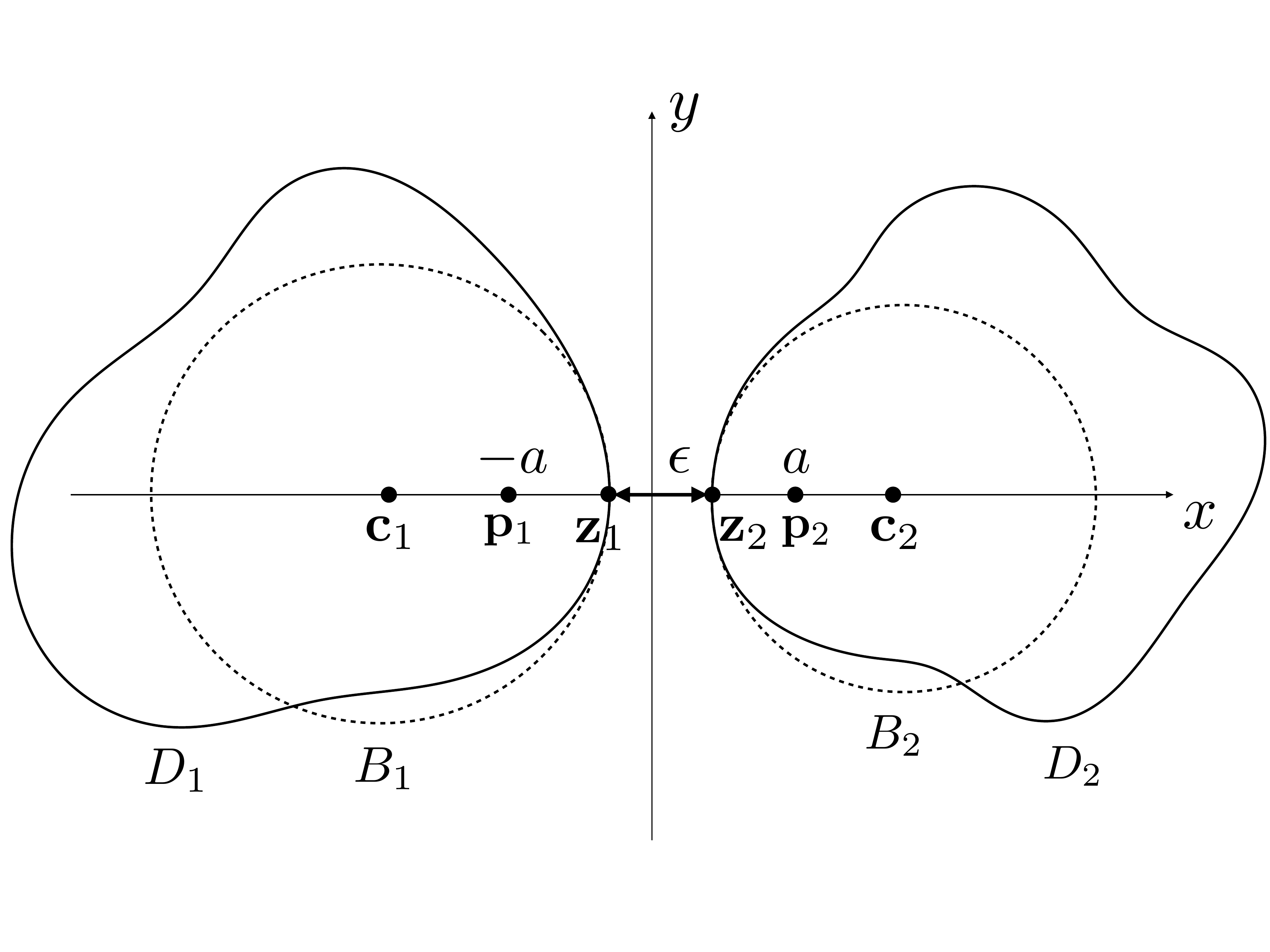,width=8cm}
\end{center}
\caption{Geometry of the two inclusions and osculating circles}
\label{fig:two_inclusions}
\end{figure*}

%%%%%%%%%%%%%%%%%%%%%%%%%%%%%%%%%%%%%%%%%%%%%%%%%%%%%%%%%%%%%%%%%%%%%%%%%%%
In this subsection, we describe geometry of two inclusions $D_1$ and $D_2$. See Figure \ref{fig:two_inclusions}.

Suppose that there are unique points $\Bz_1\in \p D_1$ and $\Bz_2\in \p D_2$ such that
\beq
|\Bz_1-\Bz_2| = \mbox{dist}(D_1, D_2).
\eeq
We assume that $D_j$ is strictly convex near $\Bz_j$, namely, there is a common neighborhood $U$ of $\Bz_1$ and $\Bz_2$ such that $D_j \cap U$ is strictly convex for $j=1,2$. Moreover, we assume that
$$
\mbox{dist}(D_1, D_2 \setminus U) \ge C \quad\mbox{and}\quad \mbox{dist}(D_2, D_1 \setminus U) \ge C
$$
for some positive constant $C$ independent of $\Ge$. This assumption says that other than neighborhoods of $\Bz_1$ and $\Bz_2$, $D_1$ and $D_2$ are at some distance to each other. We need one more assumption: the center of the circle which is osculating $D_j$ at $\Bz_j$ lies inside $D_j$ ($j=1,2$). This assumption is needed for defining the singular function $\Bq_3$ in \eqnref{q3def} later. We emphasize that strictly convex domains satisfy all the assumptions.

Let $\Gk_j$ be the curvature of $\p D_j$ at $\Bz_j$. Let $B_j$ be the disk osculating to $D_j$ at $\Bz_j$ ($j=1,2$). Then the radius $r_j$ of $B_j$ is given by  $r_j=1/\Gk_j$. Let $R_j$ be the reflection with respect to $\p B_j$ and let $\Bp_1$ and $\Bp_2$ be the unique fixed points of
the combined reflections $R_1\circ R_2$ and $R_2\circ R_1$, respectively.

Let $\Bn$ be the unit vector in the direction of $\Bp_2-\Bp_1$ and let $\Bt$ be the unit vector perpendicular to $\Bn$ such that $(\Bn,\Bt)$ is positively oriented.
We set $(x,y)\in\Rbb^2$ to be the Cartesian coordinates such that
$\Bp=(\Bp_1+\Bp_2)/2$ is the origin and the $x$-axis is parallel to $\Bn$.
Then one can see (see \cite{AKL-MA-05}) that $\Bp_1$ and $\Bp_2$ are written as
\beq\label{pjdef}
\Bp_1=(-a,0)\quad\mbox{and}\quad\Bp_2=(a,0),
\eeq
where the constant $a$ is given by
\beq\label{a_def}
a :=\frac{\sqrt\Ge\sqrt{ (2 r_1 + \Ge) (2 r_2 + \Ge) (2 r_1 + 2 r_2 +
  \Ge)}}{2 (r_1 + r_2 + \Ge)},
\eeq
from which one can infer
\beq\label{def_alpha}
a= \sqrt{\frac{2 }{\Gk_1+\Gk_2}}\sqrt{\Ge}+O(\Ge^{3/2}).
\eeq
Then the center $\Bc_i$ of $B_i$ ($i=1,2$) is given by
\beq\label{c1c2}
\Bc_i = \Big((-1)^i \sqrt{r_i^2+a^2},0\Big)=\big((-1)^i r_i +O(\Ge),0\big).
\eeq
So we have
\beq
\Bz_i = (-1)^{i+1}\Big(r_i-\sqrt{r_i^2+a^2},0\Big) = \left( (-1)^i\frac{ \Gk_i}{\Gk_1+\Gk_2} \Ge + O(\Ge^2), 0 \right).
\eeq

\begin{figure*}%[h!]
\begin{center}
\epsfig{figure=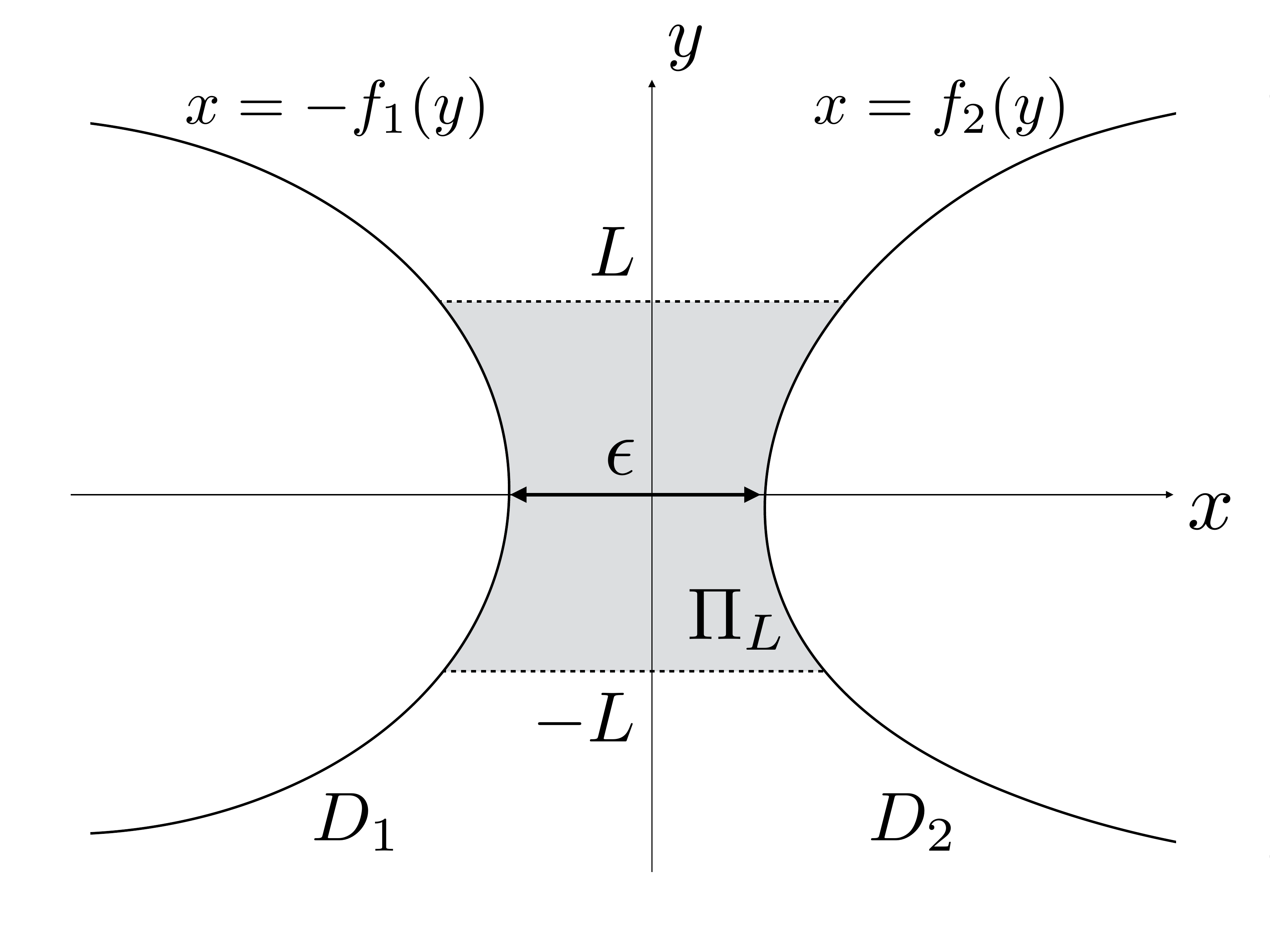,width=6cm}
\end{center}
\caption{Geometry of the narrow gap region $\Pi_L$}
\label{fig:narrow_gap}
\end{figure*}

Let us consider the narrow region between $D_1$ and $D_2$. See Figure \ref{fig:narrow_gap}.
There exists $L>0$ (independent of $\Ge$) and functions $f_1,f_2:[-L,L]\rightarrow \Rbb$ such that
\beq\label{foneftwo}
\Bz_1=\big(-f_1(0),0\big), \quad \Bz_2=\big(f_2(0),0\big), \quad f_1'(0)= f_2'(0)=0,
\eeq
and $\p D_1$ and $\p D_2$ are graphs of $-f_1(y)$ and $f_2(y)$ for $|y|<L$, {\it i.e.},
\beq\label{Bx1Bx2}
\Bx_1(y):=(-f_1(y),y)\in \p D_1 \quad\mbox{and}\quad \Bx_2(y):=(f_2(y),y)\in \p D_2.
\eeq
Since $D_i$ is strictly convex near $\Bz_j$, $f_1$ is strictly convex.
Note that, for $i=1,2$ and $|y|<L$,
\beq\label{fitaylor}
f_i(y) = \frac{\Gk_i}{\Gk_1+\Gk_2}\Ge+\frac{1}{2!}\Gk_i y^2 +\frac{1}{3!}\Go_i y^3 + O(\Ge^2+y^4)
\eeq
for some constant $\Go_i$. Let us define for later use a constant $\tau$ as
\beq\label{def_tau}
\tau=|\Gk_1-\Gk_2|+|\Go_1|+|\Go_2|.
\eeq
We denote by $\Pi_l$ for $0<l\le L$  the narrow region between $D_1$ and $D_2$ defined as
\beq\label{narrowregion}
\Pi_l = \{ (x,y) \in\Rbb^2| -f_1(y)<x<f_2(y), \, |y|<l\}.
\eeq

%%%%%%%%%%%%%%%%%%%%%%%%%%%%%%%%%%%%%%%%%%%%%%%%%%%%%%%%%%%%%%%%%%%%%%%
\section{Singular functions and their properties}\label{sec:singular}
%%%%%%%%%%%%%%%%%%%%%%%%%%%%%%%%%%%%%%%%%%%%%%%%%%%%%%%%%%%%%%%%%%%%%%%

Let $\Psi_j$, $j=1,2,3$, be the rigid motions defined in \eqnref{Psidef} and let $\Bh_j$ be the solution to the following problem:
\beq\label{hj_def}
 \ \left \{
 \begin{array} {ll}
\ds \Lcal_{\Gl,\Gm} \Bh_j= 0 \quad &\mbox{ in } D^e,\\[2mm]
\ds \Bh_j= -\frac{1}{2}\Psi_j(\Bx) \quad &\mbox{ on } \p D_1,\\[2mm]
\ds \Bh_j= \frac{1}{2}\Psi_j(\Bx)\quad &\mbox{ on } \p D_2.
 \end{array}
 \right.
\eeq

It turns out that $\Bh_j$ ($j=1,2,3$) captures the singular behavior of the solution $\Bu$ to \eqnref{elas_eqn_free}. In fact, $\Bu$ can be decomposed in the following form:
\beq
\Bu = \sum_{j=1}^3 c_j \Bh_j + \Bb
\eeq
for some constants $c_j$, where $\nabla \Bb$ is bounded in a bounded domain containing the narrow region $\Pi_L$ between $D_1$ and $D_2$. In other words, the blow-up behavior of $\nabla \Bu$ is completely characterized by that of $\sum_{j=1}^3 c_j \nabla \Bh_j$.
We emphasize that $|\Bh_j |_{\p D_1} - \Bh_j |_{\p D_2}|=1$ for $j=1,2$. So one expect that $|\nabla \Bh_j| \approx \Ge^{-1}$ in the narrow region between $D_1$ and $D_2$. The function $\Bh_3$ has a weaker singularity since $|\Bh_3 |_{\p D_1} - \Bh_3 |_{\p D_2}|=|\Bx|$.

The purpose of this section is to construct explicit singular functions, denoted by $\Bq_j$, which yield good approximations of $\Bh_j$ and to derive their important properties.

%%%%%%%%%%%%%%%%%%%%%%%%%%%%%%%%%%
\subsection{Construction of singular functions}\label{subsec:singular}
%%%%%%%%%%%%%%%%%%%%%%%%%%%%%%%%%%

We begin with a brief review of the singular function for the electro-static case.
Let $\Bp_j$ ($j=1,2$) be the fixed points of the combined reflections given in \eqnref{pjdef} and let
\beq\label{qBdef}
q_B(\Bx) = \frac{1}{2\pi} (\ln|\Bx-\Bp_1|-\ln|\Bx-\Bp_2|).
\eeq
This function was introduced in \cite{Yun-SIAP-07} and used in an essential way for characterization of the gradient blow-up in the context of electro-statics  \cite{ACKLY-ARMA-13}. The most important property of $q_B(\Bx)$ is that it takes constant values on $\p B_j$, the circles osculating to $\p D_j$ at $\Bz_j$, $j=1,2$. It is because $\p B_1$ and $\p B_2$ are circles of Apollonius of $\Bp_1$ and $\Bp_2$.

Note that $\frac{1}{2\pi}\ln|\Bx|$ is a fundamental solution of the Laplacian and represents a point source of the electric field. So it is natural to expect that, even in the linear elasticity case, the point source functions may also characterize the gradient blow-up. There are various types of point source functions in linear elasticity which are often called {\it nuclei of strain}. We will use the following nuclei of strain as basic building blocks of the singular functions:
\beq\label{nuclei}
\BGG(\Bx)\Be_1, \quad \BGG(\Bx)\Be_2, \quad
\frac{\Bx}{|\Bx|^2}, \quad \frac{\Bx^\perp}{|\Bx|^2}.
\eeq
where $\Bx^\perp=(-y,x)$ for $\Bx=(x,y)\in \Rbb^2$.
These nuclei of strain have physical meanings: the function $\BGG(\Bx)\Be_j$ represents the point force applied at the origin in the direction of $\Be_j$, and the functions ${\Bx}/{|\Bx|}$ and ${\Bx^\perp}/{|\Bx|}$ represent the point source of the pressure and that of the moment located at the origin, respectively (see, for example, \cite{Love}).
%The last two functions in \eqnref{nuclei} are derivatives of nuclei of strain.

We emphasize that the functions given in \eqnref{nuclei} are solutions to the Lam\'{e} system for $\Bx \neq 0$.
In fact, the first two are solutions since they are columns of the fundamental solution, and so are the last two because of the following relations:
\beq\label{doublet}
\begin{cases}
\ds (\Ga_1-\Ga_2) \frac{\Bx}{|\Bx|^2} = \p_1 (\BGG(\Bx) \Be_1) + \p_2 (\mathbf{\GG(\Bx)} \Be_2),
\\[1em]
\ds (\Ga_1+\Ga_2) \frac{\Bx^\perp}{|\Bx|^2} = \p_1 (\BGG(\Bx) \Be_2) - \p_2 (\mathbf{\GG(\Bx)} \Be_1),
\end{cases}
\eeq
where $\Ga_1$ and $\Ga_2$ are constants appearing in the definition \eqnref{Kelvin} of the fundamental solution. The identities in \eqnref{doublet} can be proved by straightforward computations.

The singular functions of this paper are constructed as linear combinations of functions given in \eqnref{nuclei}. To motivate the construction, we temporarily assume that two inclusions $D_1$ and $D_2$ are symmetric with respect to both $x$- and $y$-axes. If we write $\Bh_1=(h_{11}, h_{12})^T$, then thanks to the symmetry of the inclusions and boundary conditions in \eqnref{hj_def}, the following two functions are also solutions of \eqnref{hj_def} for $j=1$:
$$
\begin{bmatrix} h_{11}(x,-y) \\ -h_{12}(x,-y) \end{bmatrix}, \quad
\begin{bmatrix} -h_{11}(-x,y) \\ h_{12}(-x,y) \end{bmatrix}.
$$
By the uniqueness of the solution, we see that $\Bh_1$ has the following symmetric property with respect to $x$- and $y$-axes:
\beq\label{honesymm}
\begin{cases}
h_{11}(x,y)=h_{11}(x,-y)=-h_{11}(-x,y), \\
h_{12}(x,y)=-h_{12}(x,-y)=h_{12}(-x,y).
\end{cases}
\eeq
One can see that the following two functions have the same symmetry:
$$
\BGG(\Bx-\Bp_1)\Be_1-\BGG(\Bx-\Bp_2)\Be_1, \quad \frac{\Bx-\Bp_1}{|\Bx-\Bp_1|^2}+\frac{\Bx-\Bp_2}{|\Bx-\Bp_2|^2}.
$$
So the first singular function $\Bq_1$ is constructed as a linear combination of these functions.

On the other hand, one can see in a similar way that $\Bh_2=(h_{21}, h_{22})^T$ has the following symmetric property:
\beq\label{htwosymm}
\begin{cases}
h_{21}(x,y)=-h_{21}(x,-y)=h_{21}(-x,y), \\
h_{22}(x,y)=h_{22}(x,-y)=-h_{22}(-x,y),
\end{cases}
\eeq
and the following two functions have the same symmetry:
$$
\BGG(\Bx-\Bp_1)\Be_2-\BGG(\Bx-\Bp_2)\Be_2, \quad \frac{(\Bx-\Bp_1)^{\perp}}{|\Bx-\Bp_1|^2}+\frac{(\Bx-\Bp_2)^{\perp}}{|\Bx-\Bp_2|^2}.
$$
So $\Bq_2$ is constructed as a linear combination of these functions.

The singular functions of this paper are defined by
\beq\label{Bqone}
\Bq_1(\Bx) := \BGG(\Bx-\Bp_1)\Be_1-\BGG(\Bx-\Bp_2) \Be_1 + {\Ga_2 a} \left( \frac{\Bx-\Bp_1}{|\Bx-\Bp_1|^2}+\frac{\Bx-\Bp_2}{|\Bx-\Bp_2|^2}\right),
\eeq
and
\beq\label{Bqtwo}
\Bq_2(\Bx) :=  \BGG(\Bx-\Bp_1)\Be_2 -\BGG(\Bx-\Bp_2)  \Be_2 - {\Ga_2 a} \left( \frac{(\Bx-\Bp_1)^{\perp}}{|\Bx-\Bp_1|^2}+\frac{(\Bx-\Bp_2)^{\perp}}{|\Bx-\Bp_2|^2}\right),
\eeq
where $a$ is the number appearing in \eqnref{pjdef}. We emphasize that $a$ depends on $\Ge$. In fact, we repeatedly use the fact that $a \approx \sqrt{\Ge}$.
The functions $\Bq_j$ satisfy $\Lcal_{\Gl, \Gm} \Bq_j=0$ in $\Rbb^2 \setminus \{ \Bp_1, \Bp_2 \}$, and
\beq
\Bq_j(\Bx)= O(|\Bx|^{-1}) \quad\mbox{as } |\Bx| \to \infty,
\eeq
as one can easily see. We emphasize that the symmetry of $D_1\cup D_2$ is not assumed here.

It will be proved later in Proposition \ref{prop_Brj} that
\beq\label{BhjBqj}
\Bh_j \approx \frac{m_j}{\sqrt\Ge} \Bq_j
\eeq
where $m_1$ and $m_2$ are constants defined by
\beq\label{m2def}
m_1 := \big[(\Ga_1 - \Ga_2)\sqrt{2(\Gk_1+\Gk_2)}\,\big]^{-1}, \quad
m_2 := \big[(\Ga_1 + \Ga_2)\sqrt{2(\Gk_1+\Gk_2)}\,\big]^{-1}.
\eeq
So blow-up of $\nabla \Bh_j$ is captured by an explicit function $\frac{m_j}{\sqrt\Ge} \nabla \Bq_j$. This is a crucial fact for investigating blow-up of $\nabla \Bu$ in this paper.

We now construct the third singular function $\Bq_3$ which approximates $\Bh_3$.
For that we introduce
$\BGG^\perp = \big( \GG_{ij}^\perp \big)_{i, j = 1}^2$ which is defined by
\beq\label{Kelvin_modified}
  \BGG^\perp(\Bx) =
    \ds \Ga_1  \arg{(\Bx)} \begin{bmatrix}1 & 0 \\ 0 & 1\end{bmatrix}
     - \frac{\Ga_2}{|\Bx| ^2}  \begin{bmatrix} -x_1 x_2 & -x_2^2 \\ x_1^2 & x_1 x_2 \end{bmatrix}.
\eeq
We emphasize that $\BGG^\perp$ is a multi-valued function since $\arg{(\Bx)}$ is. So $\BGG^\perp$ is defined in $\Rbb^2$ except a branch-cut starting from the origin.

Note that
\beq\label{Gej}
\BGG^\perp(\Bx)\Be_j =\Ga_1 \arg(\Bx) \Be_j - \Ga_2 x_j \nabla (\arg(\Bx)), \quad j=1,2.
\eeq
Since $\arg(\Bx)$ is a harmonic function, we infer from Lemma \ref{harLame} that $\BGG^\perp(\Bx)\Be_j$ is a solution to the Lam\'{e} system (except on the branch-cut).

We now define the singular function $\Bq_3$ by
\begin{align}
\Bq_3(\Bx) &= m_3\left( \BGG^\perp(\Bx-\Bp_1) - \BGG^\perp(\Bx-\Bc_1) \right) \Be_1 + m_3\left( \BGG^\perp(\Bx-\Bp_2) - \BGG^\perp(\Bx-\Bc_2) \right) \Be_1 \nonumber \\
&\quad +m_3\Ga_2 a \left( \frac{(\Bx-\Bp_1)^\perp}{|\Bx-\Bp_1|^2} - \frac{(\Bx-\Bp_2)^\perp}{|\Bx-\Bp_2|^2} \right), \label{q3def}
\end{align}
where
\beq\label{m3def}
m_3:= \big[{(\Ga_1-\Ga_2)(\Gk_1+\Gk_2)}\big]^{-1},
\eeq
and $\Bc_1$ and $\Bc_2$ are centers of the osculating disks $B_1$ and $B_2$, respectively. It is worth mentioning that $\BGG^\perp(\Bx-\Bp_j) - \BGG^\perp(\Bx-\Bc_j)$ is well-defined in $\Rbb^2$ except a branch-cut connecting $\Bp_j$ and $\Bc_j$. So, $\Bq_3$ is well-defined and a solution of the Lam\'e system in $D^e$.
We will show in Lemma \ref{lem_q3_asymp} that $\Bq_3$ has the same local behavior as $\Bh_3$.

In subsections to follow we derive technical estimates of $\Bq_j$ and its derivatives which will be used in later sections.

%%%%%%%%%%%%%%%%%%%%%%%%%%%%%%%%%%%%%%%%%%%%%%%%%%%%%%%%%%%%%%%%%%%%%%%
\subsection{Estimates of the function $\Gz$}\label{sec:zeta}
%%%%%%%%%%%%%%%%%%%%%%%%%%%%%%%%%%%%%%%%%%%%%%%%%%%%%%%%%%%%%%%%%%%%%%%

We show in the next subsection that the singular function $\Bq_j$ ($j=1,2$) can be nicely represented using the function $q_B$ given in \eqnref{qBdef}. In fact, it is slightly more convenient to use the function $\Gz(\Bx)$ defined by
\beq
\Gz(\Bx) = 2\pi q_B(\Bx).
\eeq
The following lemma collects estimates for the function $\Gz$ to be used in the next subsection. Some of the estimates are essentially proved in \cite{ACKLY-ARMA-13}. However, in that paper the estimates are not explicitly written and derivations of estimates are smeared in other proofs. So we include proofs.

\begin{lemma}
\begin{itemize}
\item[(i)] Let $\Pi_{L}$ be the narrow region defined in \eqnref{narrowregion}. It holds that
\beq\label{GzBxest}
|\Gz(\Bx)| \lesssim \sqrt{\Ge}, \quad \Bx\in \Pi_{L},
\eeq
and
\beq\label{Gzgradest}
|\p_1 \Gz(\Bx) | \lesssim \frac{\sqrt\Ge}{\Ge+y^2}, \quad
|\p_2 \Gz(\Bx) | \lesssim \frac{\sqrt{\Ge}|y|}{\Ge+y^2}, \quad \Bx=(x,y) \in \Pi_{L}.
\eeq

\item[(ii)] Let $\Bx_j(y)$ be the defining functions for $\p D_j$ for $j=1,2$ as given in \eqnref{Bx1Bx2}. For $|y|<L$ and $j=1,2$, we have
\begin{align}
|\Gz(\Bx_j(y))-\Gz|_{\p B_j}| &\lesssim \sqrt{\Ge} (|\Go_j||y|+|y|^2), \label{difference} \\
\Big| \frac{d}{dy} \Gz(\Bx_j(y))\Big| &\lesssim \sqrt{\Ge}, \label{eqn_bdry_deri1_1i} \\
\Big|\frac{d^2}{dy^2}\Gz(\Bx_j(y))\Big| &\lesssim \frac{\sqrt\Ge}{{\Ge}+ y^2}. \label{eqn_bdry_deri2_1i}
\end{align}
\end{itemize}
\end{lemma}

\pf
(i) Since $\Bp_1=(-a,0)$ and $\Bp_2=(a,0)$, we can rewrite $\Gz(\Bx)$ as
\beq\label{xi_cartesian}
\Gz(\Bx)=\frac{1}{2}\ln \frac{(x+a)^2+y^2}{(x-a)^2+y^2} = \frac{1}{2} \ln \left( 1+ \frac{4 a x}{(x-a)^2+y^2} \right).
\eeq
Since $a \approx \sqrt\Ge$, $\Ge+y^2\approx (x \pm a)^2+y^2$, and $|x|\lesssim \Ge+y^2$ for $(x,y)\in\Pi_{L}$,
we obtain
$$
|\Gz(\Bx)| =\frac{1}{2}\ln(1+O(\sqrt\Ge)) =  O( \sqrt{\Ge}) \quad\mbox{for } \Bx\in\Pi_{L},
$$
which yields \eqnref{GzBxest}.

Assume $\Bx=(x,y)\in \Pi_L$.
Since
$$
|x|\lesssim \Ge+y^2 \quad\mbox{and}\quad \Ge+y^2\lesssim (x\pm a)^2+y^2,
$$
we have from the first identity in \eqnref{xi_cartesian} that
\begin{align*}
|\p_1 \Gz(\Bx)| &= \Big| \frac{x+a}{(x+a)^2+y^2} - \frac{x-a}{(x-a)^2+y^2}\Big|
\\
&=\Big|\frac{2 a (a^2-x^2+y^2) }{((x+a)^2+y^2) ((x+a)^2+y^2)} \Big|
\lesssim \frac{\sqrt\Ge(\Ge+y^2)}{(\Ge+y^2)^2} \lesssim \frac{\sqrt\Ge }{\Ge+y^2},
\end{align*}
and
\begin{align*}
|\p_2 \Gz(\Bx)|
&= |y| \Big| \frac{1}{(x+a)^2+y^2} - \frac{1}{(x-a)^2+y^2}\Big|
\\
&\lesssim |y| \Big|\frac{4 a x }{((x+a)^2+y^2) ((x-a)^2+y^2)} \Big|
\lesssim \frac{\sqrt{\Ge} |y| |x|}{(\Ge+y^2)^2} \lesssim\frac{\sqrt{\Ge}|y|}{\Ge+y^2}.
\end{align*}
So \eqnref{Gzgradest} is proved.

(ii) We now prove \eqnref{difference}. For simplicity, we assume $j=1$.
Let us write the boundary $\p B_1$ of the osculating disk $B_1$ as $(-f_B(y),y)$ for $|y|<L$. Recall from \eqnref{fitaylor} that
\beq\label{f1fB_diff}
f_1(y) - f_B(y) = \frac{1}{3!} \omega_1 y^3 + O(y^4).
\eeq
From \eqnref{xi_cartesian}, we have
\begin{align}
|\Gz(\Bx_1(y))-\Gz|_{\p B_1} |&=\frac{1}{2}
\Big|
\ln \Big(1-\frac{4af_1(y)}{(f_1(y)+a)^2+y^2}\Big)
-
\ln \Big(1-\frac{4af_B(y)}{(f_B(y)+a)^2+y^2}\Big)
\Big|
\nonumber
\\
&=
\frac{1}{2}
\Big|
\ln \Big(1- 4a  \frac{\eta_1(y)}{\eta_2(y)} \Big)\Big|,
\label{xi1y_xiB_diff}
\end{align}
where
\begin{align}
\eta_1(y) &= \frac{f_1(y)}{(f_1(y)+a)^2+y^2} - \frac{f_B(y)}{(f_B(y)+a)^2+y^2},
\\
\eta_2(y) &=1-\frac{4af_B(y)}{(f_B(y)+a)^2+y^2}.
\end{align}
Since $a\approx \sqrt{\Ge}$, $f_B(y) \approx \Ge+y^2$, and $(f_B(y)+a)^2+y^2 \approx \Ge+y^2$,
we see that
\beq\label{G2_estim}
|\eta_2(y)|\approx 1 \quad \mbox{for }|y|<L.
\eeq
From \eqnref{f1fB_diff} and the facts that $f_1 \approx \Ge+y^2$, $f_B \approx \Ge+y^2$ and $a\approx \sqrt\Ge$, we have, for $|y|<L$,
\begin{align}
|\eta_1(y)| &= \left| \frac{(f_1(y)-f_B(y))(y^2-f_1(y)f_B(y)+a^2)}{((f_1(y)+a)^2+y^2)((f_B(y)+a)^2+y^2)} \right|
\nonumber
\\
&\lesssim \frac{(\omega_1 |y|^3 +y^4) (y^2+(\Ge+y^2)^2+\Ge)}{(\Ge+y^2)^2}
\nonumber
\\
&\lesssim |\omega_1| |y| + y^2.
\label{G1_estim}
\end{align}
Since $a\approx \sqrt\Ge$, it follows from \eqref{xi1y_xiB_diff}, \eqnref{G2_estim}, and \eqnref{G1_estim} that
$$
|\Gz(\Bx_1(y))-\Gz|_{\p B_1} | \lesssim \sqrt\Ge (|\omega_1||y|+y^2).
$$
Therefore \eqnref{difference} is proved.

We now prove \eqnref{eqn_bdry_deri1_1i} and \eqnref{eqn_bdry_deri2_1i} for $j=1$. The cases for $j=2$ can be handled similarly.
In view of \eqnref{xi_cartesian}, we have
\begin{align}
\frac{d}{dy} \Gz(\Bx_1(y))
&=\frac{d}{dy} \Gz(-f_1(y),y)\nonumber\\
&=\Big(\frac{-f_1(y)+ a}{(-f_1(y)+ a)^2+y^2}
-\frac{-f_1(y)- a}{(-f_1(y)- a)^2+y^2}\Big)(-f_1'(y))
\nonumber\\
&\quad\quad + \frac{y}{(-f_1(y)+ a)^2+y^2}
-\frac{y}{(-f_1(y)- a)^2+y^2}
\nonumber\\
&=2 a \frac{N(y)}{D_+(y) D_-(y)},
\label{xi1_bdry_deri_NDD}
\end{align}
where $N$ and $D_\pm$ are given by
\begin{align*}
 N(y)&:=(-1)( a^2-f_1(y)^2+y^2)f_1'(y)+2 f_1(y) y,
 \\
 D_\pm(y) &:= (-f_1(y)\pm  a)^2+y^2.
\end{align*}

It is easy to see that
\beq\label{eqn_Dy_estim}
D_\pm(y) \approx \Ge+y^2.
\eeq
As consequences of \eqnref{def_alpha} and \eqnref{fitaylor}, we have
$$
a^2=\frac{2\Ge}{\Gk_1+\Gk_2}+O(\Ge^2), \quad f_1(y)=\frac{\Gk_1}{\Gk_1+\Gk_2}\Ge+\frac{1}{2}\Gk_1 y^2+ O(\Ge^2+y^3),
$$
and hence
\begin{align}
N(y) &= (-1)\Big(\frac{2\Ge}{\Gk_1+\Gk_2}+y^2\Big)\Gk_1 y
\nonumber
\\
&\qquad +2\Big(\frac{\Gk_1 }{\Gk_1+\Gk_2}\Ge  + \frac{1}{2}\Gk_1 y^2\Big)y +O(\Ge y^3 + \Ge^2 y + y^4)
\nonumber
\\
&= O(\Ge y^2 + \Ge^2 y + y^4).
\label{eqn_Ny_estim}
\end{align}
Then, from \eqnref{xi1_bdry_deri_NDD} and the fact that $ a \approx \sqrt{\Ge}$, we have
\begin{align*}
\Big|\frac{d}{dy} \Gz(\Bx_1(y))\Big| \lesssim  a \frac{\Ge y^2 + \Ge^2 y+y^4 }{(\Ge+y^2)^2} \lesssim  \sqrt\Ge.
\end{align*}
So, \eqnref{eqn_bdry_deri1_1i} is proved.

Now let us consider \eqnref{eqn_bdry_deri2_1i}. We have
\begin{align}
\frac{d^2}{dy^2} \Gz(\Bx_1(y))
&=\frac{d^2}{dy^2} \Gz(-f_1(y),y)
=\frac{d}{dy}\Big( 2 a \frac{N(y)}{D_+(y) D_-(y)}\Big)
\nonumber\\
&=2 a\frac{N'(y)}{D_+(y) D_-(y)}-2 a\frac{N(y) (D_+)'(y)}{(D_+(y))^2 D_-(y)} -2 a
\frac{N(y) (D_-)'(y)}{ D_+(y) (D_-(y))^2}
\label{xi1_bdry_deri2_ND}.
\end{align}
We also have
\begin{align*}
N'(y) &= - a^2 f_1'' + 2f_1 (f_1')^2 + (f_1)^2 f_1''-y^2 f_1''+2 f_1,
\\
(D_\pm)'(y) &= 2(-f_1\pm  a)(-f_1')+2y.
\end{align*}
Since $f_1(y)\approx \Ge+y^2$, $f_1'(y)=O(y)$ and $f_1''(y)=O(1)$, we have
$$
|N'(y)| \lesssim \Ge+y^2, \quad |(D_\pm)'(y)|\lesssim|y|.
$$
Then, from \eqnref{eqn_Dy_estim}-\eqnref{xi1_bdry_deri2_ND} and the fact that $a\approx \sqrt\Ge$, we obtain
$$
\Big|\frac{d^2}{dy^2} \Gz(\Bx_1(y)) \Big| \lesssim  a \bigg( \frac{\Ge+y^2}{(\Ge+y^2)^2}+\frac{(\Ge y^2+\Ge^2 y+y^4)|y|}{(\Ge+y^2)^3}\bigg) \lesssim \frac{\sqrt{\Ge}}{\Ge+y^2}.
$$
The proof is completed.
\qed

%%%%%%%%%%%%%%%%%%%%%%%%%%%%%%%%%%%%%%%%%%%%%%%%%%%%%%%%%%%%%%%%%%%%
\subsection{Estimates of singular functions}\label{subsec:estim_singular}
%%%%%%%%%%%%%%%%%%%%%%%%%%%%%%%%%%%%%%%%%%%%%%%%%%%%%%%%%%%%%%%%%%%%

This subsection is to derive estimates of the singular function $\Bq_j$ in the narrow region $\Pi_{L}$ and on $\p D_1 \cup \p D_2$, which will be used in the later part of the paper.

We begin by showing that singular functions can be explicitly represented by the function $\Gz$ introduced in the previous subsection. Set
\beq\label{AGzBxdef}
A_\Gz(\Bx):= \left[1-\frac{\sinh^2\Gz(\Bx)}{a^2} y^2 \right]^{1/2}, \quad \Bx=(x,y).
\eeq
If $\Bx=(x,y)\in \Pi_{L}$, then it holds by \eqnref{GzBxest} that $|\Gz(\Bx)| \lesssim \sqrt{\Ge}$. Since $a\approx \sqrt{\Ge}$ by \eqnref{def_alpha}, we have
$$
\frac{\sinh^2\Gz(\Bx)}{a^2} \lesssim 1.
$$
So there exists a constant $0<L_0<L$ (independent of $\Ge$) such that
\beq\label{1sinh2Gza212}
 1-\frac{\sinh^2\Gz(\Bx)}{a^2} y^2 \geq \frac{1}{2}, \quad \Bx =(x,y) \in \Pi_{L_0}.
\eeq
Note that
\beq\label{AGzBx}
A_\Gz(\Bx) = 1 + O(y^2), \quad \Bx =(x,y) \in \Pi_{L_0}.
\eeq

Singular functions $\Bq_j$ can be represented in terms of $\Gz$ as follows

\begin{prop}\label{h_S1_cartesian}
Let $\Bq_i (\Bx)= (q_{i1}(\Bx), q_{i2}(\Bx))^T$ for $i=1,2$. If $\Bx\in\Pi_{L_0}$, then $q_{ij}$ are given by
\begin{align}
\ds q_{11}(\Bx) &= \Ga_1  \Gz(\Bx) - \Ga_2 A_\Gz(\Bx) \sinh\Gz(\Bx), \label{qoneone} \\
\ds q_{12}(\Bx)= q_{21}(\Bx) &= \Ga_2 a^{-1}y \sinh^2\Gz(\Bx), \label{qonetwo} \\
\ds q_{22}(\Bx) &= \Ga_1 \Gz(\Bx) + \Ga_2 A_\Gz(\Bx) \sinh\Gz(\Bx). \label{qtwotwo}
\end{align}
\end{prop}

\pf
From the definition \eqnref{Bqone} of $\Bq_1$ and the first identity in \eqnref{xi_cartesian}, we have
\begin{align}
q_{11}(\Bx) &= {\Ga_1} \Gz(\Bx)- {\Ga_2} \bigg[ \frac{(x+a)^2}{(x+a)^2+y^2}-\frac{(x-a)^2}{(x-a)^2+y^2}\bigg] \nonumber\\
&\quad + {\Ga_2 a} \bigg[ \frac{x+a}{(x+a)^2+y^2}+\frac{x-a}{(x-a)^2+y^2}\bigg] \nonumber\\
&={\Ga_1} \Gz(\Bx) - {\Ga_2} \bigg[ \frac{x(x+a)}{(x+a)^2+y^2}-\frac{x(x-a)}{(x-a)^2+y^2}\bigg] \nonumber\\
&= {\Ga_1} \Gz(\Bx) - {\Ga_2} \frac{2a x (a^2-x^2+y^2)}{((x-a)^2+y^2)((x+a)^2+y^2)}
\label{q11_xi_cartesian_rep}
\end{align}
for $\Bx\in \mathbb{R}^2\setminus\{\Bp_1,\Bp_2\}$.

Thanks to the first identity in \eqnref{xi_cartesian} again, we have
\begin{align}
\sinh\Gz(\Bx)&=\frac{1}{2}\left( \sqrt{\frac{(x+a)^2+y^2}{(x-a)^2+y^2}}- \sqrt{\frac{(x-a)^2+y^2}{(x+a)^2+y^2}}\right)
\nonumber
\\
&=\frac{2a x}{\sqrt{(x-a)^2+y^2}\sqrt{(x+a)^2+y^2}}
\label{sinh_1i_cartes}
\end{align}
for $\Bx\in \mathbb{R}^2\setminus\{\Bp_1,\Bp_2\}$. Then straightforward computations yield
\beq\label{A_1i_cartes}
A_\Gz(\Bx) =\frac{a^2-x^2+y^2}{\sqrt{((x+a)^2+y^2)((x-a)^2+y^2)}}, \quad \Bx\in\Pi_{L_0}.
\eeq
This together with \eqnref{q11_xi_cartesian_rep} yields \eqnref{qoneone}.

The identity \eqnref{qonetwo} can be proved similarly. In fact, one can see that
\begin{align}
q_{12}(\Bx) &= -{\Ga_2} \bigg[ \frac{(x+a)y}{(x+a)^2+y^2}-\frac{(x-a)y}{(x-a)^2+y^2}\bigg] \nonumber\\
& \quad + {\Ga_2 a} \bigg[ \frac{y}{(x+a)^2+y^2}+\frac{y}{(x-a)^2+y^2}\bigg] \nonumber\\
&={\Ga_2} \frac{4ax^2y}{((x+a)^2+y^2) ((x-a)^2+y^2)},
\label{q12_xi_cartesian_rep}
\end{align}
and \eqnref{qonetwo} follows from \eqnref{sinh_1i_cartes}.

Similarly one can show that $q_{21}=q_{12}$ and \eqnref{qtwotwo} hold. We omit the proof.
\qed

\medskip

Proposition \ref{h_S1_cartesian} already reveals an important property of the singular functions. They are almost constant near the points $\Bz_1$ and $\Bz_2$. This can be seen more clearly if two osculating disks have the same radii.

\begin{lemma} \label{lem:Bq_on_circle}
Assume $B_1$ and $B_2$ have the same radii $r_0$. Then
it holds for $\Bx \in \p B_i$, $i=1,2$, that
\begin{align}
\Bq_1(\Bx) &= \left(\frac{\sqrt\Ge}{m_1} +  t_1\right) \frac{(-1)^i}{2} \Psi_1+ \Ga_2\frac{a}{r_0^2}
\begin{bmatrix} x \\ y \end{bmatrix}, \label{qonesame}
\\
\Bq_2(\Bx) &= \left(\frac{\sqrt\Ge}{m_2}+t_2\right)\frac{(-1)^i}{2}\Psi_2 + \Ga_2\frac{a}{r_0^2}
\begin{bmatrix} y \\ x \end{bmatrix}, \label{qtwosame}
\end{align}
where $m_j$ are constants defined by \eqnref{m2def} and $t_j$ are constants satisfying
\beq\label{tj_estim}
|t_j| \le C (\Ga_1+\Ga_2)\Ge^{3/2}
\eeq
for some constant $C$ independent of $(\Ga_1,\Ga_2)$, or equivalently, independent of $(\Gl, \Gm)$, as well as $\Ge$.
\end{lemma}

\pf
If $r_1=r_2=r_0$, \eqnref{a_def} reads
\beq\label{asame}
a= \frac{\sqrt{\Ge(4r_0+\Ge)}}{2} .
\eeq
We see from \eqnref{xi_cartesian} that the constant value $\Gz(\Bx)$ on $\p B_j$ are as follows:
\beq\label{zeta_Bi_s}
\zeta|_{\p B_i} = (-1)^i \sinh^{-1} (a/r_0), \quad i=1,2.
\eeq
Let $s=\sinh^{-1} (a/r_0)$. Note that $a=r_0\sinh s$.
 Then it follows from \eqnref{asame} that
$$
r_0\cosh s = {\sqrt{r_0^2+a^2}} = {r_0+\Ge/2}.
$$
Since the center of $\p B_i$ is  $(-1)^i(r_0 + {\Ge}/{2}, 0)$, we have
$$
\p B_i = \big\{ (x,y)\in \Rbb^2: \big(x - (-1)^i r_0 \cosh s\big)^2 +y^2 = r_0^2 \big\}.
$$
So, for $(x,y) \in \p B_i$, we obtain
\begin{align}
(x\pm a)^2 +  y^2
&=x^2\pm 2x r_0 \sinh s + r_0^2\sinh^2 s +y^2
\nonumber
\\
&=(x - (-1)^i r_0 \cosh s)^2 +y^2 -r_0^2 +(-1)^i 2x r_0 \cosh s  \pm 2x r_0 \sinh s
\nonumber
\\
&= ((-1)^i\cosh s\pm \sinh s)2r_0x,
\label{xpma2py2}
\end{align}
and
\begin{align*}
a^2 -x^2 +y^2 &=
(x - (-1)^i r_0 \cosh s)^2+y^2-r_0^2 -2x^2 +(-1)^i 2xr_0 \cosh s
\\
&=((-1)^i r_0 \cosh s-x)2x =\big((-1)^i (r_0 + \Ge/2)-x\big)2x.
\end{align*}
Then, for $(x,y) \in \p B_i$, we have %from %\eqnref{A_1i_cartes} that
\begin{align}
%{A}_\Gz(\Bx)
% = \frac{((-1)^i (r_0 +\ds \frac{\Ge}{2})-x)2x}{\sqrt{(\cosh^2 s - \sinh^2 s) 4 r_0^2 x^2}}
\frac{2a x}{\sqrt{(x-a)^2+y^2}\sqrt{(x+a)^2+y^2}}
&=(-1)^i \frac{a}{r_0},
\nonumber
\\
\frac{a^2-x^2+y^2}{\sqrt{((x+a)^2+y^2)((x-a)^2+y^2)}}
&= 1+\frac{\Ge}{2r_0}-(-1)^i\frac{x}{r_0}.
\label{twofractions_pBi}
\end{align}

Note that $\zeta|_{\p B_i}=(-1)^is$.
So we obtain from \eqnref{q11_xi_cartesian_rep} and \eqnref{twofractions_pBi} that
$$
q_{11}(\Bx) = (-1)^i\Ga_1 s -\Ga_2 (-1)^i  (1+\frac{\Ge}{2r_0})\frac{a}{r_0} + \Ga_2\frac{a}{r_0}\frac{x}{r_0}, \quad \Bx \in \p B_i.
$$
Since $a=\sqrt{r_0 \Ge} + O(\Ge^{3/2})$, $s=\sqrt{\Ge/r_0} + O(\Ge^{3/2})$, and $m_1=\sqrt{r_0}/[2(\Ga_1-\Ga_2)]$, it follows that
\begin{align*}
q_{11}(\Bx) &=(-1)^i \left( (\Ga_1 - \Ga_2)
%\frac{\sqrt\Ge}{\sqrt{r_0}}
\sqrt{\Ge/r_0}
%\frac{\sqrt\Ge}{\sqrt{r_0}}
%\sqrt{\Ge/r_0}
 + (\Ga_1+\Ga_2)O(\Ge^{3/2}) \right) + \Ga_2\frac{a}{r_0^2}x \\
&=(-1)^i \Big( \frac{1}{2}\frac{\sqrt\Ge}{m_1} + (\Ga_1+\Ga_2)O(\Ge^{3/2}) \Big) + \Ga_2\frac{a}{r_0^2}x, \quad \Bx \in \p B_1.
\end{align*}
We also obtain from \eqnref{q12_xi_cartesian_rep} and \eqnref{xpma2py2} that
$$
q_{12}(\Bx) = \Ga_2 \frac{a}{r_0^2}y, \quad \Bx \in \p B_1 \cup \p B_2.
$$
This proves \eqnref{qonesame}. One can prove \eqnref{qtwosame} similarly.
\qed

We have the following lemma in $\Pi_{L_0}$.

\begin{lemma}\label{lem_hS_grad_estim}
We have, for $\Bx \in \Pi_{L_0}$,
\beq\label{poneest}
|\p_1 q_{11}(\Bx)| + |\p_1 q_{22}(\Bx)| \lesssim \frac{\sqrt\Ge}{\Ge+y^2},
\eeq
and
\beq\label{otherest}
|\p_2 q_{11} (\Bx)| + |\p_2 q_{22} (\Bx)| + | \nabla q_{12}(\Bx) | \lesssim \frac{\sqrt\Ge|y|}{\Ge+y^2}+ \sqrt\Ge.
\eeq
\end{lemma}

\pf
We only consider $\Bq_1$. Estimates for $\Bq_2$ can be obtained similarly.

First we consider $\p_1 q_{11}(\Bx)$. By \eqnref{qoneone}, we have
\begin{align}
 \p_1 q_{11} (\Bx) &= \Ga_1 \p_1 \Gz(\Bx) - \Ga_2 \cosh(\Gz(\Bx)) A_\Gz(\Bx) \p_1 \Gz(\Bx) \nonumber \\
&\quad + \Ga_2 \sinh^2(\Gz(\Bx))
\cosh(\Gz(\Bx)) \frac{y^2}{a^2}
A_\Gz(\Bx)^{-1}  \p_1 \Gz(\Bx) .
\label{eqn_dx_hS_1x}
\end{align}
Thanks to \eqnref{1sinh2Gza212}, we have
\beq\label{AGz_approx_1}
|A_\zeta(\Bx)|\approx 1.
\eeq
Then, using the fact that $a\approx\sqrt\Ge$, we obtain
$$
 |\p_1 q_{11} (\Bx)| \lesssim |\p_1 \Gz(\Bx)| + |\Gz(\Bx)|^2 \frac{y^2}{\Ge} |\p_1 \Gz(\Bx)|.
$$
We then infer from \eqnref{GzBxest} and \eqnref{Gzgradest} that
$$
|\p_1 q_{11}(\Bx)| \lesssim  \frac{\sqrt{\Ge}}{\Ge+y^2}
+ \Ge \frac{y^2}{\Ge}\frac{\sqrt{\Ge}}{\Ge+y^2}
\lesssim \frac{\sqrt\Ge}{\Ge+y^2}.
$$
This proves \eqnref{poneest}.

To prove \eqnref{otherest} for $\Bq_1$, we compute using \eqnref{qoneone} and \eqnref{qonetwo}
\begin{align*}
\p_2 q_{11} (\Bx) &= \Ga_1 \p_2 \Gz(\Bx)
- \Ga_2 \cosh(\Gz(\Bx)) A_\Gz(\Bx) \p_2 \Gz(\Bx)
\\
&\quad + \Ga_2 \sinh(\Gz(\Bx)) A_\Gz(\Bx)^{-1}
\Big(\sinh \Gz(\Bx) \p_2 \Gz(\Bx) \frac{y^2}{a^2}+\sinh^2\Gz(\Bx) \frac{y}{a^2}\Big),
\end{align*}
and
\begin{align*}
 \p_1 q_{12}(\Bx) &= 2\Ga_2 a y \sinh\Gz(\Bx) \cosh\Gz(\Bx) \p_1 \Gz(\Bx), \\
 \p_2 q_{12}(\Bx) &= 2\Ga_2 a y \sinh\Gz(\Bx) \cosh\Gz(\Bx) \p_2 \Gz(\Bx) + \Ga_2 a \sinh^2\Gz(\Bx).
\end{align*}
So, \eqnref{otherest} can be proved in the same way as above.
The proof is completed.
\qed

\smallskip

Let $\Bh_1=(h_{11},h_{12})^T$ be the solution to \eqnref{hj_def} for $j=1$. Then,
we have $\Bh_{1}(\Bx_2(y))-\Bh_{1}(\Bx_1(y)) = (1,0)^T$. Since $|\Bx_2(0)- \Bx_1(0)|=\Ge$,
one can expect
\beq\label{h1_temp_estim}
\begin{cases}
\p_1 h_{11}(0,0) = \Ge^{-1} + O(1), \\
|\p_2 h_{11} (0,0)| + |\p_1 h_{12} (0,0)| + |\p_2 h_{12} (0,0)| \lesssim 1.
\end{cases}
\eeq
One can expect a similar behavior for $\Bh_2$ as well.
We now show that $\frac{m_j}{\sqrt\Ge}\Bq_j$ has the exactly same behavior as $\Ge \to 0$.

\begin{lemma}\label{singular_q_origin}
It holds for small $\Ge>0$ that
\beq\label{p1q11}
\begin{cases}
\p_1 q_{11} (0,0) = \ds \frac{1}{m_1\sqrt\Ge} + O(\sqrt\Ge), \\
|\p_2 q_{11} (0,0)| + |\p_1 q_{12} (0,0)| + |\p_2 q_{12} (0,0)| \lesssim \sqrt\Ge,
\end{cases}
\eeq
and
\beq\label{p1q22}
\begin{cases}
\p_1 q_{22} (0,0) = \ds \frac{1}{m_2\sqrt\Ge} + O(\sqrt\Ge), \\
|\p_1 q_{21} (0,0)| + |\p_2 q_{21} (0,0)| + |\p_2 q_{22} (0,0)| \lesssim \sqrt\Ge.
\end{cases}
\eeq
\end{lemma}
\pf
Since $\p_1 \Gz(0,0) = 2/a$ and $\Gz(0,0)=0$, it follows from \eqnref{eqn_dx_hS_1x} that
$$
\p_1 q_{11} (0,0) = \frac{2(\Ga_1-\Ga_2)}{a}.
$$
Since $a = \sqrt{2\Ge}/\sqrt{(\Gk_1+\Gk_2)}+O(\Ge^{3/2})$, the first equality in \eqnref{p1q11} follows.
From \eqnref{otherest}, we have
$$
|\p_2 q_{11} (0,0)|+| \p_1 q_{12}(0,0) |+| \p_2 q_{12}(0,0) | \lesssim \sqrt\Ge.
$$
This proves \eqnref{p1q11}. \eqnref{p1q22} can be proved similarly.
\qed

\begin{lemma}\label{lem_qj_far_estim}
For $j=1,2$, we have
\beq\label{Bqjest}
\|\Bq_j\|_{L^\infty(D^e\setminus \Pi_{L_0})}+\| \nabla \mathbf{q}_{j} \|_{L^\infty(D^e\setminus \Pi_{L_0})} \lesssim \sqrt\Ge.
\eeq
\end{lemma}
\pf
We only prove \eqnref{Bqjest} for $j=1$. The same proof applies to the case when $j=2$.

Recall that
$$
\Bq_1(\Bx) = \BGG(\Bx-\Bp_1)\Be_1-\BGG(\Bx-\Bp_2) \Be_1 + {\Ga_2 a} \left( \frac{\Bx-\Bp_1}{|\Bx-\Bp_1|^2}+\frac{\Bx-\Bp_2}{|\Bx-\Bp_2|^2}\right).
$$
Note that if $\Bx\in D^e \setminus \Pi_{L_0}$, then $1 \lesssim |\Bx-\Bp|$ for all $\Bp$ lying on the line segment $\ol{\Bp_1\Bp_2}$.
Since $a \approx \sqrt{\Ge}$, the second term on the right-hand side of the above and its derivative is less than $\sqrt{\Ge}$.

One can easily show that the first term also satisfies the same estimate. In fact, by the mean value theorem, we have
\begin{align*}
|\BGG(\Bx-\Bp_1)-\BGG(\Bx-\Bp_2)| \lesssim |\nabla \BGG(\Bx-\Bp_*)| |\Bp_1-\Bp_2|
\end{align*}
for some $\Bp_*$ on $\ol{\Bp_1\Bp_2}$. We also have
\begin{align*}
|\nabla(\BGG(\Bx-\Bp_1)-\BGG(\Bx-\Bp_2))| \lesssim |\nabla^2 \BGG(\Bx-\Bp_{**})| |\Bp_1-\Bp_2|
\end{align*}
for some $\Bp_{**}$ on $\ol{\Bp_1\Bp_2}$. Since $|\Bp_1-\Bp_2|=2a\approx \sqrt\Ge$, \eqnref{Bqjest} follows.
\qed

\medskip

As a corollary, we have the following estimate for $\nabla\Bq_j$.

\begin{cor}\label{cor_hS_blowup_estim}
For $j=1,2$, we have
\beq\label{Bq12est}
\|\nabla \Bq_j\|_{L^\infty(D^e)} \approx \Ge^{-1/2}.
\eeq
\end{cor}

\pf
The upper estimate $\|\nabla \Bq_j\|_{L^\infty(D^e)} \lesssim \Ge^{-1/2}$ is a consequence of Lemma \ref{lem_hS_grad_estim} and \ref{lem_qj_far_estim}, and the lower one is that of Lemma \ref{singular_q_origin}.
\qed

\medskip

We have the following lemma on $\p D_1\cup \p D_2$.

\begin{lemma}\label{lem_hS_bdry_estim}
Let $\Bx_k(y)$ be the defining functions for $\p D_k$ for $k=1,2$ as given in \eqnref{Bx1Bx2}. For $|y|<L_0$, the following holds:
\begin{align}
q_{11} (\Bx_k(y)) &= (-1)^k {(\Ga_1-\Ga_2) \Gk_k a} + O \left( E \right), \label{q11asym}
\\
q_{12}(\Bx_k(y))=q_{21}(\Bx_k(y)) &= {\Ga_2 \Gk_k^2 ay} + O \left( |y| E \right), \label{q12asym} \\
q_{22} (\Bx_k(y)) &= (-1)^k {(\Ga_1+\Ga_2) \Gk_k a} + O \left( E \right), \label{q22asym}
\end{align}
where
\beq\label{Edef}
E:= \Ge^{3/2} + \sqrt\Ge y^2 + \tau \sqrt\Ge |y|.
\eeq
\end{lemma}

\pf
We see from \eqnref{xi_cartesian} that
$$
\Gz|_{\p B_k} =(-1)^i\sinh^{-1} (\Gk_k a)= (-1)^i \Gk_k a + O(a^3).
$$
Since $a\approx \sqrt{\Ge}$, we infer from \eqnref{difference} that
\beq\label{101}
\Gz(\Bx_k(y))=  (-1)^i\Gk_k a +O \left(  E \right),
\eeq
and
\beq\label{102}
\sinh \Gz(\Bx_k(y))=  (-1)^i\Gk_k a +O \left( E \right).
\eeq
Combining \eqnref{AGzBx}, \eqnref{101} and \eqnref{102}, one can see that  \eqnref{q11asym}, \eqnref{q12asym} and \eqnref{q22asym} follow from \eqnref{qoneone}, \eqnref{qonetwo}, and \eqnref{qtwotwo}, respectively.
\qed

\medskip

Then, using \eqnref{def_alpha} and the definitions \eqnref{m2def} of $m_1$ and $m_2$, we immediately obtain the following corollary.
\begin{cor}\label{cor_hS_D1D2_diff_estim}
For $|y|<L_0$, we have
\begin{align}
q_{11}(\Bx_2(y)) - q_{11}(\Bx_1(y)) &=  m_1^{-1} \sqrt\Ge +O(E),
\\
q_{12}(\Bx_2(y)) - q_{12}(\Bx_1(y)) &= q_{21}(\Bx_2(y)) - q_{21}(\Bx_1(y))= O(\sqrt\Ge\tau|y|), \\
q_{22}(\Bx_2(y)) - q_{22}(\Bx_1(y)) &= m_2^{-1}\sqrt\Ge +O(E),
\end{align}
where $E$ is given by \eqnref{Edef}.
\end{cor}

\smallskip
We then obtain the following lemma for estimates of the derivatives of $\Bq_j$.

\begin{lemma} \label{lem_hS_bdry_deri_estim}
Let $\Bx_k(y)$ be the defining functions for $\p D_k$ for $k=1,2$ as given in \eqnref{Bx1Bx2}. For $|y|<L_0$, the following holds:
\begin{align}
\Big| \frac{d}{dy} \Bq_{j}(\Bx_k(y)) \Big| & \lesssim \sqrt\Ge, \label{qderiest1}
\\
\Big|\frac{d^2}{dy^2} q_{11}(\Bx_k(y)) \Big| + \Big|\frac{d^2}{dy^2} q_{22} (\Bx_k(y))\Big| & \lesssim \frac{\sqrt\Ge}{{\Ge}+ y^2}, \label{qderiest2} \\
\Big|\frac{d^2}{dy^2} q_{12} (\Bx_k(y))\Big| & \lesssim \sqrt\Ge. \label{qderiest3}
\end{align}
\end{lemma}

\pf
We only prove inequalities corresponding to $\Bq_1 (\Bx_1(y))$. Those for other cases, namely, $\Bq_1 (\Bx_2(y))$ and $\Bq_2 (\Bx_k(y))$, can be treated similarly.

For ease of notation, let us define
$\Gvf(y)$ and $\Phi(y)$ by
$$
\Gvf(y):=\Gz(\Bx_1(y)), \quad
\Phi(y):= A_\Gz (\Bx_1(y)).
$$
We see from \eqnref{difference}-\eqnref{eqn_bdry_deri2_1i} that
\beq\label{Gvfest}
|\Gvf(y)| \lesssim \sqrt\Ge, \quad |\Gvf'(y)| \lesssim\sqrt\Ge, \quad |\Gvf''(y)| \lesssim \frac{\sqrt\Ge}{\Ge+y^2}.
\eeq
We also have
\beq\label{Phiest}
|\Phi(y)| \approx 1, \quad |\Phi'(y)| \lesssim |y|, \quad |\Phi''(y)| \lesssim 1.
\eeq
The first estimate in the above is an immediate consequence of \eqnref{AGz_approx_1}, and the last two can be proved using the definition \eqnref{AGzBxdef} of $A_\Gz(\Bx)$. In fact, straightforward computations yield
$$
\ds\Phi'(y)=-\frac{1}{2a^2 \Phi} \left(y^2 \Gvf' \sinh 2\Gvf + 2y\sinh^2 \Gvf\right),
$$
and
\begin{align*}
\Phi''(y)&=
\ds-\frac{1}{2a^2 \Phi} \big( 4 y \Gvf' \sinh 2\Gvf + y^2 \Gvf'' \sinh 2\Gvf + 2y^2(\Gvf')^2\cosh 2\Gvf + 2\sinh^2 \Gvf \big)
\\
&\qquad
+\frac{\Phi'}{2a^2 \Phi^2} \left(y^2 \Gvf' \sinh 2\Gvf + 2y\sinh^2 \Gvf\right).
%+\frac{\Phi'}{2a^2 \Phi^2} \left(y^2 \Gvf' \sinh 2\Gvf + 2y\sinh^2 \Gvf\right).
\end{align*}
Then, using \eqnref{Gvfest} and the fact that $|\Phi|\approx 1$, we obtain
\begin{align*}
|\Phi'(y)| &\lesssim \frac{1}{\Ge}( y^2\Ge + |y|{\Ge})
\lesssim |y|,
\\[0.5em]
|\Phi''(y)| &\lesssim \frac{1}{\Ge} \Big(|y|\Ge + y^2 \frac{\sqrt\Ge}{\Ge+y^2}\sqrt\Ge + y^2 \Ge + \Ge \Big) + \frac{|y|}{\Ge} \big(y^2 \Ge + |y| \Ge \big)\lesssim 1.
\end{align*}

We have from \eqnref{qoneone} and \eqnref{qonetwo} that
\begin{align*}
 \frac{d}{dy} q_{11} (\Bx_1(y))&= \Ga_1 \Gvf' - \Ga_2 \left( \Phi \Gvf'\cosh \Gvf  + \Phi'\sinh \Gvf  \right),
\\
 \frac{d}{dy} q_{12} (\Bx_1(y))&= \Ga_2 a^{-1} (\sinh^2 \Gvf + y \Gvf' \sinh 2 \Gvf),
\\
 \frac{d^2}{dy^2} q_{11}(\Bx_1(y))&=  \Ga_1 \Gvf'' -
 \Ga_2\big( (\Phi' \Gvf'  + \Phi \Gvf''  + \Phi' \Gvf')\cosh\Gvf+ (\Phi \Gvf'^2 + \Phi'') \sinh\Gvf \big),
% \\
% &\quad\frac{\Ga_2}{4\Phi^{3/2}} \sinh \Gvf \Big( (2 \Phi \Gvf')^2 +2 \Phi \Phi''-(\Phi')^2\Big) \\& \qquad -\frac{\Ga_2}{\Phi^{1/2}} \cosh \Gvf \big( \Phi \Gvf'' + \Phi' \Gvf'\big),
\\
 \frac{d^2}{dy^2} q_{12} (\Bx_1(y))
&= \Ga_2 a^{-1} (2\Gvf'\sinh 2\Gvf + y\Gvf'' \sinh2\Gvf + 2y\Gvf'^2  \cosh 2\Gvf ).
\end{align*}
Since $a \approx \sqrt{\Ge}$, \eqnref{qderiest1}-\eqnref{qderiest3} now follow from \eqnref{Gvfest} and \eqnref{Phiest}.
\qed

\medskip

We now estimate $\Bq_3$ whose behavior resembles that of the solution $\Bh_3=(h_{31},h_{32})^T$ to \eqnref{hj_def} for $j=3$.
Since $\Bh_3|_{\p D_i} = \frac{(-1)^i}{2}(-y,x)^T$ for $i=1,2$, we see that
$$
\Bh_{3}(\Bx_2(y))-\Bh_{3}(\Bx_1(y)) = \left( -y, \frac{f_1(y) + f_2(y)}{2} \right).
$$
Since $|\Bx_2(y)-\Bx_1(y)|=f_1(y)+f_2(y)$,
one can expect that the following holds for small $\Ge>0$ and for $(x,y)$ near the origin:
\beq\label{h3_temp_estim}
\p_1 h_{31} \approx \frac{-y}{f_1(y)+f_2(y)} + O(1) =  -\frac{y}{\Ge+\frac{1}{2}(\Gk_1+\Gk_2)y^2} + O(1)
\eeq
and
\beq
|\p_2  h_{31}(x,y)|+|\p_1  h_{32}(x,y)|+|\p_2  h_{32}(x,y)|  \lesssim 1.
\eeq
The following lemma shows that $\Bq_3$ has the exactly same local behavior.

\begin{lemma} \label{lem_q3_asymp}
For $\Bx =(x,y)\in \Pi_L$, we have
\beq
\p_1  q_{31}(\Bx) =   -\frac{y}{\Ge+\frac{1}{2}(\Gk_1+\Gk_2)y^2} +O(1),
\eeq
and
\beq
|\p_2  q_{31}(\Bx)|+|\p_1  q_{32}(\Bx)|+|\p_2  q_{32}(\Bx)|  \lesssim 1.
\eeq
\end{lemma}
\pf
Let $m_3$ be the number defined by \eqnref{m3def}.
For ease of computation, we decompose $\Bq_3$ as $\Bq_3=\widetilde{\Bq}_3+\Bw$ where $\Bw(\Bx) = -m_3(\BGG^\perp(\Bx-\Bc_1)+\BGG^\perp(\Bx-\Bc_2))$. It is clear that $|\nabla \Bw(\Bx)| \lesssim 1$ for $\Bx\in\Pi_L$.

Now we consider $\widetilde{\Bq}_3=(\tilde{q}_{31},\tilde{q}_{32})^T$, which is given by
$$
\widetilde{\Bq}_3=m_3\left(\BGG^\perp(\Bx-\Bp_1)+\BGG^\perp(\Bx-\Bp_2)\right)+m_3 \Ga_2 a \left( \frac{(\Bx-\Bp_1)^\perp}{|\Bx-\Bp_1|^2} - \frac{(\Bx-\Bp_2)^\perp}{|\Bx-\Bp_2|^2} \right).
$$
From the definition \eqnref{Kelvin_modified} of $\BGG^\perp$, we have
\begin{align*}
\tilde{q}_{31}(\Bx) &=  m_3\Ga_1 \sum_{i=1}^2 \arg\big(x-(-1)^{i}a+iy\big)
\\
&\quad -m_3\Ga_2 \sum_{i=1}^2\frac{(x-(-1)^{i}a)(-y) }{(x-(-1)^i a)^2+y^2}  + m_3\Ga_2 a \sum_{i=1}^2 \frac{(-1)^{i+1} (-y)}{(x-(-1)^ia)^2+y^2},
\\
\tilde{q}_{32}(\Bx) &= -m_3\Ga_2 \sum_{i=1}^2 \frac{(x-(-1)^{i}a)^2 }{(x-(-1)^{i}a)^2+y^2} + m_3\Ga_2 a \sum_{i=1}^2 \frac{(-1)^{i+1}(x-(-1)^{i}a)}{(x-(-1)^{i}a)^2+y^2}.
\end{align*}
Straightforward computations yield
\begin{align*}
\p_1 \tilde{q}_{31}(\Bx) &= -(\Gk_1+\Gk_2)^{-1} f(\Bx)
- 2m_3\Ga_2 xy \big[x h_+(\Bx) + a h_-(\Bx)\big], \\
\p_2 \tilde{q}_{31}(\Bx) &=
m_3(\Ga_1+\Ga_2)x g_+(\Bx)+m_3 \Ga_1 a g_-(\Bx)-2m_3\Ga_2 x y^2 h_+(\Bx),
%%%
%m_3(\Ga_1+2\Ga_2)a  g_-(\Bx) + m_3(\Ga_1+\Ga_2)x g_+(\Bx)
%\\
%& \qquad -  2m_3\Ga_2 y^2  \big[x h_+(\Bx)+2a h_-(\Bx)\big],
\\
\p_1 \tilde{q}_{32}(\Bx) &= -2m_3 \Ga_2 x g_+(\Bx)-m_3\Ga_2 a g_-(\Bx)
\\
&\qquad
 + 2m_3 \Ga_2 x \big[ (x^2+a^2) h_+(\Bx)+2ax h_-(\Bx)\big],
%%% - 2m_3\Ga_2 y^2 \big[x  h_+(\Bx) + 2a  h_-(\Bx) \big],
\\
\p_2 \tilde{q}_{32}(\Bx) &= 2m_3\Ga_2 xy \big[x h_+(\Bx) + a h_-(\Bx)\big],
\end{align*}
where $f$, $g_\pm$ and $h_\pm$ are defined by
\begin{align*}
f(\Bx)&=\frac{y}{(x+a)^2+y^2}+\frac{y}{(x-a)^2+y^2},
\\
g_\pm(\Bx)&=\frac{1}{(x+a)^2 + y^2}\pm\frac{1}{(x-a)^2 + y^2},
\\
h_\pm(\Bx)&=\frac{1}{((x+a)^2 + y^2)^2}\pm\frac{1}{((x-a)^2 + y^2)^2}.
\end{align*}

Since $a\approx \sqrt\Ge$, $|x|\lesssim \Ge+y^2$ and $(x\pm a)+y^2 \approx \Ge+y^2$, we see that
\begin{align*}
&|g_+(\Bx)| \lesssim \frac{1}{\Ge+y^2}, \quad |h_+(\Bx)| \lesssim \frac{1}{(\Ge+y^2)^2},
\\
&|g_-(\Bx)| = \left| \frac{4ax}{((x+a)^2+y^2)((x-a)^2+y^2)} \right| \lesssim \frac{\sqrt\Ge}{\Ge+y^2},
\\
&|h_-(\Bx)| = \left| \frac{4ax((x+a)^2+(x-a)^2+2y^2)}{((x+a)^2+y^2)^2((x-a)^2+y^2)^2} \right| \lesssim \frac{\sqrt\Ge}{(\Ge+y^2)^2}.
\end{align*}
Therefore, we obtain
\begin{align*}
|&\p_1 \tilde{q}_{31}(\Bx) + (\Gk_1+\Gk_2)^{-1} f(\Bx)| \lesssim \frac{(\Ge+y^2)^2 y}{(\Ge+y^2)^2} + \frac{\Ge(\Ge+y^2)y}{(\Ge+y^2)^2}\lesssim 1,
\\
|&\p_2 \tilde{q}_{31}(\Bx)|  \lesssim  \frac{\Ge+y^2}{\Ge+y^2}+\sqrt\Ge\frac{ \sqrt\Ge}{\Ge+y^2} +
 + \frac{(\Ge+y^2)y^2}{(\Ge+y^2)^2}  \lesssim 1,
\\
|&\p_1 \tilde{q}_{32}(\Bx)|  \lesssim  \frac{\Ge+y^2}{\Ge+y^2} +\sqrt\Ge\frac{ \sqrt\Ge}{\Ge+y^2}
  + \frac{(\Ge+y^2)^3 +\Ge (\Ge+y^2)}{(\Ge+y^2)^2}+\frac{\Ge(\Ge+y^2)}{(\Ge+y^2)^2} \lesssim 1, \\
|&\p_2 \tilde{q}_{32}(\Bx)| \lesssim \frac{(\Ge+y^2)^2 y}{(\Ge+y^2)^2} + \frac{\Ge(\Ge+y^2)y}{(\Ge+y^2)^2}\lesssim 1.
\end{align*}

Now it remains to show that
\beq\label{f_asymp}
f(\Bx)  = \frac{(\Gk_1+\Gk_2) y}{\Ge+\frac{1}{2}(\Gk_1+\Gk_2)y^2} + O(1).
\eeq
Since
\begin{align*}
\left| \frac{y}{(x\pm a)^2 + y^2}-\frac{ y}{a^2+y^2} \right|
&\lesssim \left| \frac{y(x^2\pm 2ax)}{((x\pm a)^2 + y^2)(a^2+y^2)} \right|
\\
&\lesssim \frac{y((\Ge+y^2)^2+\sqrt\Ge(\Ge+y^2))}{(\Ge+y^2)^2}\lesssim 1,
\end{align*}
we see that
$$
f(\Bx) = \frac{2y}{a^2+y^2} + O(1).
$$
Since
$$
\frac{y}{a^2+y^2}=\frac{y}{{2\Ge}/{(\Gk_1+\Gk_2)}+O(\Ge^2)+y^2}
=\frac{\frac{1}{2}(\Gk_1+\Gk_2) y}{\Ge+\frac{1}{2}(\Gk_1+\Gk_2)y^2}+O(1),
$$
the desired estimate \eqnref{f_asymp} follows. This completes the proof.
\qed

\begin{lemma}\label{cor_q3_near_origin}
The following holds:
\beq\label{q3est1}
|\nabla \Bq_3(0,0)| \lesssim 1,
\quad
\nabla \Bq_3(0,a) = \frac{\sqrt{\Gk_1+\Gk_2}}{\sqrt{2\Ge}} \Be_1\otimes\Be_1 + O(1).
\eeq
Moreover, we have
\beq\label{Bq3est2}
\|\nabla \Bq_{3}\|_{L^\infty(D^e \setminus \Pi_L)} \lesssim 1,
\eeq
and
\beq\label{Bq3est}
\|\nabla \Bq_{3}\|_{L^\infty(D^e)} \approx \frac{ 1}{\sqrt{\Ge}}.
\eeq
\end{lemma}

\pf
The estimates in \eqnref{q3est1} are consequences of Lemma \ref{lem_q3_asymp}. The estimate \eqnref{Bq3est} is a consequence of Lemma \ref{lem_q3_asymp} and \eqnref{Bq3est2}.

To prove \eqnref{Bq3est2}, recall that
\begin{align*}
\Bq_3(\Bx) &= m_3\left( \BGG^\perp(\Bx-\Bp_1) - \BGG^\perp(\Bx-\Bc_1) \right) \Be_1 + m_3\left( \BGG^\perp(\Bx-\Bp_2) - \BGG^\perp(\Bx-\Bc_2) \right) \Be_1 \nonumber \\
&\quad -m_3\Ga_2 a \left( \frac{(\Bx-\Bp_1)^\perp}{|\Bx-\Bp_1|^2} - \frac{(\Bx-\Bp_2)^\perp}{|\Bx-\Bp_2|^2} \right).
\end{align*}
If $\Bx\in D^e \setminus \Pi_{L}$, then $1 \lesssim |\Bx-\Bc|$ for all $\Bc$ on the line segment $\ol{\Bc_1\Bc_2}$.
Note that $\Bp_1$ and $\Bp_2$ are on $\ol{\Bc_1\Bc_2}$. So, all the terms in parentheses above and their gradients are bounded. So, \eqnref{Bq3est2} follows.
\qed

%%%%%%%%%%%%%%%%%%%%%%%%%%%%%%%%%%%%%%%%%%%%%%%%%%%%%%%%%%%%%%%
\subsection{Approximations by singular functions}\label{subsec:approx_energy_rjrj}
%%%%%%%%%%%%%%%%%%%%%%%%%%%%%%%%%%%%%%%%%%%%%%%%%%%%%%%%%%%%%%%

In this section we prove \eqnref{BhjBqj}. More precisely, we prove the following proposition.

\begin{prop}\label{prop_Brj}
For $j=1,2$, let $\Bh_j$ be the solution to \eqnref{hj_def} in $\Acal^*$ and $m_j$ be the constant defined in \eqnref{m2def}. Then it holds that
\beq\label{hj_def_decomp}
\Bh_j = \frac{m_j}{\sqrt\Ge}\Bq_j + \Br_j,
\eeq
where $\nabla \Br_j$ satisfies
\beq\label{BrjestDe}
\int_{D^e} \Cbb \hatna\Br_j:\hatna\Br_j \lesssim 1.
\eeq
\end{prop}

To prove Proposition \ref{prop_Brj} we apply the variational principle. We emphasize that this is possible only because the singular function $\Bq_j$ is the solution to the Lam\'e system, namely, $\Lcal_{\Gl,\Gm} \Bq_j= 0$ in $D^e$, and so is $\Br_j$. Note that $\Br_j \in \Acal^*$ and
$$
\Br_j= (-1)^i\frac{1}{2} \Psi_j - \frac{m_j}{\sqrt\Ge}\Bq_j \quad \mbox{ on } \p D_i, \quad i=1,2.
$$
Let
\beq\label{Wjdef}
W_j = \left\{ \Bv \in \Acal^* ~|~ \Bv|_{\p D_i} = (-1)^i\frac{1}{2}\Psi_j -\frac{m_j}{\sqrt\Ge}\Bq_j \right\},
\eeq
and let $\Ecal_{D^e}$ be the energy functional defined in \eqnref{eqn_def_Ecal}.
By the variational principle \eqnref{variation}, we have
\beq\label{rj_var_prin}
\Ecal_{D^e}[\Br_j] =\min_{\Bv\in W_j} \Ecal_{D^e}[\Bv].
\eeq

We define the test function $\Br_j^K$ as follows: for $(x,y)\in \Pi_{L_0}$ let
\begin{align}\label{eqn_def_rKj_PiL0}
\ds\Br^K_j(x,y) &:=
\frac{\Br_j(\Bx_2(y))-\Br_j(\Bx_1(y))}{f_1(y)+f_2(y)}[ x+f_1(y)] + \Br_j(\Bx_1(y)).
\end{align}
Note that
\beq\label{rKbdry}
\Br^K_j = (-1)^i\frac{1}{2}\Psi_j -\frac{m_j}{\sqrt\Ge}\Bq_j = \Br_j \quad\mbox{on } \p D_i \cap \p \Pi_{L_0}, \ \ i=1,2,
\eeq
and $\Br_j^K$ is a linear interpolation of $\Br_j|_{\p D^e}$ in the $x$-direction. So, in $\Pi_{L_0}$, $\Br_j^K(x,y)$ is a linear function of $x$ for each fixed $y$.
Let $B_0$ be a disk containing $\overline{D_1\cup D_2}$, and extend $\Br^K_j$ to $D^e \setminus \Pi_{L_0}$ so that $\Br_j^K|_{\Rbb^2\setminus B_0}=0$, $\|\Br_j^K\|_{H^1(D^e\setminus \Pi_{L_0})} \lesssim 1$, and the boundary condition \eqnref{rKbdry} holds on $\p D_i$ for $i=1,2$.
Then, $\Br^K_j$ belongs to $W_j$.

We have the following lemma.
\begin{lemma}\label{lem_rjK_grad_estim}
We have, for $(x,y)\in \Pi_{L_0}$,
\beq\label{Brjest}
|\nabla \Br_j^K(x,y)| \lesssim 1+\frac{\tau|y|}{\Ge+y^2}.
\eeq
\end{lemma}

\pf
We prove \eqnref{Brjest} for $j=1$. The case for $j=2$ can be proved in a similar way.

Let us write $\Br_j^{K}(\Bx)= (r^K_{j1}(\Bx),r^K_{j2}(\Bx) )^T$. To keep the expressions simple, we introduce
\begin{align}
d(y)&:=f_1(y)+f_2(y),\nonumber
\\
 \phi(y)&:=1-\frac{m_1}{\sqrt\Ge}\big[ q_{11}(\Bx_2(y))-{q}_{11}(\Bx_1(y))\big],\nonumber
\\
 \Gn(y) &:= -\frac{1}{2}-\frac{m_1}{\sqrt\Ge} q_{11}(\Bx_1(y)).
\label{eqn_def_d_phi}
\end{align}
Then $r^K_{11}$ can be rewritten as
$$
r^K_{11}(x,y) = \frac{\phi(y)}{d(y)} x +\frac{\phi(y)f_1(y)}{d(y)} +\Gn(y), \quad (x,y)\in \Pi_{L_0}.
$$
Straightforward computations show that
\begin{align}
\p_1 r^K_{11} &= \frac{\phi}{d}, \label{dx_rK_1x}
\\[0.5em]
\p_2 r^K_{11} &= \left[\frac{\phi'}{d}-\frac{  \phi d'}{d^2} \right]x +\frac{\phi' f_1}{d}  + \frac{\phi f_1' }{d}-\frac{  \phi f_1 d'}{d^2} +\eta'.
\label{eqn_dy_rK_1x}
\end{align}

Note that
\beq\label{eqn_estim_d}
d(y) \approx \Ge +y^2, \quad  |d'(y)|\lesssim|y| ,\quad  |d''(y)|\lesssim 1.
\eeq
Note also that, from \eqnref{m2def}, Corollary \ref{cor_hS_D1D2_diff_estim} and Lemma \ref{lem_hS_bdry_deri_estim}, we have
\beq\label{eqn_estim_phi}
|\phi(y)| \lesssim \Ge+ y^2+\tau|y| ,\qquad |\phi'(y)|,|\Gn'(y)| \lesssim 1.
\eeq
From \eqnref{dx_rK_1x}-\eqnref{eqn_estim_phi} and the fact that $|x|\lesssim \Ge +y^2 $ for $(x,y)\in\Pi_L$, we have
\begin{align*}
|\p_1 r^K_{11}| &\lesssim \frac{\Ge + y^2 + \tau|y|}{\Ge+y^2} \lesssim 1+\frac{ \tau|y|}{\Ge+y^2},
\\[0.5em]
|\p_2 r^K_{11}| &\lesssim \left[ \frac{1}{\Ge+y^2}+\frac{(\Ge+|y|)|y|}{(\Ge+y^2)^2}\right] (\Ge + y^2)
\\
&\qquad
+  \frac{\Ge+y^2}{\Ge+y^2} + \frac{(\Ge+|y|) |y|}{\Ge+y^2}+\frac{(\Ge+|y|)(\Ge+y^2)|y|}{(\Ge+y^2)^2} + 1
 \lesssim 1.
\end{align*}
In a similar way, one can see that
\begin{align*}
|\p_1 r^K_{12}| &\lesssim 1+\frac{ \tau|y| }{\Ge+y^2},
\quad
|\p_2 r^K_{12}| \lesssim 1.
\end{align*}
This completes the proof.
\qed

\medskip
\noindent{\sl Proof of Proposition \ref{prop_Brj}}.
By the variational principle \eqnref{rj_var_prin}, we have
\begin{align*}
\Ecal_{D^e}[\Br_j] \le \Ecal_{D^e}[\Br_j^K]
\lesssim \|\nabla\Br_j^K \|_{L^2(D^e)}^2 .
\end{align*}
It follows from Lemma \ref{lem_rjK_grad_estim} that
\begin{align*}
\int_{D^e} |\nabla\Br^K_j|^2
&\lesssim\int_{\Pi_{L_0}} |\nabla\Br^K_j|^2 + \int_{D^e\setminus \Pi_{L_0}} |\nabla\Br^K_j|^2
\\
&\lesssim\int_{-L_0}^{L_0}\int_{-f_1(y)}^{f_2(y)} \Big(\frac{\Ge+|y|}{\Ge+y^2}\Big)^2 \,dx\, dy +1
\\
&\lesssim \int_{-L_0}^{L_0} \frac{(\Ge+|y|)^2}{\Ge+y^2} \,dy +1
\lesssim 1.
\end{align*}
So the proof is complete.
\qed

%%%%%%%%%%%%%%%%%%%%%%%%%%%%%%%%%%%%%%%%%%%%%%%%%%%%%%%%%%%%%%%%%%%%
\section{Stress concentration-boundary value problem}\label{sec:bvp}
%%%%%%%%%%%%%%%%%%%%%%%%%%%%%%%%%%%%%%%%%%%%%%%%%%%%%%%%%%%%%%%%%%%%

This section deals with the stress concentration, {\it i.e.}, the gradient blow-up of the solution to the boundary value problem \eqnref{elas_eqn_bdd}. We characterize the stress concentration in the narrow region between two inclusions in terms of the singular functions $\Bq_j$ defined in \eqnref{Bqone}, \eqnref{Bqtwo}, and \eqnref{q3def}. The main results (Theorem \ref{main_thm_2_general_bdd} and \ref{cor-bdd}) are stated and proved in subsection \ref{subsec:bvp}. Preliminary results required for proving main ones are also stated in the same subsection. Their proofs are given in subsequent subsections. At the end of subsection \ref{subsec:bvp} we include a brief comparison of this paper's method with that of \cite{BLL-ARMA-15} where the upper bound of the gradient blow-up is obtained.

%%%%%%%%%%%%%%%%%%%%%%%%%%%%%%%%%%%%%%%%%%%%%%%%%%%%%%%%%
\subsection{Characterization of stress concentration-BVP}\label{subsec:bvp}
%%%%%%%%%%%%%%%%%%%%%%%%%%%%%%%%%%%%%%%%%%%%%%%%%%%%%%%%%

We first introduce functions $\Bh_{\GO,j}$ for the boundary value problem, analogously to the functions $\Bh_j$ defined in \eqnref{hj_def} for the free space problem. They are solutions to the following problem:
\beq \label{eqn_hOmega}
 \ \left \{
 \begin{array} {ll}
\ds \Lcal_{\Gl,\Gm} \Bh_{\GO,j}= 0 \quad &\mbox{ in } \widetilde{\GO},\\[2mm]
\ds \Bh_{\GO,j}=\frac{(-1)^i}{2}\Psi_j \quad &\mbox{ on } \p D_i, \ i=1,2,
\\[2mm]
\ds \Bh_{\GO,j}=0 \quad &\mbox{ on } \p \GO.
 \end{array}
 \right.
 \eeq
One can easily see that the solution $\Bu$ to \eqnref{elas_eqn_bdd} admits the decomposition
\beq\label{eqn_decomp_u_bdd}
\Bu =  \Bv_{\GO} - \sum_{j=1}^3 (c_{1j}-c_{2j}){\Bh}_{\GO,j}  \quad \mbox{in }\widetilde\GO,
\eeq
where $\Bv_\GO$ is the solution to $\Lcal_{\Gl, \Gm} \Bv_\GO =0$ in $\widetilde{\GO}$ with the boundary condition
\beq
\Bv_\GO = \frac{1}{2} \sum_{j=1}^3 (c_{1j}+c_{2j}) \Psi_j \quad \mbox{on } \p D_1 \cup \p D_2.
\eeq
Note that
$$
\Bv_\GO|_{\p D_1} - \Bv_\GO|_{\p D_2} = \frac{1}{2} (c_{13}+c_{23}) (\Psi_3 |_{\p D_1} - \Psi_3 |_{\p D_2})= O(|\Bx|),
$$
from which one expects that $\nabla \Bv_\GO$ does not blow up even when $\Ge \to 0$. In fact, it was proved in \cite{BLL-ARMA-15} that
\beq\label{eqn_v_Omega_grad_bdd}
 \|\nabla\Bv_{\GO}\|_{L^\infty(\widetilde\GO)} \lesssim \| \Bg\|_{C^{1,\Gg}(\p\GO)}.
\eeq
So the singular behavior of $\nabla\Bu$ is determined by the function $\sum_{j=1}^3 (c_{1j}-c_{2j}){\Bh}_{\GO,j}$.

In the sequel, we investigate asymptotic behavior of $c_{1j}-c_{2j}$ and $\Bh_{\GO,j}$ as $\Ge \to 0$.
For doing so, we introduce the following boundary integrals:
\beq\label{def_Ical_Jcal}
\Ical_{jk}:=\int_{D^e} \Cbb \hatna \Bh_j :\hatna \Bh_k  \quad\mbox{ and }\quad \Jcal_{\GO,k}:=\int_{\p D^e} \frac{\p\Bh_k}{\p \Gv}\Big|_+ \cdot \BH_\GO, \quad j,k=1,2,3,
\eeq
where $\Bh_j$ is the solution to \eqnref{hj_def} in $\Acal^*$ and $\BH_\GO$ is the function defined by \eqnref{eqn_def_H_Omega}. We emphasize that $\Ical_{jk}$ is defined by $\Bh_j$, not by $\Bh_{\GO,j}$.

The relation among $c_{1j}-c_{2j}$, $\Ical_{jk}$ and $\Jcal_{\GO,k}$ is given by the following lemma.

\begin{lemma}\label{lem_represent_c1_c2_bdd}
The constants $c_{ij}$ appearing in \eqnref{elas_eqn_bdd} satisfy
\beq\label{c_diff_rep}
\begin{bmatrix}
\Ical_{11} && \Ical_{12} && \Ical_{13} \\
\Ical_{12} && \Ical_{22} && \Ical_{23} \\
\Ical_{13} && \Ical_{23} && \Ical_{33}
\end{bmatrix}
\begin{bmatrix}
c_{11}-c_{21} \\
c_{12}-c_{22} \\
c_{13}-c_{23}
\end{bmatrix}
=
\begin{bmatrix}
\Jcal_{\GO,1} \\
\Jcal_{\GO,2} \\
\Jcal_{\GO,3}
\end{bmatrix}.
\eeq
\end{lemma}

By inverting \eqnref{c_diff_rep}, we will see that the asymptotic behavior of $c_{1j}-c_{2j}$ as $\Ge \to 0$ can be described in terms of $\Kcal_{\GO,j}$ which are defined by
\beq\label{def_Kcal_Omega}
{\Kcal}_{\GO,1}={\Jcal}_{\GO,1}-\frac{{\Jcal}_{\GO,3} \Ical_{13}}{\Ical_{33}}, \quad
{\Kcal}_{\GO,2}={\Jcal}_{\GO,2}-\frac{{\Jcal}_{\GO,3} \Ical_{23}}{\Ical_{33}},\quad
{\Kcal}_{\GO,3}=\frac{{\Jcal}_{\GO,3}}{\Ical_{33}}.
\eeq
In fact, the following propositions hold. Here, we mention that they are consequences of Proposition \ref{prop_Brj}, which is proved by the variational principle and the properties of singular functions $\Bq_j$.
\begin{prop}\label{prop_Kcal_Omega_estim}
For $j=1,2,3$, we have
\beq\label{Kcaljest}
|\Kcal_{\GO,j}|\lesssim \| \Bg\|_{C^{1,\gamma}(\p\GO)}.
\eeq
\end{prop}

\begin{prop}\label{prop_diff_c1_c2_asymp_bdd}
We have
\begin{align}
c_{11}-c_{21}&=  \Kcal_{\GO,1}m_1^{-1} \sqrt\Ge +O(\sqrt\Ge\widetilde{E}), \label{c11c21}
\\
c_{12}-c_{22}&= \Kcal_{\GO,2}m_2^{-1} \sqrt\Ge +O(\sqrt\Ge\widetilde{E}), \label{c12c22}
\\
c_{13}-c_{23}&= \Kcal_{\GO,3} + O(\widetilde{E}), \label{c13c23}
\end{align}
where
\beq
\widetilde{E}:=(\sqrt\Ge+\tau \sqrt\Ge|\ln\Ge|)\| \Bg\|_{C^{1,\gamma}(\p\GO)}.
\eeq
\end{prop}

As an immediate consequence of Propositions \ref{prop_Kcal_Omega_estim} and \ref{prop_diff_c1_c2_asymp_bdd}, we obtain the following corollary.

\begin{cor}\label{cor_cij_estim}
We have
\beq
|c_{11}-c_{21}|+ |c_{12}-c_{22}| \lesssim \sqrt\Ge \| \Bg\|_{C^{1,\gamma}(\p\GO)}
\eeq
and
\beq
|c_{13}-c_{23}| \lesssim  \| \Bg\|_{C^{1,\gamma}(\p\GO)}.
\eeq
\end{cor}

Regarding the asymptotic behavior of $\Bh_{\GO,j}$, we obtain the following proposition.

\begin{prop}\label{hGOjandBq}
Let $m_j$, $j=1,2$, be the constant defined by \eqnref{m2def}. We have for $j=1,2$
\beq\label{eqn_hOj_qj_rOj_decomp}
\Bh_{\GO,j} = \frac{m_j}{\sqrt\Ge} \Bq_j +\Br_{\GO,j},
\eeq
where $\Br_{\GO,j}$ satisfies
\beq\label{rGOj}
\begin{cases}
\ds|\nabla \Br_{\GO,j}(\Bx) | \lesssim 1+\frac{\tau |y|}{\Ge+y^2} &\quad \mbox{for }\Bx\in \Pi_{L_0},
\\[1em]
\ds|\nabla \Br_{\GO,j}(\Bx) | \lesssim 1 &\quad \mbox{for }\Bx\in \widetilde\GO\setminus{\Pi_{L_0}}.
\end{cases}
\eeq
Here, $\tau$ is the constant defined by \eqnref{def_tau}. We also have
\beq\label{eqn_hOj_qj_rOj_decomp2}
\Bh_{\GO,3} =  \Bq_3+ \Br_{\GO,3} \quad\mbox{in } \widetilde\GO,
\eeq
where $\Br_{\GO,3}$ satisfies
\beq\label{rGO3}
| \nabla \Br_{\GO,3}| \lesssim 1 \quad \mbox{in } \widetilde\GO.
\eeq
\end{prop}

It is worth emphasizing that if two inclusions are symmetric with respect to both $x$- and $y$-axes, then $\tau=0$. So we have $|\nabla \Br_{\GO,j}|\lesssim 1$ in $\widetilde \GO$ for $j=1,2$ as well.

With help of preliminary results presented above, we are now able to state and prove the main results of this section.

\begin{theorem}\label{main_thm_2_general_bdd}
Let $\Bu$ be the solution to \eqnref{elas_eqn_bdd} for some $\Bg \in C^{1,\Gg}(\p\GO)$.
The following decomposition holds
\beq
\Bu(\Bx) = \Bb_\GO(\Bx)- \sum_{j=1}^3 \big({\Kcal}_{\GO,j}+s_{\GO,j}\big)\Bq_j(\Bx)
 , \quad\Bx\in \widetilde\GO,
\eeq
where ${\Kcal}_{\GO,j}$ are the constants defined by \eqnref{def_Kcal_Omega} (so satisfies \eqnref{Kcaljest}), $s_{\GO, j}$ are constants satisfying
\beq
|s_{\GO,j}| \lesssim \tau\sqrt\Ge|\ln\Ge| \|\Bg\|_{C^{1,\Gg}(\p\GO)},
\eeq
and the function $\Bb_\GO$ satisfies
\beq
\|\nabla \Bb_\GO\|_{L^\infty(\widetilde\GO)} \lesssim \|\Bg\|_{C^{1,\Gg}(\p\GO)}.
\eeq
\end{theorem}

\begin{theorem}\label{cor-bdd}
It holds that
\beq\label{gradestbdd}
 \frac{ \sum_{j=1}^2 |\Kcal_{\GO,j}| }{\sqrt{\Ge}}\lesssim\| \nabla \Bu\|_{L^\infty(\widetilde\GO)} \lesssim \frac{\|\Bg\|_{C^{1,\Gg}(\p\GO)}}{\sqrt{\Ge}}.
\eeq
\end{theorem}

The upper estimate in \eqnref{gradestbdd} was proved in \cite{BLL-ARMA-15}. The lower estimate shows that $\Ge^{-1/2}$ is also the lower bound on the blow-up rate of $\nabla \Bu$ as $\Ge \to 0$, provided that
\beq
1 \lesssim \sum_{j=1}^2 |\Kcal_{\GO,j}|.
\eeq
We will show in some special cases that this is the case (see section \ref{sec:symmetric_case}).

\medskip
\noindent{\sl Proof of Theorem \ref{main_thm_2_general_bdd}}.
According to Proposition \ref{prop_diff_c1_c2_asymp_bdd}, $c_{1j}-c_{2j}$ can be written as
\begin{align*}
c_{11}-c_{21}&=  \Kcal_{\GO,1}m_1^{-1} \sqrt\Ge + m_1^{-1}\sqrt\Ge (s_{\GO, 1} +s_{\GO, 1}'),
\\
c_{12}-c_{22}&= \Kcal_{\GO,2}m_2^{-1} \sqrt\Ge + m_2^{-1}\sqrt\Ge (s_{\GO, 2} +s_{\GO, 2}'),
\\
c_{13}-c_{23}&= \Kcal_{\GO,3} + s_{\GO, 3} +s_{\GO, 3}',
\end{align*}
where the constants $s_{\GO,j}$ and $s_{\GO,j}'$ satisfy
\begin{align}
|s_{\GO,j}|&\lesssim \tau \sqrt\Ge |\ln\Ge|\| \Bg\|_{C^{1,\gamma}(\p\GO)}, \label{s_GO_j_estim_temp1}
\\
 |s_{\GO,j}'|&\lesssim \sqrt\Ge\| \Bg\|_{C^{1,\gamma}(\p\GO)}. \label{s_GO_j_estim_temp2}
\end{align}
By substituting \eqnref{eqn_hOj_qj_rOj_decomp} and above three identities into \eqnref{eqn_decomp_u_bdd}, we have
\begin{align}
 \Bu
 &= \Bv_\GO - \bigg(\sum_{j=1}^2 (c_{1j}-c_{2j} )\Big(\frac{m_j}{\sqrt\Ge}\Bq_j + \Br_{\GO,j}\Big)\bigg)
 -(c_{13}-c_{23}) (\Bq_3 + \Br_{\GO,3})\nonumber
 \\
 &= \Bv_\GO - \sum_{j=1}^3 (\Kcal_{\GO,j} + s_{\GO,j} + s_{\GO,j}' )\Bq_j - \sum_{j=1}^3(c_{1j}-c_{2j}) \Br_{\GO,j}.
 \label{eqn_u_decomp_temp1}
\end{align}

Let
$$
\Bb_\GO := \Bu + \sum_{j=1}^3 (\Kcal_{\GO,j} + s_{\GO,j})\Bq_j.
$$
Then, from \eqnref{eqn_u_decomp_temp1}, we have
\begin{align*}
\nabla \Bb_\GO &=  \nabla \Bv_\GO -\sum_{j=1}^3 s_{\GO,j}'\nabla\Bq_j -\sum_{j=1}^3 (c_{1j}-c_{2j}) \nabla \Br_{\GO,j}
=: I_1+I_2+I_3.
\end{align*}

We now prove that $I_j$ are bounded. That $|I_1| \lesssim \| \Bg\|_{C^{1,\gamma}(\p\GO)}$ is already mentioned in \eqnref{eqn_v_Omega_grad_bdd}.
By \eqnref{Bq12est} and \eqnref{Bq3est}, we have
\beq\label{Bq123est}
\| \nabla \Bq_j\|_{L^\infty(\widetilde\GO)} \lesssim \Ge^{-1/2}, \quad j=1,2,3.
\eeq
So, by \eqnref{s_GO_j_estim_temp2},  we have
$$
|I_2| \lesssim  \| \Bg\|_{C^{1,\gamma}(\p\GO)}.
$$
We have from \eqnref{rGOj} and \eqnref{rGO3} that
$$
\| \nabla \Br_{\GO,1}\|_{L^\infty(\widetilde\GO)}+\| \nabla \Br_{\GO,2}\|_{L^\infty(\widetilde\GO)} \lesssim 1+\tau /\sqrt\Ge
$$
and
$$
\| \nabla \Br_{\GO,3}\|_{L^\infty(\widetilde\GO)} \lesssim 1.
$$
Therefore, it follows from Corollary \ref{cor_cij_estim} that
\begin{align*}
|I_3| &\leq \left| \sum_{j=1}^2 (c_{1j}-c_{2j}) \nabla\Br_{\GO,j} \right| + \left| (c_{13}-c_{23}) \nabla \Br_{\GO,3} \right|
\\
&\lesssim
\big(\sqrt\Ge (1+\tau /\sqrt\Ge) + 1\big)\| \Bg\|_{C^{1,\gamma}(\p\GO)}\lesssim \| \Bg\|_{C^{1,\gamma}(\p\GO)}.
\end{align*}
The proof is complete.
\qed

\medskip

\noindent{\sl Proof of Theorem \ref{cor-bdd}}.
The upper estimate in \eqnref{gradestbdd} is a consequence of Proposition \ref{prop_Kcal_Omega_estim}, Theorem \ref{main_thm_2_general_bdd}, and \eqnref{Bq123est}.
To derive the lower estimate, we consider $\nabla \Bu(0,0)$.
It follows from Lemma \ref{singular_q_origin}, Lemma \ref{cor_q3_near_origin} and Theorem \ref{main_thm_2_general_bdd} that
$$
\nabla \Bu(0,0) =-\frac{\mathcal{K}_{\GO,1}}{m_1\sqrt\Ge}\Be_1\otimes\Be_1
- \frac{\mathcal{K}_{\GO,2}}{m_2\sqrt\Ge}\Be_2\otimes\Be_1 + O(1+\tau \ln\Ge).
$$
So we obtain the lower estimate.
\qed

\medskip

As mentioned earlier, the upper bound in \eqnref{gradestbdd} was proved in \cite{BLL-ARMA-15}. So, it is helpful to compare the method of that paper with that of this paper. In fact, some of results obtained in \cite{BLL-ARMA-15} will be used for proofs in this section. There, the solution $\Bu$ to \eqnref{elas_eqn_bdd} is expressed as follows:
\beq\label{yyli_decomp}
\Bu=\sum_{i=1}^2\sum_{k=1}^3 c_{ik} \Bv_{ik} + \Bv_3,
\eeq
where $\Bv_{ik}$ is the solution to
\beq\label{vikyan}
\begin{cases}
\Lcal_{\Gl,\Gm}\Bv_{ik}=0 &\quad \mbox{in } \widetilde \GO,
\\
\Bv_{ik} = \Psi_k &\quad \mbox{on }\p D_i,
\\
\Bv_{ik} = 0 &\quad \mbox{on } \p D_j \cup \p\GO, \ j\neq i,
\end{cases}
\eeq
and $\Bv_3$ is the solution to
$$
\begin{cases}
\Lcal_{\Gl,\Gm}\Bv_3=0 &\quad \mbox{in } \widetilde \GO,
\\
\Bv_3 = 0 &\quad \mbox{on }\p D_1\cup\p D_2,
\\
\Bv_3 = \Bg &\quad \mbox{on } \p\GO.
\end{cases}
$$
Note that
\beq\label{eqn_hGOk_v1k_v2k}
\Bh_{\GO,k}=-\frac{1}{2}\Bv_{1k}+\frac{1}{2}\Bv_{2k}, \quad k=1,2,3.
\eeq

The $6\times 6$ linear system of equations for $c_{ik}$ is derived using \eqnref{int_zero}. The linear system is truncated to a $3\times 3$ one and then the difference $c_{1j}-c_{2j}$ is expressed using the following integrals:
\begin{align*}
a_{jk} &:= \int_{\p D_1} \p_\nu \Bv_{1j}|_+ \cdot \Psi_k = \int_{\widetilde\GO} \Cbb \hatna \Bv_{1j}: \hatna \Bv_{1k},
\\
b_k &:= \int_{\p D_1} \p_\nu \Bv_3|_+ \cdot \Psi_k,
\end{align*}
for $j,k=1,2,3$. Note that the integral $a_{jk}$ is similar to the quantity $\Ical_{jk}$ of this paper. The difference lies in that $\Ical_{jk}$ is defined using the free space solution $\Bh_j$.

To investigate asymptotic behavior of $a_{jk}$ and $b_k$ as $\Ge\rightarrow 0$,
the function $\Bv_{ik}$ is approximated by $\Bv_{ik}^K$, which is defined by
\beq\label{BvikK}
\begin{cases}
\ds \Bv_{1k}^K(x,y)= \frac{-x+f_2(y)}{f_1(y)+f_2(y)} \Psi_k, \\
\nm
\ds \Bv_{2k}^K(x,y)= \frac{x+f_1(y)}{f_1(y)+f_2(y)} \Psi_k,
\end{cases}
(x,y) \in \Pi_L, \quad k=1,2,3.
\eeq
In fact, it is proved that
\begin{align}
\nabla \Bv_{ik}(x,y) &= \nabla \Bv_{ik}^K(x,y) + O \left( 1+ \frac{y}{\Ge+y^2} \right) \quad\mbox{for } k=1,2, \label{Bao1} \\
\nabla \Bv_{i3}(x,y) &= \nabla \Bv_{i3}^K(x,y) + O (1).  \label{Bao2}
\end{align}
From these approximations which are derived using a new iteration technique, the upper bound on the blow-up rate, $\Ge^{-1/2}$, of $|\nabla \Bu|$ is obtained in \cite{BLL-ARMA-15}. However, a lower bound has not been obtained. It is partly because the functions $\Bv_{ik}^K$ are {\it not} solutions of the Lam\'e system.

In this paper, we introduce new singular functions $\Bq_j$, which are solutions of the Lam\'{e} system, as explained in section \ref{sec:singular}. Using singular functions, we are able to derive precise asymptotic formulas for $\nabla \Bu$ as $\Ge \to 0$. As a consequence we are able to reprove that $\Ge^{-1/2}$ is indeed an upper bound on the blow-up rate. Moreover, the asymptotic formulas enable to show that $\Ge^{-1/2}$ is a lower bound on the blow-up rate as well in some cases, as presented in section \ref{sec:symmetric_case}. We emphasize that the asymptotic formulas are obtained using the variational principle, which is possible only because $\Bq_j$ are solutions of the Lam\'{e} system.

%%%%%%%%%%%%%%%%%%%%%%%%%%%%%%%%%%%%%%%%%%%%%%%%%%%%%%%%%%%%%
\subsection{Preliminary estimates of boundary integrals}
%%%%%%%%%%%%%%%%%%%%%%%%%%%%%%%%%%%%%%%%%%%%%%%%%%%%%%%%%%%%%

In this subsection, we characterize asymptotic behaviors of the following boundary integrals as $\Ge \to 0$:
$$
\int_{\p D_i} {\p_\nu \Bq_j}\cdot \Psi_k,
\quad
\int_{\p D^e} {\p_\nu \Bq_j}\cdot \Bq_k.
$$
These integrals appear in later sections.

We first prove the following lemma.

\begin{lemma}\label{lem_hS_bdry_int12}
\begin{itemize}
\item[(i)] For $k=1,2$, we have
\beq\label{qjpsik}
\ds\int_{\p D_i} {\p_\nu \Bq_j}\cdot \Psi_k =(-1)^{i+1}\Gd_{jk} , \quad i, j=1,2.
\eeq
\item[(ii)] For $k=3$, we have
\begin{align}
\ds\int_{\p D_i} {\p_\nu \Bq_1}\cdot \Psi_3 &= 0,
\label{int_h1S_bdry_Psi3}
\\
\ds\int_{\p D_i} {\p_\nu \Bq_2}\cdot \Psi_3 &= (-1)^{i+1} {a} (-1 + 4\pi \Ga_2 \Gm),
\label{int_h2S_bdry_Psi3}
\end{align}
for $i=1,2$, where $a$ is the constant defined by \eqnref{a_def}.
\end{itemize}
\end{lemma}

\pf
Suppose that $k=1,2$.
Since $\Lcal_{\Gl,\Gm}\BGG(\Bx-\Bp_l)\Be_j = \Gd(\Bx-\Bp_l) \Be_j$, Green's formula yields
\begin{align*}
\int_{\p D_i} \p_{\nu_\Bx}\BGG(\Bx-\Bp_l) \Be_j\cdot \Psi_k  \,d\Gs(\Bx) = \int_{D_i} \Lcal_{\Gl,\Gm}\BGG(\Bx-\Bp_l)\Be_j \cdot \Psi_k  \,d\Gs(\Bx) = \Gd_{il}\Gd_{jk}.
\end{align*}
Green's formula also yields
\begin{align*}
\int_{\p D_i} \p_{\nu_\Bx}\left(\frac{\Bx-\Bp_l}{|\Bx-\Bp_l|^2}\right)\cdot \Psi_k  \,d\Gs(\Bx)=0.
\end{align*}
In fact, if $i=l$, then we apply Green's formula to $\Rbb^2 \setminus D_i$, and to $D_i$ if $i \neq l$. So, \eqnref{qjpsik} follows from \eqnref{Bqone} and \eqnref{Bqtwo}.

We now prove \eqnref{int_h1S_bdry_Psi3} and \eqnref{int_h2S_bdry_Psi3} when $i=1$. The case when $i=2$ can be proved in the same way.
Let us prove \eqnref{int_h2S_bdry_Psi3} first.
In view of the definition \eqnref{Bqtwo} of $\Bq_2$, we have
\beq\label{350}
\int_{\p D_1} {\p_\nu \Bq_2}\cdot \Psi_3 =  \int_{\p D_1} \p_{\nu_\Bx} (\BGG (\Bx-\Bp_1)\Be_2) \cdot \Psi_3 -
{\Ga_2 a} \int_{\p D_1} \p_{\nu_\Bx}\left( \frac{(\Bx-\Bp_1)^\perp}{|\Bx-\Bp_1|^2}\right) \cdot \Psi_3.
\eeq
Since $\Lcal_{\Gl,\Gm}(\BGG (\Bx-\Bp_1) \Be_2)= \Gd_{\Bp_1}(\Bx) \Be_2$, one can see that
\beq\label{351}
\int_{\p D_1} \p_{\nu_\Bx} (\BGG (\Bx-\Bp_1)\Be_2) \cdot \Psi_3 = \Be_2 \cdot \Psi_3(\Bp_1)= -a,
\eeq
where the last equality holds because $\Bp_1=(-a,0)$.

By using a change of variables $\Bx\rightarrow \Bx+\Bp_1$ and the fact that $\Psi_3(\Bx+\Bp_1)=\Psi_3(\Bx) - a\Psi_2$, we obtain
\begin{align}
& \int_{\p D_1} \p_{\nu_\Bx}\left( \frac{(\Bx-\Bp_1)^\perp}{|\Bx-\Bp_1|^2}\right) \cdot \Psi_3(\Bx)
=\int_{\p D_1-\Bp_1} \p_{\nu_\Bx}\left( \frac{\Bx^\perp}{|\Bx|^2}\right) \cdot \Psi_3(\Bx+\Bp_1) \nonumber \\
&= -a \int_{\p D_1-\Bp_1} \p_{\nu_\Bx}\left( \frac{\Bx^\perp}{|\Bx|^2}\right) \cdot \Psi_2
+\int_{\p D_1-\Bp_1} \p_{\nu_\Bx}\left( \frac{\Bx^\perp}{|\Bx|^2}\right) \cdot \Psi_3(\Bx). \label{int_dipole_Psi3_decomp}
\end{align}
One can show as before that
\beq\label{int_dipole_Psi3_prep}
\int_{\p D_1-\Bp_1} \p_{\nu_\Bx}\left( \frac{\Bx^\perp}{|\Bx|^2}\right) \cdot \Psi_2 \, d\Gs(\Bx) =0.
\eeq
Let $B$ be a disk centered at $0$ such that $\p D_1-\Bp_1 \subset B$. Then Green's formula yields
$$
\int_{\p D_1-\Bp_1} \p_{\nu_\Bx}\left( \frac{\Bx^\perp}{|\Bx|^2}\right) \cdot \Psi_3(\Bx) =
\int_{\p B} \p_{\nu_\Bx}\left( \frac{\Bx^\perp}{|\Bx|^2}\right) \cdot \Psi_3(\Bx).
$$
Straightforward computations show that
$$
\p_{\nu_\Bx}\left( \frac{\Bx^\perp}{|\Bx|^2}\right) = (-2\mu) \frac{\Bx^\perp}{|\Bx|^3} \quad \mbox{ for }\Bx \in \p B.
$$
So, we have
$$
\int_{\p D_1-\Bp_1} \p_{\nu_\Bx}\left( \frac{\Bx^\perp}{|\Bx|^2}\right) \cdot \Psi_3(\Bx) =
\int_{\p B} \p_{\nu_\Bx}\left( \frac{\Bx^\perp}{|\Bx|^2}\right) \cdot \Bx^{\perp}
= \int_{\p B} (-2\mu) \frac{1}{|\Bx|}=-4\pi \Gm .
$$
It then follows from \eqnref{int_dipole_Psi3_decomp} and \eqnref{int_dipole_Psi3_prep} that
\beq\label{352}
\int_{\p D_1} \p_{\nu_\Bx} \left( \frac{(\Bx-\Bp_1)^\perp}{|\Bx-\Bp_1|^2}\right) \cdot \Psi_3(\Bx) =-4\pi \Gm.
\eeq
Combining \eqnref{350}-\eqnref{352}, we obtain \eqnref{int_h2S_bdry_Psi3}.

We now prove \eqnref{int_h1S_bdry_Psi3}. Like \eqnref{351} we have
$$
\int_{\p D_1} \p_{\nu_\Bx} (\BGG (\Bx-\Bp_1)\Be_1) \cdot \Psi_3 = \Be_1 \cdot \Psi_3(\Bp_1)= 0.
$$
In the same way to show \eqnref{352} one can show that
$$
\int_{\p D_1} \p_{\nu_\Bx}\left( \frac{\Bx-\Bp_1}{|\Bx-\Bp_1|^2}\right) \cdot \Psi_3=0.
$$
Therefore, from the definition \eqnref{Bqone} of $\Bq_1$, we have \eqnref{int_h1S_bdry_Psi3}, and the proof is completed.
\qed

\begin{lemma}\label{lem_dhS_hS_int12}
We have
\begin{align}
&\ds\int_{\p D^e} {\p_{\nu} \Bq_j}
\cdot \Bq_j
= - {m_j^{-1}}\sqrt\Ge + O(\tau\Ge|\ln\Ge|+\Ge), \quad j=1,2,
\label{qjqj}
\\
& \ds\int_{\p D^e} {\p_\nu \Bq_1} \cdot \Bq_2 = O(\tau\Ge|\ln\Ge|+\Ge).
\label{q1q2}
\end{align}
\end{lemma}

Before proving  Lemma \ref{lem_dhS_hS_int12}, we need to estimate the conormal derivatives $\p_{\nu} \Bq_j$ on $\p D_1 \cup \p D_2$. We have the following lemma.
\begin{lemma} \label{lem_conormal_estim}
For $\Bx=(x,y)\in (\p D_1 \cup \p D_2) \cap \p\Pi_{L_0}$, we have
\beq\label{4200}
\big|{\p_\nu \Bq_1}(\Bx)\cdot \Be_1\big| \lesssim \frac{\sqrt\Ge}{\Ge+y^2},
\quad
\big|{\p_\nu \Bq_1}(\Bx)\cdot \Be_2\big| \lesssim \frac{\sqrt\Ge|y|}{\Ge+y^2}+\sqrt\Ge,
\eeq
and
\beq\label{4300}
\big|{\p_\nu \Bq_2}(\Bx)\cdot \Be_1\big| \lesssim \frac{\sqrt\Ge|y|}{\Ge+y^2}+\sqrt\Ge,
\quad
\big|{\p_\nu \Bq_2}(\Bx)\cdot \Be_2\big| \lesssim \frac{\sqrt\Ge}{\Ge+y^2}.
\eeq
\end{lemma}

\pf We prove \eqnref{4200} only. \eqnref{4300} can be proved similarly.
Let $\BGs^1= (\Gs^1_{ij})_{i,j=1}^2$ be the stress tensor of $\Bq_1$, namely, $\BGs^1 := \Cbb \hatna\Bq_1$. According to \eqnref{stress_cartesian}, the entries of $\BGs^1$ can be written as
\begin{align*}
\Gs^1_{11} &= (\Gl+2\Gm) \p_1 q_{11} + \Gl\p_2 q_{12}, \\
\Gs^1_{22} &= \Gl \p_1 q_{11} + (\Gl+2\Gm) \p_2 q_{12}, \\
\Gs^1_{12} &= \Gs^k_{21} = \Gm (\p_2 q_{11} + \p_1 q_{12}).
\end{align*}
Thus we have the following estimates from Lemma \ref{lem_hS_grad_estim}:
$$
|\Gs^1_{11}| + |\Gs^1_{22}|\lesssim \frac{\sqrt\Ge}{\Ge+y^2},
\qquad
|\Gs^1_{12}| \lesssim \frac{\sqrt\Ge|y|}{\Ge+y^2} + \sqrt\Ge \qquad \mbox{for } (x,y)\in\Pi_{L_0}.
$$

Note that $\p_\nu \Bq_1=\BGs^1 \Bn$ and the outward unit normal vector $\Bn$ on $\p D_i \cap \p\Pi_{L_0}$ is given as follows:
$$
\Bn = \frac{1}{\sqrt{1+(f_i'(y))^2}} \big( (-1)^{i+1},f_i'(y) \big), \quad i=1,2.
$$
Moreover, we have $|f_i'(y)|\lesssim |y|$. Therefore, we obtain
\begin{align*}
|{\p_\nu \Bq_1}(x,y)\cdot \Be_1|
&=|(\BGs^1 \Bn )_1|
= \left| \frac{1}{\sqrt{1+(f_i'(y))^2}} ((-1)^{i+1}\Gs^1_{11}  + f_i'(y)\Gs^1_{12}) \right|
\\
&\lesssim \frac{\sqrt\Ge}{\Ge+y^2} + \frac{\sqrt\Ge y^2}{\Ge+y^2} \lesssim \frac{\sqrt\Ge}{\Ge+y^2},
\end{align*}
and
\begin{align*}
|{\p_\nu \Bq_1}(x,y)\cdot \Be_2|
&=
|(\BGs^1 \Bn )_2|
= \left| \frac{1}{\sqrt{1+(f_i'(y))^2}} ((-1)^{i+1}\Gs^1_{12}  + f_i'(y)\Gs^1_{22}) \right|
\\
&\lesssim   \frac{\sqrt\Ge |y|}{\Ge+y^2} + \sqrt\Ge,
\end{align*}
for $\Bx=(x,y)\in \p D_i \cap \p \Pi_{L_0}$ and $i=1,2$.
The proof is completed.
\qed

\medskip

\noindent{\sl Proof of Lemma \ref{lem_dhS_hS_int12}}. To prove \eqnref{qjqj}, we write
\begin{align*}
\int_{\p D_1\cup \p D_2} {\p_\nu \Bq_1}\cdot \Bq_1
&= \sum_{i=1}^2 (-1)^i  (\Ga_1-\Ga_2)\Gk_i a \int_{\p D_i} {\p_\nu \Bq_1}\cdot\Psi_1 \\
&\qquad + \sum_{i=1}^2 \int_{\p D_i} {\p_\nu \Bq_1}\cdot \big[\Bq_1 -
(-1)^i (\Ga_1-\Ga_2) \Gk_i a \, \Psi_1 \big].
\end{align*}
By Lemma \ref{lem_hS_bdry_int12} (i), we have
$$
\sum_{i=1}^2 (-1)^i  (\Ga_1-\Ga_2)\Gk_i a \int_{\p D_i} {\p_\nu \Bq_1}\cdot\Psi_1 = -(\Ga_1-\Ga_2)(\Gk_1+\Gk_2) a.
$$
Then \eqnref{def_alpha} and \eqnref{m2def} yield
$$
\sum_{i=1}^2 (-1)^i  (\Ga_1-\Ga_2)\Gk_i a \int_{\p D_i} {\p_\nu \Bq_1}\cdot\Psi_1 = -m_1^{-1} \sqrt{\Ge} + O(\Ge^{3/2}).
$$
It then remains to show that
\beq\label{500}
\int_{\p D_i} {\p_\nu \Bq_1}\cdot \big[\Bq_1 -
(-1)^i (\Ga_1-\Ga_2) \Gk_i a \, \Psi_1 \big] = O(\tau\Ge|\ln\Ge|+\Ge), \quad i=1,2.
\eeq

To prove \eqnref{500}, let us write
\begin{align*}
\int_{\p D_1} {\p_\nu \Bq_1}\cdot \big[\Bq_1 -
(-1)  (\Ga_1-\Ga_2) \Gk_1 a \, \Psi_1
\big] =
\int_{\p D_1 \cap \p \Pi_{L_0} }
+ \int_{\p D_1  \setminus \p\Pi_{L_0} }
:=I_1 + I_2.
\end{align*}
From Lemma \ref{lem_qj_far_estim} and the fact that $a\approx\sqrt{\Ge}$, we see that $|I_2|\lesssim \Ge$.
Note that
$$
|I_1| \le \int_{\p D_1\cap\p \Pi_{L_0}} \big|{\p_\nu \Bq_1}\cdot\Be_1 \big(q_{11} +  (\Ga_1-\Ga_2) \Gk_1 a
\big)\big| + \big|{\p_\nu \Bq_1}\cdot\Be_2 \,q_{12}\big|.
$$
From \eqnref{q11asym} and \eqnref{q12asym}, we see that
$$
| q_{11} +  (\Ga_1-\Ga_2) \Gk_1 a | \lesssim \Ge^{3/2} + \sqrt\Ge y^2+ \tau \sqrt\Ge |y|
$$
and
$$
|q_{12}|\lesssim \sqrt\Ge |y|.
$$
It then follows from Lemma \ref{lem_conormal_estim} that
\begin{align*}
|I_1| &\lesssim \int_{\p D_1\cap\p \Pi_{L_0}} \frac{\sqrt\Ge}{\Ge+y^2} (\Ge^{3/2} + \sqrt\Ge y^2+ \tau \sqrt\Ge |y|) + \Big(\frac{ \sqrt\Ge |y|}{\Ge+y^2} + \sqrt\Ge \Big)  \sqrt\Ge|y|
\\
&\lesssim \int_{-L_0}^{L_0} \frac{\Ge \tau|y|}{\Ge+y^2} \,dy + \Ge
\lesssim \tau \Ge |\ln \Ge|+ \Ge.
\end{align*}
This proves \eqnref{500} for $i=1$. The case for $i=2$ can be proved in the same way. So,
\eqnref{qjqj} is proved.

Next we prove \eqnref{q1q2}. Thanks to \eqnref{qjpsik} with $j=1$ and $k=2$, we can write
\begin{align}
\int_{\p D^e} {\p_\nu \Bq_1}\cdot \Bq_2
&=  \sum_{i=1}^2 \Ga_2 \Gk_i^2 a \int_{\p D_i} \p_\nu \Bq_1 \cdot y\Psi_1 \nonumber \\
&\qquad + \sum_{i=1}^2 \int_{\p D_i} {\p_\nu \Bq_1}\cdot \big[\Bq_2 -
\Ga_2 \Gk_i^2 a y \,\Psi_1 -(-1)^i(\Ga_1+\Ga_2) \Gk_i a \, \Psi_2) \big]. \label{5000}
\end{align}
Green's formula yields
\begin{align*}
\int_{\p D_i}\p_\nu \Bq_1 \cdot y\Psi_1 &= \int_{\p D_i}\p_\nu \Bq_1 \cdot y\Psi_1
-\int_{\p D_i}\p_\nu (y\Psi_1) \cdot \Bq_1 +\int_{\p D_i}\p_\nu (y\Psi_1) \cdot \Bq_1 \\
&= \int_{ D_i}\Lcal_{\Gl,\Gm} \Bq_1 \cdot y\Psi_1
-\int_{ D_i} \Lcal_{\Gl,\Gm}(y\Psi_1) \cdot \Bq_1 +\int_{\p D_i}\p_\nu (y\Psi_1) \cdot \Bq_1 \\
&= \int_{ D_i}\Lcal_{\Gl,\Gm} \Bq_1 \cdot y\Psi_1 +\int_{\p D_i}\p_\nu (y\Psi_1) \cdot \Bq_1 .
\end{align*}
Observe from \eqnref{doublet} and the definition \eqnref{Bqone} of $\Bq_1$ that
\beq\label{qonedirac}
\Lcal_{\Gl,\Gm}\Bq_1 = ({\Gd_{\Bp_1}-\Gd_{\Bp_2})\Be_1} + \sum_{j=1}^2 \frac{\Ga_2 a}{\Ga_1-\Ga_2}
\big( \p_1 \Gd_{\Bp_j} \Be_1 + \p_2 \Gd_{\Bp_j} \Be_2 \big),
\eeq
where $\Gd_{\Bp_j}$ denotes the Dirac delta at $\Bp_j$. So, we see that
$$
\int_{ D_i}\Lcal_{\Gl,\Gm} \Bq_1 \cdot y\Psi_1 =0.
$$
It follows from Lemma \ref{lem_qj_far_estim} and Lemma \ref{lem_hS_bdry_estim} that $\| \Bq_1\|_{L^\infty(\p D_i)} \lesssim \sqrt{\Ge}$ for $i=1,2$. So we have
$$
\int_{\p D_i}\p_\nu (y\Psi_1) \cdot \Bq_1 = O(\sqrt\Ge),
$$
and hence
\beq\label{5100}
\int_{\p D_i}\p_\nu \Bq_1 \cdot y\Psi_1 = O(\sqrt\Ge), \quad i=1,2.
\eeq

Let
\begin{align*}
\int_{\p D_1} {\p_\nu \Bq_1}\cdot \big[\Bq_2 -
\Ga_2 \Gk_1^2 a y \,\Psi_1 +(\Ga_1+\Ga_2) \Gk_1 a \, \Psi_2)
\big] &=
\int_{\p D_1 \cap \p \Pi_{L_0} }
+
\int_{\p D_1  \setminus \p\Pi_{L_0} }
\\
&:=J_1+J_2.
\end{align*}
As before, from Lemma \ref{lem_qj_far_estim} and the fact that $a\approx\sqrt{\Ge}$, we see that $|J_2| \lesssim \Ge$.
From  Lemma \ref{lem_hS_bdry_estim}, Lemma \ref{lem_conormal_estim} and the fact that $a\approx\sqrt\Ge$, we have
\begin{align*}
|J_1|&\lesssim \int_{\p D_1\cap\p \Pi_{L_0}} \big|{\p_\nu \Bq_1}\cdot\Psi_1 \, \big( q_{21} - \Ga_2 \Gk_1^2 a y \big) \big| + \big|{\p_\nu \Bq_1}\cdot\Psi_2 \,\big( q_{22} +  (\Ga_1+\Ga_2) \Gk_1 a
\big)\big|
\\
&\lesssim \int_{\p D_1\cap\p \Pi_{L_0}} \Big( \frac{\sqrt\Ge}{\Ge+y^2}   + \frac{ \sqrt\Ge |y|}{\Ge+y^2} + \sqrt\Ge \Big)
(\Ge^{3/2} + \sqrt\Ge y^2+ \tau \sqrt\Ge |y|)
\\
&\lesssim \int_{-L_0}^{L_0} \frac{\tau\Ge |y|}{\Ge+y^2} \,dy + \Ge
\\
&\lesssim \tau \Ge |\ln \Ge|+ \Ge.
\end{align*}
So we obtain
\beq\label{5200}
\left| \int_{\p D_1} {\p_\nu \Bq_1}\cdot \big(\Bq_2  - \Ga_2 \Gk_1^2 a y \,\Psi_1 +  (\Ga_1+\Ga_2) \Gk_1 a \,\Psi_2 \big) \right| \lesssim\tau \Ge |\ln \Ge|+ \Ge.
\eeq
Similarly, one can see that
\beq\label{5300}
\left| \int_{\p D_2} {\p_\nu \Bq_1}\cdot \big(\Bq_2  - \Ga_2 \Gk_2^2 a y \,\Psi_1 -  (\Ga_1+\Ga_2) \Gk_2 a \,\Psi_2 \big) \right| \lesssim\tau \Ge |\ln \Ge|+ \Ge.
\eeq
Since $a\approx \sqrt\Ge$, \eqnref{q1q2} follows from \eqnref{5000} and \eqnref{5100}-\eqnref{5300}.
The proof is completed.
\qed

%%%%%%%%%%%%%%%%%%%%%%%%%%%%%%%%%%%%%%%%%%%%%%%%%%%%%%%%%%%%%%%%%%%
\subsection{Proof of Lemma \ref{lem_represent_c1_c2_bdd}}
%%%%%%%%%%%%%%%%%%%%%%%%%%%%%%%%%%%%%%%%%%%%%%%%%%%%%%%%%%%%%%%%%%%%

We first show that
\beq\label{Ijk_different_rep}
\Ical_{jk} =\int_{\p D_1} {\p_\nu \Bh_j} \cdot \Psi_k =- \int_{\p D_2} {\p_\nu \Bh_j} \cdot \Psi_k, \quad
j=1,2,\,k=1,2,3.
\eeq
In fact, we see from Lemma \ref{cor_betti} and the boundary conditions of $\Bh_j$ that
$$
\Ical_{jk} = -\int_{\p D^e} {\p_\nu \Bh_j} \cdot \Bh_k
= \frac{1}{2}\int_{\p D_1} {\p_\nu \Bh_j} \cdot \Psi_k -
\frac{1}{2}\int_{\p D_2} {\p_\nu \Bh_j} \cdot \Psi_k.
$$
Since $\hatna \Psi_k=0$, we obtain using Lemma \ref{cor_betti} that
\beq\label{eqn_int_hj_D1_D2_zero}
\int_{\p D^e} {\p_\nu \Bh_j} \cdot \Psi_k = -\int_{D^e} \Cbb \hatna \Bh_j :\hatna \Psi_k =0.
\eeq
So, \eqnref{Ijk_different_rep} follows.

Since $\Lcal_{\Gl,\Gm}\BH_\GO=0$ in $D_i$, we have
$$
\int_{\p D_i} {\p_\nu \BH_\GO}\cdot \Bh_j = (-1)^i \frac{1}{2} \int_{\p D_i} {\p_\nu\BH_\GO}\cdot \Psi_j=0.
$$
Thus we have
\beq\label{400}
\Jcal_{\GO,j} = \int_{\p D^e}   {\p_\nu  \Bh_j} \big|_+ \cdot \BH_\GO =
\int_{\p D^e}   {\p_\nu \Bh_j} \big|_+ \cdot \BH_\GO
-   {\p_\nu \BH_\GO}\cdot \Bh_j .
\eeq
One can easily see from \eqnref{rep_bdd} and Lemma \ref{lem:Acal} (i) that $\Bu-\BH_\GO$ can be extended to $D^e$ so that the extended function, still denoted by $\Bu-\BH_\GO$, satisfies $\Lcal_{\Gl,\Gm}(\Bu-\BH_\GO)=0$ in $D^e$ and $\Bu-\BH_\GO\in \Acal$. Therefore, we have from  Lemma \ref{cor_betti}
$$
\int_{\p D^e}
(\Bu-\BH_\GO)\cdot {\p_\nu \Bh_j} \big|_+ - {\p_\nu(\Bu-\BH_\GO)} \big|_+ \cdot\Bh_j
=0.
$$
We then infer from \eqnref{int_zero} and \eqnref{400} that
$$
\Jcal_{\GO,j} = \int_{\p D^e}  {\p_\nu  \Bh_j} \big|_+ \cdot \Bu - {\p_\nu  \Bu} \big|_+ \cdot \Bh_j = \int_{\p D^e}  {\p_\nu  \Bh_j} \big|_+ \cdot \Bu.
$$
Then the boundary condition in \eqnref{elas_eqn_bdd} and \eqnref{Ijk_different_rep} yield
\begin{align*}
\Jcal_{\GO,j} &= \sum_{k=1}^3 c_{1k} \int_{\p D_1} {\p_\nu  \Bh_j} \big|_+ \cdot  \Psi_k  +
c_{2k} \int_{\p D_2}  {\p_\nu  \Bh_j} \big|_+ \cdot \Psi_k \\
& = \sum_{k=1}^3 (c_{1k} - c_{2k}) \Ical_{jk}.
\end{align*}
So, \eqnref{c_diff_rep} follows.
\qed

%%%%%%%%%%%%%%%%%%%%%%%%%%%%%%%%%%%%%%%%%%%%%%%%%%%%%%%%%%%%%%
\subsection{Estimates of integrals $\Ical_{jk}$ and $\Jcal_k$ and proof of Proposition \ref{prop_Kcal_Omega_estim}}
%%%%%%%%%%%%%%%%%%%%%%%%%%%%%%%%%%%%%%%%%%%%%%%%%%%%%%%%%%%%%%

In this subsection we derive estimates of the integrals $\Ical_{jk}$ and $\Jcal_{\GO,k}$, and prove
Proposition \ref{prop_Kcal_Omega_estim} as a consequence. Some of estimates obtained in this subsection will be used in the next subsection as well.

\begin{lemma}\label{prop_I11_I22_I12_asymp}
The following holds:
\begin{align}
\Ical_{11} &= m_1 \Ge^{-1/2} + O(\tau|\ln\Ge| + 1), \label{Ical11} \\
\Ical_{12} &= O(\tau|\ln\Ge| + 1), \label{Ical12} \\
\Ical_{22} &= m_2 \Ge^{-1/2} + O(\tau|\ln\Ge| + 1), \label{Ical22}
\end{align}
as $\Ge \to 0$.
\end{lemma}

\pf
According to \eqnref{hj_def_decomp}, we have
\begin{align*}
\Ical_{jk} &= \int_{D^e} \Cbb \hatna \Bh_j :\hatna \Bh_k
=
\int_{D^e} \Cbb\hatna (\frac{m_j}{\sqrt\Ge}\Bq_j + \Br_j):\hatna \Bh_k
\\
&= \frac{m_j}{\sqrt\Ge}\int_{D^e} \Cbb\hatna \Bq_j:\hatna \Bh_k
+\int_{D^e} \Cbb\hatna \Br_j:\hatna \Bh_k
\\
&=\frac{m_j}{\sqrt\Ge}\int_{D^e} \Cbb\hatna \Bq_j:\hatna \Bh_k
+\int_{D^e} \Cbb\hatna \Br_j:\hatna (\frac{m_k}{\sqrt\Ge}\Bq_k + \Br_k)
\\
&=\frac{m_j}{\sqrt\Ge}\int_{D^e} \Cbb\hatna \Bq_j:\hatna \Bh_k
+\frac{m_k}{\sqrt\Ge}\int_{D^e} \Cbb\hatna \Br_j:\hatna \Bq_k
+ \int_{D^e} \Cbb\hatna \Br_j:\hatna\Br_k .
\end{align*}
Since
\begin{align*}
\frac{m_k}{\sqrt\Ge}\int_{D^e} \Cbb\hatna \Br_j:\hatna \Bq_k &=
\frac{m_k}{\sqrt\Ge}\int_{D^e} \Cbb\hatna (\Bh_j  - \frac{m_j}{\sqrt\Ge}\Bq_j):\hatna \Bq_k
\\
&=
\frac{m_k}{\sqrt\Ge}\int_{D^e} \Cbb\hatna \Bq_k:\hatna \Bh_j
-\frac{m_j m_k}{\Ge}\int_{D^e} \Cbb\hatna \Bq_j:\hatna \Bq_k,
\end{align*}
it follows that
\begin{align*}
\Ical_{jk}&=
\frac{m_j}{\sqrt\Ge}\int_{D^e} \Cbb\hatna \Bq_j:\hatna \Bh_k
+\frac{m_k}{\sqrt\Ge}\int_{D^e} \Cbb\hatna \Bq_k:\hatna \Bh_j \\
&\qquad
-\frac{m_j m_k}{\Ge}\int_{D^e} \Cbb\hatna \Bq_j:\hatna \Bq_k
+ \int_{D^e} \Cbb\hatna \Br_j:\hatna\Br_k.
\end{align*}
Then Lemma \ref{cor_betti} yields
\begin{align*}
\Ical_{jk}
&= -\frac{m_j}{\sqrt\Ge} \int_{\p D^e}{\p_\nu  \Bq_j}\cdot \Bh_k
- \frac{m_k}{\sqrt\Ge} \int_{\p D^e}{\p_\nu \Bq_k}\cdot \Bh_j \\
&\qquad +
\frac{m_j m_k}{\Ge}\int_{\p D^e} {\p_\nu  \Bq_j}\cdot \Bq_k
+\int_{D^e} \Cbb\hatna \Br_j:\hatna\Br_k.
\end{align*}

Now, \eqnref{Ical11}-\eqnref{Ical22} follow from Lemma \ref{lem_hS_bdry_int12} and \ref{lem_dhS_hS_int12}. In fact, we have from Proposition \ref{prop_Brj} that
$$
\left| \int_{D^e}\Cbb \hatna \Br_j:\hatna \Br_k \right|
 \lesssim \Ecal[\Br_j]^{1/2} \Ecal[\Br_k]^{1/2} \lesssim 1.
$$
Since $\Bh_j = (-1)^i \frac{1}{2} \Psi_j$ on $\p D_i$, we have
\begin{align*}
\Ical_{jk}
&= -\frac{m_j}{\sqrt\Ge}\sum_{i=1}^2\frac{(-1)^i}{2}\int_{\p D_i}{\p_\nu  \Bq_j}\cdot \Psi_k
\\
&\quad -\frac{m_k}{\sqrt\Ge}\sum_{i=1}^2\frac{(-1)^i}{2}\int_{\p D_i}{\p_\nu  \Bq_k}\cdot \Psi_j
+ \frac{m_j m_k}{\Ge}\int_{\p D^e} {\p_\nu  \Bq_j}\cdot \Bq_k + O(1).
\end{align*}
Then, from \eqnref{qjpsik}, \eqnref{qjqj} and \eqnref{q1q2}, we have
\begin{align*}
\Ical_{11}&=  \frac{m_1^2}{\Ge} \Big(-\frac{\sqrt\Ge}{m_1} + O(\tau\Ge|\ln\Ge|+\Ge)\Big) + 2 \frac{m_1}{\sqrt\Ge} +O(1)
%\\&
= \frac{m_1}{\sqrt\Ge} + O(\tau|\ln\Ge| + 1),
\\
\Ical_{22}&=  \frac{m_2^2}{\Ge} \Big(-\frac{\sqrt\Ge}{m_2} + O(\tau\Ge|\ln\Ge|+\Ge)\Big) + 2 \frac{m_2}{\sqrt\Ge} +O(1)
%\\&
= \frac{m_2}{\sqrt\Ge} + O(\tau|\ln\Ge| + 1),
\\
\Ical_{12}&=  \frac{m_1 m_2}{\Ge} O(\tau\Ge|\ln\Ge|+\Ge) +O(1)=  O(\tau|\ln\Ge| + 1).
\end{align*}
This completes the proof.
\qed

\medskip
\begin{lemma}\label{lem_I13_I23_estim}
We have
\beq\label{Icalj3}
|\Ical_{13}|, |\Ical_{23}| \lesssim 1,
\eeq
and
\beq\label{Ical33}
 \Ical_{33}\approx 1.
\eeq
\end{lemma}
\pf
We prove \eqnref{Ical33} first. For that we closely follow the proof of (4.12) in \cite{BLL-ARMA-15}.
Let $\Bh_3^K$ be the function defined as follows: for $(x,y)\in \Pi_L$,
\beq\label{eqn_h3K_def1}
\Bh^K_3(x,y)=\frac{ x+f_1(y) }{f_2(y)+f_1(y)} \Psi_3 +\frac{(-1)}{2}\Psi_3.
\eeq
We emphasize that
\beq\label{heKv3K}
\Bh_3^K = -\frac{1}{2} \Bv_{13}^K + \frac{1}{2} \Bv_{23}^K,
\eeq
where $\Bv_{i3}$ is defined by \eqnref{BvikK}.
We then extend $\Bh_3^K$ to $D^e \setminus \Pi_{L}$ so that
\begin{align}
\begin{cases}
\ds\Bh^K_3 = (-1)^i\frac{1}{2} \Psi_3 \quad &\mbox{ on } \p D_i, i=1,2,
\\[3mm]
\ds \Bh_3^K|_{\Rbb^2\setminus B_0}=0,
\\[0.5em]
\ds\|\Bh_3^K\|_{H^1(D^e\setminus \Pi_{L})} \lesssim 1,
\end{cases}
\label{eqn_h3K_def2}
\end{align}
where $B_0$ is a disk which contains $\overline{D_1\cup D_2}$.

It is easy to see that, for $(x,y)\in \Pi_{L}$,
\beq \label{eqn_h3K_asymp}
\p_1 h^K_{31} = -\frac{y}{\Ge+\frac{1}{2}(\Gk_1+\Gk_2) y^2} + O(1),
\eeq
and
\beq \label{eqn_h3K_asymp2}
\p_2 h^K_{31}, \,
\p_1 h^K_{32}, \,
\p_2 h^K_{32} = O(1).
\eeq
We mention that these estimates together with Lemma \ref{lem_q3_asymp} show that $\nabla \Bh_3^K$ and $\nabla \Bq_3$ have the same behavior in $\Pi_L$. In fact, we have
\beq\label{q3h3}
|\nabla \Bq_3 - \nabla \Bh_3^K| \lesssim 1 \quad \mbox{in }\Pi_L.
\eeq
This estimate will be used in the proof of Proposition \ref{hGOjandBq}.

By Lemma \ref{lem_var_principle}, we have
\begin{align*}
\Ical_{33}&=\Ecal_{D^e}[\Bh_3] \leq \Ecal_{D^e} [\Bh_3^K]
\lesssim
\int_{D^e} |\nabla\Bh_3^K|^2
\\
&\lesssim\int_{\Pi_{L}} |\nabla\Bh^K_3|^2 + \int_{D^e\setminus \Pi_{L}} |\nabla\Bh^K_3|^2
\\
&\lesssim\int_{-L}^{L}\int_{-f_1(y)}^{f_2(y)} \Big(\frac{|y|}{\Ge+y^2}\Big)^2 \,dx\, dy +1
\\
&\lesssim \int_{-L}^{L} \frac{y^2}{\Ge+y^2} \,dy +1
\lesssim 1,
\end{align*}
where the second to last inequality holds since $f_2(y)+f_1(y) \lesssim \Ge + y^2$.

To prove the opposite inequality, we invoke a result in \cite{BLL-ARMA-15}: For any $\Bv\in H^1(\Pi_{L}\setminus\Pi_{L_0})$ satisfying $\Bv=0$ on $\p D_1 \cap \p (\Pi_{L}\setminus\Pi_{L_0})$, it holds
\beq\label{Korn_ver_h3}
\int_{\Pi_{L}\setminus\overline{\Pi_{L_0}}} |\nabla \Bv|^2
\lesssim \int_{\Pi_{L}\setminus\overline{\Pi_{L_0}}} |\hatna \Bv|^2.
\eeq
(See the proof of (4.12) in \cite{BLL-ARMA-15}.)

Let $\widetilde{\Bh}_3 := \Bh_3 + \frac{1}{2}\Psi_3$. Then $\widetilde{\Bh}_3=0$ on $\p D_1 \cap \p(\Pi_{L}\setminus\Pi_{L_0})$ and
$\widetilde{\Bh}_3=\Psi_3$ on $\p D_2 \cap \p(\Pi_{L}\setminus\Pi_{L_0})$. Therefore, using \eqnref{Korn_ver_h3}, we have
\begin{align}
\Ical_{33} = \Ecal_{D^e}[\Bh_3]=\Ecal_{D^e}[\widetilde{\Bh}_3] \gtrsim  \int_{\Pi_{L}\setminus\overline{\Pi_{L_0}}} |\nabla \widetilde{\Bh}_3|^2 \gtrsim 1.
\end{align}
So, \eqnref{Ical33} is proved.

To prove \eqnref{Icalj3}, let $j=1$ or $2$. From Lemma \ref{cor_betti} and \eqnref{hj_def_decomp}, we have
\begin{align*}
\Ical_{j3}
&=\frac{m_j}{\sqrt\Ge}\int_{D^e}  \hatna \Bq_j:\hatna \Bh_3 +\int_{D^e} \Cbb \hatna \Br_j:\hatna \Bh_3 \\
&=-\frac{m_j}{\sqrt\Ge}\int_{\p D^e}  \p_\nu \Bq_j \cdot\Bh_3 +\int_{D^e} \Cbb \hatna \Br_j:\hatna \Bh_3 \\
&=\frac{1}{2} \frac{m_j}{\sqrt\Ge} \Big(\int_{\p D_1}  \p_\nu \Bq_j \cdot\Psi_3
- \int_{\p D_2}  \p_\nu \Bq_j \cdot\Psi_3\Big)
 +\int_{D^e} \Cbb \hatna \Br_j:\hatna \Bh_3 \\
&=: I+II.
\end{align*}
From Lemma \ref{lem_hS_bdry_int12} (ii) and the fact that $a \approx \sqrt\Ge$, we have
$$
I=
 %(-1)\frac{\Ga^2\pi}{\Ge}\frac{\Gm}{\Gl+2\Gm}
  \frac{m_j}{\sqrt\Ge} \Gd_{2j} {a(-1 + 4\pi \Ga_2 \Gm)} =O(1).
$$
It is clear from Proposition \ref{prop_Brj} and \eqnref{Ical33} that
$$
|II| \lesssim {\Ecal_{D^e}[\Br_j]}^{1/2} {\Ical_{33}}^{1/2}\lesssim 1.
$$
This proves \eqnref{Icalj3}.
\qed

\medskip
\begin{lemma}\label{lem_Jk_Omega_estim}
We have
$$
|\Jcal_{\GO,k}|  \lesssim \| \Bg\|_{C^{1,\Gg}(\p\GO)}, \quad k=1,2,3.
$$
\end{lemma}

Before proving Lemma \ref{lem_Jk_Omega_estim}, let us make a short remark on regularity of $\BH_\GO$.
Recall that $\BH_\GO$ is defined by
$$
\BH_\GO = -\Scal_{\p \GO}\big[\p_\nu \Bu\big|_{\p\GO}\big]+\Dcal_{\p \GO}[\Bg]\quad \mbox{in }\GO.
$$
As shown in \cite{BLL-ARMA-15}, we have
\beq\label{Omega_gradu_bdd}
\| \nabla \Bu\|_{L^\infty(\widetilde\GO \setminus \Pi_{L})} \lesssim \| \Bg\|_{C^{1,\Gg}(\p\GO)}.
\eeq
In particular, we have
$$
\| \p_\Gv \Bu \|_{L^\infty(\p\GO)}\lesssim \| \Bg\|_{C^{1,\Gg}(\p\GO)}.
$$
So, for any $\GO_1$ such that $\ol{\GO_1} \subset \GO$, we have
\beq\label{HOmega_regular}
 \|\BH_\GO \|_{C^2(\GO_1)}\lesssim \| \Bg\|_{C^{1,\Gg}(\p\GO)}.
\eeq
We also have
\beq\label{HOmega_regular2}
 \|\BH_\GO \|_{H^1(\p\GO)} \lesssim \| \Bg\|_{C^{1,\Gg}(\p\GO)}.
\eeq
Importance of these inequalities is that they hold independently of $\Ge$.

\medskip

\noindent{\sl Proof of Lemma \ref{lem_Jk_Omega_estim}}.
Let us first consider the case when $k=1,2$. For simplicity, we assume $\| \Bg\|_{C^{1,\Gg}(\p\GO)}=1$.
Since $\int_{\p D^e} {\p_\Gv \BH_\GO} \cdot \Bh_k=0$, we have
$$
\Jcal_{\GO,k} = \int_{\p D^e} {\p_\Gv\Bh_k} \cdot \BH_\GO
- \int_{\p D^e} {\p_\Gv \BH_\GO} \cdot \Bh_k.
$$
Then \eqnref{hj_def_decomp} yields
\begin{align*}
\Jcal_{\GO,k} &=
\frac{m_k}{\sqrt\Ge}\Big(\int_{\p D^e} {\p_\Gv \Bq_k} \cdot \BH_\GO
- {\p_\nu \BH_\GO} \cdot \Bq_k\Big)
 + \int_{\p D^e} {\p_\nu \Br_k} \cdot \BH_\GO + (-1)\int_{\p D^e} {\p_\nu \BH_\GO} \cdot \Br_k
\\
&:=I_k + II_k + III_k.
\end{align*}

Green's formula for the Lam\'e system and  \eqnref{qonedirac} yield
\begin{align}
&\frac{\sqrt\Ge}{m_1} I_1 = \int_{D_1\cup D_2} \Lcal_{\Gl,\Gm}\Bq_1 \cdot \BH_\GO \nonumber \\
&={ (\BH_\GO(\Bp_1)-\BH_\GO(\Bp_2))\cdot\Be_1 }
- \sum_{j=1}^2 \frac{ \Ga_2 a}{(\Ga_1-\Ga_2)} { (\p_1 \BH_\GO(\Bp_j)\cdot\Be_1+\p_2 \BH_\GO(\Bp_j)\cdot \Be_2)}. \label{m1I1}
\end{align}
Since $a \approx \sqrt\Ge$ and \eqnref{HOmega_regular} holds, we have
$$
\sum_{j=1}^2 \frac{ \Ga_2 a}{(\Ga_1-\Ga_2)} { (\p_1 \BH_\GO(\Bp_j)\cdot\Be_1+\p_2 \BH_\GO(\Bp_j)\cdot \Be_2)} = O(\sqrt{\Ge}).
$$
Since $\Bp_1 = (-a,0)$ and $\Bp_2=(a,0)$, the mean value theorem shows that there is a point, say $\Bp_*$, on the line segment $\ol{\Bp_1\Bp_2}$ such that
$$
|(\BH_\GO(\Bp_1)-\BH_\GO(\Bp_2))\cdot\Be_1| \le 2a |\p_1 \BH_\GO(\Bp_*)\cdot \Be_1|.
$$
So, we have
$$
|(\BH_\GO(\Bp_1)-\BH_\GO(\Bp_2))\cdot\Be_1| \lesssim \sqrt\Ge.
$$
Therefore, from \eqnref{m1I1}, we obtain
$$
I_1 = O(1).
$$
Similarly, one can show
$$
I_2 = O(1).
$$

Next we estimate $II_k$. Let $\Bv\in \Acal^*$ be the solution to the following exterior Dirichlet problem
\beq\label{eqn_vH}
 \ \left \{
 \begin{array} {ll}
\ds \Lcal_{\Gl,\Gm} \Bv= 0 \quad &\mbox{ in } D^e,\\[2mm]
\ds \Bv= {\BH_\GO} \quad &\mbox{ on } \p D^e.
 \end{array}
 \right.
\eeq
From \eqnref{variation}, we have
\beq
\Ecal_{D^e}[\Bv] =\min_{\Bw\in W} \Ecal_{D^e}[\Bw],
\eeq
where
$$
W: = \big\{ \Bw\in \Acal^* : \Bw|_{\p D^e} = \BH_\GO\big\}.
$$
Let $\Bw$ be a function such that
\beq
\begin{cases}
\Bw|_{\Pi_L \cup \p D^e}=\BH_\GO,
\\
\Bw|_{\Rbb^2\setminus B_0}=0,
\\
\| \Bw \|_{H^1(\Rbb^2)} \lesssim  \| \Bg\|_{C^{1,\Gg}(\p\GO)},
\end{cases}
\eeq
where $B_0$ is a disk which contains $\overline{D_1\cup D_2}$. It is worth to mention that the third condition in the above is fulfilled thanks to \eqnref{HOmega_regular}.
Since $\Bw\in W$, we have
\beq\label{vH}
\Ecal_{D^e}[\Bv] \leq \Ecal_{D^e}[\Bw] \lesssim 1.
\eeq
Then, using Lemma \ref{cor_betti} and Proposition \ref{prop_Brj}, we obtain
\begin{align}
|II_k| &=
\left| \int_{\p D^e} {\p_\nu \Br_k} \cdot \Bv \right| \nonumber \\
&= \left| \int_{D^e} \Cbb \hatna \Br_k : \hatna\Bv \right|
\lesssim {\Ecal_{D^e}[\Br_k]}^{1/2}{\Ecal_{D^e}[\Bv]}^{1/2} \lesssim 1, \quad k=1,2.
\label{eqn_JkII_1}
\end{align}

Let us now consider $III_k$. We see from Lemma \ref{lem_qj_far_estim} and \ref{lem_hS_bdry_estim} that
\begin{align*}
|\Br_k|_{\p D_i}| & = \left| \frac{(-1)^i}{2}\Psi_k- \frac{m_k}{\sqrt\Ge}\Bq_k|_{\p D_i} \right|
\lesssim 1.
\end{align*}
So it follows from \eqnref{HOmega_regular2} that
$$
|III_k| \lesssim \| \p_\Gv \BH_\GO \|_{L^2(\p\GO)} \lesssim 1, \quad k=1,2.
$$
Therefore, we have
$$
|\Jcal_{\GO,k}|\leq |I_k|+|II_k|+|III_k| \lesssim 1, \quad k=1,2.
$$

To deal with the case when $k=3$, let $\Bv$ be the solution to \eqnref{eqn_vH} as before. It follows from \eqnref{Ical33} and \eqnref{vH} that
\begin{align*}
|\Jcal_{\GO,3}| &= \left| \int_{\p D^e} {\p_\nu \Bh_3} \cdot \Bv \right|
\\
&= \left| \int_{D^e} \Cbb \hatna \Bh_3 : \hatna\Bv \right|
\lesssim \Ical_{33}^{1/2} \Ecal_{D^e}[\Bv]^{1/2} \lesssim 1.
\end{align*}
The proof is completed.
\qed

\medskip
Proposition \ref{prop_Kcal_Omega_estim} follows from Lemma \ref{prop_I11_I22_I12_asymp}, \ref{lem_I13_I23_estim} and \ref{lem_Jk_Omega_estim}.

%%%%%%%%%%%%%%%%%%%%%%%%%%%%%%%%%%%%%%%%%%%%%%%%%%%%%%%%%%%%%%%%%%%
\subsection{Proof of Proposition \ref{prop_diff_c1_c2_asymp_bdd}}
%%%%%%%%%%%%%%%%%%%%%%%%%%%%%%%%%%%%%%%%%%%%%%%%%%%%%%%%%%%%%%%%%%%%

Set
$$
\Ical:=
\begin{bmatrix}
\Ical_{11} && \Ical_{12} && \Ical_{13} \\
\Ical_{12} && \Ical_{22} && \Ical_{23} \\
\Ical_{13} && \Ical_{23} && \Ical_{33}
\end{bmatrix}.
$$
From  Lemma \ref{prop_I11_I22_I12_asymp} and \ref{lem_I13_I23_estim}, we have
\begin{align}
\mbox{det}\, \Ical &=
{\Ical_{11}}   \big(\Ical_{22} \Ical_{33}- \Ical_{23}^2\big)
-
{\Ical_{12}}  \big(\Ical_{12} \Ical_{33}- \Ical_{13}\Ical_{23}\big)
+
{\Ical_{13}}   \big(\Ical_{12} \Ical_{23}- \Ical_{13}\Ical_{22}\big)
\notag
\\
&= \Ical_{11} \Ical_{22} \Ical_{33} + O(\Ge^{-1/2})
=  \frac{m_1 m_2}{\Ge} \Ical_{33}(1+O(\sqrt\Ge)).
\label{detI_asymp}
\end{align}
So, by \eqnref{Ical33}, the matrix $\Ical$ is invertible for sufficiently small $\Ge$.

By Lemma \ref{lem_represent_c1_c2_bdd} and Cramer's rule, we have
\begin{align}
c_{11}-c_{21} &=
\frac{\Jcal_{\GO,1}}{\mbox{det}\,\Ical}   \big(\Ical_{22} \Ical_{33}- \Ical_{23}^2\big)
-
\frac{\Jcal_{\GO,2}}{\mbox{det}\,\Ical}   \big(\Ical_{12} \Ical_{33}- \Ical_{13}\Ical_{23}\big)
\notag
\\
&\quad +
\frac{\Jcal_{\GO,3}}{\mbox{det}\,\Ical}   \big(\Ical_{12} \Ical_{23}- \Ical_{13}\Ical_{22}\big).
\label{c11_c21_rep}
\end{align}

Recall from
Lemma \ref{prop_I11_I22_I12_asymp}, \ref{lem_I13_I23_estim} and \ref{lem_Jk_Omega_estim} that
$$
\Ical_{11}, \Ical_{22} \approx \Ge^{-1/2}, \quad
|\Ical_{12}| \lesssim 1+\tau |\ln\Ge|,
\quad
|\Ical_{j3}|\lesssim 1, \quad \Ical_{33}\approx 1,
$$
and
$$
|\Jcal_{\GO,j}| \lesssim \|\Bg\|_{C^{1,\gamma}(\p\GO)}, \quad
$$
for $j=1,2,3$.
For simplicity, we may assume $\|\Bg\|_{C^{1,\gamma}(\p\GO)}$=1.
Then, from \eqnref{detI_asymp} and \eqnref{c11_c21_rep}, we can easily check that
\begin{align*}
c_{11}-c_{21}&=
\frac{\Jcal_{\GO,1}}{\mbox{det}\,\Ical}  \Ical_{22}\Ical_{33} -
\frac{\Jcal_{\GO,3}}{\mbox{det}\,\Ical}  \Ical_{13}\Ical_{22}
+ O(\Ge+\tau \Ge|\ln\Ge|).
\end{align*}
Hence, by applying \eqnref{Ical11} and the second equality in \eqnref{detI_asymp}, we obtain
$$
c_{11}-c_{21}
= \frac{\sqrt\Ge}{m_1}\Big(\Jcal_{\GO,1}-\frac{\Jcal_{\GO,3} \Ical_{13}}{\Ical_{33}}\Big) + O(\Ge+\tau \Ge|\ln\Ge|).
$$

Similarly, we have
\begin{align*}
c_{12}-c_{22}&=\frac{\sqrt\Ge}{m_2}\Big(\Jcal_{\GO,2}-\frac{\Jcal_{\GO,3} \Ical_{23}}{\Ical_{33}}\Big) + O(\Ge+\tau \Ge|\ln\Ge|),
\\
c_{13}-c_{23}&= \frac{\Jcal_{\GO,3}}{\Ical_{33}} + O(\sqrt\Ge + \tau \sqrt\Ge|\ln\Ge|).
\end{align*}
Finally, the definition \eqnref{def_Kcal_Omega} of $\Kcal_{\GO,j}$ yield \eqnref{c11c21}-\eqnref{c13c23}. The proof of Proposition \ref{prop_diff_c1_c2_asymp_bdd} is completed.
\qed

%%%%%%%%%%%%%%%%%%%%%%%%%%%%%%%%%%%%%%%%%%%%%%%%%%%%%%%%%%%%%
\subsection{Proof of Proposition \ref{hGOjandBq}}
%%%%%%%%%%%%%%%%%%%%%%%%%%%%%%%%%%%%%%%%%%%%%%%%%%%%%%%%%%%%%

To prove \eqnref{rGOj} we modify the function $\Br_j^K$ introduced in \eqnref{eqn_def_rKj_PiL0}. Let
$\Br_{\GO,j}^K$ be a function in $C^2(\Rbb^2)$  such that
\beq\label{6000}
\begin{cases}
\ds\Br_{\GO,j}^K = \Br_j^K|_{\Pi_{L_0}} \quad & \mbox{in } \Pi_{L_0},
\\[2mm]
\ds\Br_{\GO,j}^K = \frac{(-1)^i}{2} \Psi_j - \frac{m_j}{\sqrt\Ge}\Bq_j\quad &\mbox{on } \p D_i, \ i=1,2,
\\[2mm]
\ds \Br_{\GO,j}^K =-\frac{m_j}{\sqrt\Ge}\Bq_j \quad & \mbox{on } \p \GO,
\end{cases}
\eeq
and
\beq\label{rK_Omega_C2}
\|\Br^K_{\GO,j} \|_{H^1(\Rbb^2\setminus \Pi_{L_0})} \lesssim 1.
\eeq
We emphasize that $\Br_{\GO,j}^K= \Br_{j}^K$ is a linear interpolation in the gap region $\Pi_{L_0}$. Note that $\nabla \Br_j^K|_{\Pi_{L_0}}$ is already estimated in Lemma \ref{lem_rjK_grad_estim}.

Let
\beq\label{eqn_decomp_r_Omega_j}
\Bw_j := \Br_{\GO,j}- \Br^K_{\GO,j}, \quad j=1,2 ,
\eeq
where $\Br_{\GO,j}$ is the function defined by \eqnref{eqn_hOj_qj_rOj_decomp}. We see that the function $\Bw_j$ is the solution to the following problem:
\beq\label{eqn_wj}
 \ \left \{
 \begin{array} {ll}
\ds \Lcal_{\Gl,\Gm} \Bw_j= -\Lcal_{\Gl,\Gm}\Br_{\GO,j}^K \quad &\mbox{ in } \widetilde\GO,\\[2mm]
\ds \Bw_j= 0 \quad &\mbox{ on } \p D_1 \cup \p D_2 \cup \p \GO.
 \end{array}
 \right.
\eeq

The following lemma can be proved by arguments parallel to the proof of Lemma 3.6 in \cite{BLL-ARMA-15}. So, we omit the proof.

\begin{lemma}\label{lem_grad_estim_YYLi}
Let $\Bv$ be a solution to
$\Lcal_{\Gl,\Gm} \Bv= -\Lcal_{\Gl,\Gm}\Bf$ in $\widetilde{\GO}$ with
$\Bv=0$ on $\p D_1\cup \p D_2 \cup \p\GO$, where $\Bf$ is a given function belongs to $C^2(\Rbb^2)$.
Assume that the following conditions hold:
\begin{itemize}
\item[(i)]
The function $\Bv$ satisfies
$$
\int_{\widetilde{\GO}} |\nabla \Bv|^2\lesssim 1.
$$
\item[(ii)] The function $\Bf$ satisfies
$$
|(\Lcal_{\Gl,\Gm}\Bf)(x,y)| \lesssim\frac{1}{\Ge+y^2} \quad \mbox{ for }(x,y)\in \Pi_{L}.
$$
\end{itemize}
Then we have, for $0<L'<L$,
$$
\|\nabla \Bv\|_{L^\infty(\Pi_{L'})}\lesssim 1.
$$
\end{lemma}

\begin{lemma}\label{cor_estim_wj}
For $j=1,2$, let $\Bw_j$ be the solution to \eqnref{eqn_wj}.
Then we have
\beq\label{nablaBw}
|\nabla \Bw_j(\Bx)|\lesssim 1 \quad \mbox{ \rm for }\Bx\in \Pi_{L_0}.
\eeq
\end{lemma}
\pf
It suffices to show that the hypotheses (i) and (ii) of Lemma \ref{lem_grad_estim_YYLi} are fulfilled, namely,
\beq\label{hypo1}
\int_{\widetilde\GO} |\nabla \Bw_j|^2\lesssim 1.
\eeq
and
\beq\label{hypo2}
|(\Lcal_{\Gl,\Gm} \Br_{\GO,j}^K)(x,y)|\lesssim \frac{1}{\Ge+y^2}
\quad \mbox{for } (x,y)\in \Pi_{L}.
\eeq

By the first Korn's inequality, the variational principle and Lemma \ref{lem_rjK_grad_estim},  we have
\begin{align*}
\int_{\widetilde\GO} |\nabla \Bw_j|^2 &\leq
2 \int_{\widetilde\GO} |\hatna \Bw_j|^2
\lesssim\int_{\widetilde\GO} \Cbb \hatna \Bw_j: \hatna \Bw_j
\\
&\lesssim \int_{\widetilde\GO} \Cbb \hatna \Br_{\GO,j}: \hatna \Br_{\GO,j}
+ \int_{\widetilde\GO} \Cbb \hatna \Br_{\GO,j}^K: \hatna \Br_{\GO,j}^K
\\
&\lesssim \int_{\widetilde\GO} \Cbb \hatna \Br_{\GO,j}^K: \hatna \Br_{\GO,j}^K
\\
&\lesssim
\int_{\Pi_{L_0}} |\nabla \Br_{j}^K|^2 + \int_{\widetilde{\GO}\setminus\Pi_{L_0}} |\nabla \Br_{j}^K|^2
\lesssim\int_{-L_0}^{L_0} \int_{-f_1(y)}^{f_2(y)} \frac{\Ge+|y|}{\Ge+y^2} \,dx dy +1
\\
&\lesssim
\int_{-L_0}^{L_0} (\Ge+|y|)\,dy
+1
\lesssim 1.
\end{align*}
So we obtain \eqnref{hypo1}.

We now prove \eqnref{hypo2}. Let $d$, $\phi$ and $\Gn$ be the function defined in \eqnref{eqn_def_d_phi}. It follows from \eqnref{dx_rK_1x} and \eqnref{eqn_dy_rK_1x} that, for $(x,y)\in\Pi_{L_0}$,
\begin{align}
\p_{11} r^K_{\GO,11}=\p_{11} r^K_{11} &=0,
\nonumber\\[0.5em]
\p_{12} r^K_{\GO,11}=\p_{12} r^K_{11} &= \frac{\phi'}{d}-\frac{  \phi d'}{d^2},
\nonumber\\[0.5em]
\p_{22} r^K_{\GO,11}=\p_{22} r^K_{11} &=
\left[
\frac{\phi''}{d} -\frac{2\phi' d' }{d^2} -\frac{\phi d''}{d^2}
+\frac{2 \phi d'^2}{d^3} \right]x
\nonumber\\
&\quad \quad
+\frac{\phi'' f_1}{d} + 2\frac{\phi'f_1'}{d} -2\frac{\phi' f_1 d'}{d^2}  + \frac{\phi f_1''}{d}
\nonumber\\
&\quad\quad
-2 \frac{\phi f_1' d'}{d^2} - \frac{\phi f_1 d''}{d^2}+ 2 \frac{\phi f_1 d'^2 }{d^3}
+\Gn''. \label{dyy_rK_1x}
\end{align}
In addition to \eqnref{eqn_estim_d} and \eqnref{eqn_estim_phi}, we have
\beq\label{eqn_estim_phi_recall_modified}
|\phi''(y)|,|\Gn''(y)|\lesssim \frac{1}{\Ge+y^2}.
\eeq
Then, using \eqnref{eqn_estim_d}, \eqnref{eqn_estim_phi}, \eqnref{dyy_rK_1x}, \eqnref{eqn_estim_phi_recall_modified} and the fact that $|x|\lesssim\Ge+y^2$ for $(x,y)\in\Pi_{L_0}$, we have
\begin{align*}
 |\p_{12} r^K_{\GO,11}| &\lesssim \frac{1}{\Ge+y^2} + \frac{(\Ge+y^2)|y|}{(\Ge+y^2)^2}
\lesssim \frac{1}{\Ge+y^2},
\end{align*}
and
\begin{align*}
|\p_{22} r^K_{\GO,11}| &\lesssim \left[\frac{1}{\Ge+y^2} + \frac{|y|}{(\Ge+y^2)^2}+\frac{\Ge+y^2
}{(\Ge+y^2)^2}+ \frac{(\Ge+y^2)y^2}{(\Ge+y^2)^3}\right] (\Ge+y^2)
\\[0.5em]
&\qquad + \frac{1}{{\Ge}+y^2}\frac{\Ge+y^2}{\Ge+y^2} + \frac{|y|}{\Ge+y^2} +\frac{(\Ge+y^2)|y|}{(\Ge+y^2)^2}+ \frac{\Ge+y^2}{\Ge+y^2}
\\
&\qquad +
  \frac{(\Ge+y^2)y^2}{(\Ge+y^2)^2} + \frac{(\Ge+y^2)(\Ge+y^2)}{(\Ge+y^2)^2} + \frac{(\Ge+y^2)(\Ge+y^2)y^2}{(\Ge+y^2)^3} + \frac{1}{{\Ge}+y^2}
\\
&
\lesssim \frac{1}{{\Ge}+y^2}.
\end{align*}
The proof is completed.
\qed

Now we are ready to prove Proposition \ref{hGOjandBq}.
\smallskip

\noindent{\sl Proof of Proposition \ref{hGOjandBq}}. Let us look into estimates in $\widetilde{\GO}\setminus\Pi_{L}$ first. Let $\Bv_{ij}$ be the function defined in \eqnref{vikyan}. It is proved in \cite{BLL-ARMA-15} that
$$
\| \nabla \Bv_{ij} \|_{L^\infty(\widetilde{\GO}\setminus\Pi_{L})} \lesssim 1, \quad i=1,2, \ \ j=1,2,3.
$$
Since $\Bh_{\GO,j}=-\frac{1}{2}\Bv_{1j}+\frac{1}{2}\Bv_{2j}$, we have
$$
\| \nabla \Bh_{\GO,j}\|_{L^\infty(\widetilde{\GO}\setminus\Pi_{L_0})} \lesssim 1, \quad j=1,2,3.
$$
This estimate together with \eqnref{Bqjest} and \eqnref{Bq3est2} yields the second part of \eqnref{rGOj} and \eqnref{rGO3} on $\widetilde{\GO}\setminus\Pi_{L_0}$.

By \eqnref{Brjest} and the first line of \eqnref{6000}, we have
$$
|\nabla \Br_{\GO,j}^K(\Bx) | = |\nabla \Br_j^K (\Bx)| \lesssim 1+\frac{\tau |y|}{\Ge+y^2} \quad \mbox{for }\Bx\in \Pi_{L_0}.
$$
Then, the first part of \eqnref{rGOj} follows from \eqnref{eqn_decomp_r_Omega_j} and \eqnref{nablaBw}.

The estimate \eqnref{rGO3} on $\Pi_{L}$ follows from \eqnref{eqn_hGOk_v1k_v2k}, \eqnref{Bao2}, \eqnref{heKv3K} and \eqnref{q3h3}. In fact, we have on $\Pi_L$
\begin{align*}
\nabla \Bh_{\GO,3} &=-\frac{1}{2} \nabla \Bv_{13} + \frac{1}{2} \nabla \Bv_{23} \\
&=-\frac{1}{2} \nabla \Bv_{13}^K + \frac{1}{2} \nabla \Bv_{23}^K + O (1) \\
&= \nabla \Bh_3^K + O (1) = \nabla \Bq_3 + O (1).
\end{align*}
This completes the proof.
\qed

%%%%%%%%%%%%%%%%%%%%%%%%%%%%%%%%%%%%%%%%%%%%%%%%%%%%%%%%%%%%%%%%%%%%
\section{Stress concentration-the free space problem}\label{sec:free}
%%%%%%%%%%%%%%%%%%%%%%%%%%%%%%%%%%%%%%%%%%%%%%%%%%%%%%%%%%%%%%%%%%%%

In this section we consider the free space problem \eqnref{elas_eqn_free} and characterize the singular behavior of its solution.

Analogously to $\Jcal_{\GO,j}$ in \eqnref{def_Ical_Jcal}, we define
\beq\label{def_Jcal_free}
\Jcal_{j}=\int_{\p D^e} {\p_\nu \Bh_j} \cdot \BH, \quad j=1,2,3,
\eeq
where $\BH$ is the background solution of the problem \eqnref{elas_eqn_free}. It is worth emphasizing that
$\Jcal_j$ is defined using $\BH$ while $\Jcal_{\GO,j}$ uses $\BH_\GO$. Analogously to $\Kcal_{\GO,j}$ in \eqnref{def_Kcal_Omega}, we define
\beq\label{def_Kcal}
\Kcal_1=\Jcal_1-\frac{\Jcal_3 \Ical_{13}}{\Ical_{33}}, \quad
\Kcal_2=\Jcal_2-\frac{\Jcal_3 \Ical_{23}}{\Ical_{33}},\quad
\Kcal_3=\frac{\Jcal_3}{\Ical_{33}}.
\eeq
Then the constants $\Kcal_j$ are bounded regardless of $\Ge$ (see \eqnref{KGOcaljest}).

The following is the main result of this section

\begin{theorem}\label{main_thm_2_general_free}
Let $\Bu$ be the solution to \eqnref{elas_eqn_free}.
Then we have the following decomposition of $\Bu-\BH$:
\beq\label{Bu-BH=}
(\Bu-\BH)(\Bx) = \Bb(\Bx)- \sum_{j=1}^3 \big({\Kcal}_{j}+s_j\big)\Bq_j(\Bx)
 , \quad\Bx\in D^e,
\eeq
where the constants $s_j$, $j=1,2,3$, satisfy
\beq\label{|s_j|}
|s_j|\lesssim \tau\sqrt\Ge|\ln\Ge| \|\BH \|_{H^1(B)},
\eeq
and the function $\Bb$ satisfies
\beq\label{nablaBb}
\|\nabla \Bb\|_{L^\infty(D^e)} \lesssim \|\BH \|_{H^1(B)}.
\eeq
Here, $B$ is a disk containing $\ol{D_1 \cup D_2}$.
\end{theorem}

By the proof analogous to that of Theorem \ref{cor-bdd}, we can derive the following theorem from Theorem \ref{main_thm_2_general_free}.

\begin{theorem}
It holds that
\beq
 \frac{\sum_{j=1,2}|\Kcal_{j}|}{\sqrt{\Ge}}\lesssim\| \nabla (\Bu-\BH)\|_{L^\infty(D^e)} \lesssim
 \frac{\|\BH \|_{H^1(B)}}{\sqrt{\Ge}}.
\eeq
\end{theorem}

We prove Theorem \ref{main_thm_2_general_free} based on Theorem \ref{main_thm_2_general_bdd}. Let $B$ be a disk containing $\overline{D_1\cup D_2}$. We assume for convenience that the center of $B$ is $0$.
Then the solution $\Bu$ to \eqnref{elas_eqn_free} is the solution to \eqnref{elas_eqn_bdd} with $\GO=B$ and $\Bg = \Bu|_{\p B}$.
So we obtain the following decomposition of the solution $\Bu$ in $B$ by applying Theorem \ref{main_thm_2_general_bdd}:
\beq\label{eqn_u_free_decomp_BR}
\Bu = \Bb_{B} - \sum_{j=1}^3 (\Kcal_{B,j}+s_{B,j}) \Bq_j \quad \mbox{in } D^e \cap B,
\eeq
where the constants $s_{B,j}$ and the function $\Bb_{B}$ satisfy
\beq\label{eqn_s_BR_estim}
|s_{B,j}| \lesssim \tau \sqrt\Ge |\ln\Ge| \| \Bu\|_{C^{1,\Gg}(\p B)}.
\eeq
and
\beq\label{eqn_nabla_b_BR_estim}
\|\nabla\Bb_{B}\|_{L^\infty(D^e\cap B)} \lesssim \| \Bu\|_{C^{1,\Gg}(\p B)},
\eeq

Although \eqnref{eqn_u_free_decomp_BR} looks similar to \eqnref{Bu-BH=},
there are three things to be clarified. First, the coefficient of $\Bq_j$ in \eqnref{eqn_u_free_decomp_BR} is given by $\Kcal_{B,j}$, not by $\Kcal_j$.
Second,  the right-hand sides of \eqnref{eqn_s_BR_estim} and \eqnref{eqn_nabla_b_BR_estim} depend on $\Ge$ since $\|\Bu\|_{C^{1,\Gg}(\p B)}$ does. We need to prove the $\|\Bu\|_{C^{1,\Gg}(\p B)}$ is bounded regardless of $\Ge$.
Third, the decomposition is valid only in $B$, not in the whole region $D^e$. In the following we elaborate on these issues to show that \eqnref{eqn_u_free_decomp_BR}-\eqnref{eqn_nabla_b_BR_estim} actually yield \eqnref{Bu-BH=}-\eqnref{nablaBb}.

\begin{lemma}\label{Kcalj_lem}
$\Kcal_{B,j} = \Kcal_j$ for $j=1,2,3$.
\end{lemma}
\pf
By Green's formula and the fact that $\Bu-\BH \in \Acal$, we have
$$
-\Scal_{\p B}\big[\p_\Gv(\Bu-\BH) |_{\p B}\big](\Bx)+\Dcal_{\p B}[(\Bu-\BH)|_{\p B}](\Bx)
=0, \quad \Bx \in B.
$$
Then, by Green's formula again, we have
\begin{align*}
\BH_{B}(\Bx) &= -\Scal_{\p B}\big[{\p_\nu \Bu}\big|_{\p B}\big](\Bx)+\Dcal_{\p B}[\Bu|_{\p B}](\Bx)
\\
&= -\Scal_{\p B}\big[{\p_\nu \BH}\big|_{\p B}\big](\Bx)+\Dcal_{\p B}[\BH|_{\p B}](\Bx)
=\BH(\Bx),
\end{align*}
for $\Bx\in B$. Therefore, Lemma \ref{Kcalj_lem} follows from \eqnref{def_Ical_Jcal}, \eqnref{def_Kcal_Omega} and \eqnref{def_Kcal}.
\qed

\begin{lemma}\label{lem_uBnorm_estim}
Let $B_1$ be a disk containing $\ol{B}$. We have
\beq\label{uB_norm_estim}
\| \Bu\|_{C^{1,\Gg}(\p B)} \lesssim \| \BH\|_{H^1(B_1)},
\eeq
and
\beq\label{KGOcaljest}
|\Kcal_j| \lesssim \| \BH\|_{H^1(B_1)}, \quad j=1,2,3.
\eeq
\end{lemma}
\pf
By Proposition \ref{freeest}, we have
$$
\| \Bu\|_{C^{1,\Gg}(\p B)} \leq \| \Bu-\BH\|_{C^{1,\Gg}(\p B)} + \| \BH\|_{C^{1,\Gg}(\p B)} \lesssim \| \BH\|_{H^1(B_1)}.
$$
By Proposition \ref{prop_Kcal_Omega_estim} and Lemma \ref{Kcalj_lem}, we have
$$
|\Kcal_{j}|=|\Kcal_{B,j}|\lesssim \| \Bu\|_{C^{1,\Gg}(\p B)}.
$$
So, \eqnref{KGOcaljest} follows from \eqnref{uB_norm_estim}.
\qed

\medskip

\noindent{\sl Proof of Theorem \ref{main_thm_2_general_free}}.
Let $s_j:=s_{B,j}$. Then it follows from \eqnref{eqn_s_BR_estim} and \eqnref{uB_norm_estim} that
\beq\label{eqn_sj_estim}
|s_{j}| \lesssim \tau \sqrt\Ge |\ln\Ge| \| \BH\|_{H^1(B_1)}.
\eeq

Let
\beq\label{Bb}
\Bb := \Bu-\BH+\sum_{j=1}^3 (\Kcal_{j}+s_{j}) \Bq_j.
\eeq
To estimate $\Bb$ in $D^e$, we split the region $D^e$ into $D^e\cap B$ and $D^e\setminus B$.
Using \eqnref{eqn_u_free_decomp_BR}, Lemma \ref{Kcalj_lem} and the fact that $s_j=s_{B,j}$, we have
$$
\Bb = \Big(\Bu+\sum_{j=1}^3 (\Kcal_{B,j}+s_{B,j}) \Bq_j\Big) - \BH =\Bb_{B}-\BH \quad\mbox{in } D^e \cap B.
$$
So, we infer from \eqnref{eqn_nabla_b_BR_estim} and \eqnref{uB_norm_estim} that
\begin{align*}
\|\nabla\Bb\|_{L^\infty(D^e\cap B)}&=
 \|\nabla\Bb_{B}\|_{L^\infty(D^e \cap B)}+ \|  \nabla \BH \|_{L^\infty(D^e \cap B)}
 \\
 &\lesssim \| \Bu\|_{C^{1,\Gg}(\p B)} + \|\nabla \BH \|_{L^\infty(D^e \cap B)}
 \lesssim \| \BH\|_{H^1(B_1)}.
\end{align*}

Let us now consider estimates on the region $D^e\setminus B$. From Proposition \ref{freeest} and \eqnref{Bb}, we see that
\begin{align*}
\| \nabla\Bb\|_{L^\infty(D^e\setminus B)}&
= \Big\|\nabla(\Bu-\BH)+\sum_{j=1}^3 (\Kcal_{j}+s_{j}) \nabla\Bq_j \Big\|_{L^\infty(D^e\setminus B)} \\
&\leq
 \|\nabla(\Bu- \BH) \|_{L^\infty(D^e\setminus B)}
  + \Big\| \sum_{j=1}^3 (\Kcal_{j}+s_{j}) \nabla\Bq_j \Big\|_{L^\infty(D^e\setminus B)}
 \\
 &\leq \| \BH\|_{H^1(B_1)} +  \Big\| \sum_{j=1}^3 (\Kcal_{j}+s_{j}) \nabla\Bq_j \Big\|_{L^\infty(D^e\setminus B)}.
\end{align*}
Lemma \ref{lem_qj_far_estim} and \eqnref{Bq3est2} show that
$$
\| \nabla\Bq_j\|_{L^\infty(D^e\setminus B)} \lesssim 1, \quad j=1,2,3.
$$
Therefore, from \eqnref{KGOcaljest} and \eqnref{eqn_sj_estim}, we obtain
\begin{align*}
\|\nabla\Bb\|_{L^\infty(D^e\setminus B)}
 &\leq \| \BH\|_{H^1(B_1)} +  (1+\sqrt\Ge|\ln\Ge|)\| \BH\|_{H^1(B_1)} \leq \|\BH\|_{H^1(B_1)}.
\end{align*}
So we have
$$
\|\nabla\Bb\|_{L^\infty(D^e)}\lesssim \| \BH\|_{H^1(B_1)},
$$
and the proof of Theorem \ref{main_thm_2_general_free} is completed (with $B$ replaced by $B_1$).
\qed

%%%%%%%%%%%%%%%%%%%%%%%%%%%%%%%%%%%%%%%%%%%%%%%%%%
\section{Symmetric inclusions and optimality of the blow-up rate}\label{sec:symmetric_case}
%%%%%%%%%%%%%%%%%%%%%%%%%%%%%%%%%%%%%%%%%%%%%%%%%%

In this section, we show that \eqnref{Bu-BH=} can be further simplified by assuming some symmetry of the inclusions $D_1$ and $D_2$. More importantly, we show that the blow-up rate $\Ge^{-1/2}$ of $|\nabla\Bu|$ is \emph{optimal} by considering two circular inclusions. The singular functions $\Bq_j$ play an essential role here as well.

Let us first assume that the background field $\BH$ can be decomposed as
\beq\label{H_symm_decomp}
\BH = \BH_e + \BH_o,
\eeq
where $\BH_e=(H_{e1}, H_{e2})^T$ and $\BH_o=(H_{o1}, H_{o2})^T$ respectively have the following symmetric properties:
\beq\label{Hesymm}
\begin{cases}
H_{e1}(x,y)=H_{e1}(x,-y)=-H_{e1}(-x,y), \\
H_{e2}(x,y)=-H_{e2}(x,-y)=H_{e2}(-x,y),
\end{cases}
\eeq
and
\beq\label{Hosymm}
\begin{cases}
H_{o1}(x,y)=-H_{o1}(x,-y)=H_{o1}(-x,y), \\
H_{o2}(x,y)=H_{o2}(x,-y)=-H_{o2}(-x,y).
\end{cases}
\eeq

If $\BH$ is a uniform loading, that is, $\BH(\Bx) = (Ax,B y)^T + C(y, x)^T$ for some real coefficients $A,B$ and $C$, then we may take $\BH_e = (Ax,By)^T$ and $\BH_o = C(y,x)^T$.
%%%%%%%%%%%%%%%%%%%%%%%%%%%%%%%%%%%%%%%%%%%%%%%%
\subsection{Symmetric inclusions}
%%%%%%%%%%%%%%%%%%%%%%%%%%%%%%%%%%%%%%%%%%%%%%%%

Let us assume that $D_1 \cup D_2$ is symmetric with respect to both $x$- and $y$-axes. Then we have the following theorem.

\begin{theorem} \label{thm_symm_charact_free}
Let $\Bu$ be the solution to \eqnref{elas_eqn_free} under the assumption that
$D_1 \cup D_2$ is symmetric with respect to both $x$- and $y$-axes. Then, it holds that
\beq\label{BuBHBq}
(\Bu-\BH)(\Bx) = \Bb(\Bx)+  {\Jcal}_{1}\mathbf{q}_1(\Bx)
+ {\Jcal}_{2}\mathbf{q}_2(\Bx)
 , \quad\Bx\in D^e,
\eeq
where the function $\Bb$ satisfies
\beq\label{nablabH}
|\nabla \Bb(\Bx)| \lesssim \| \BH\|_{H^1(B)} \quad\mbox{for all } \Bx \in D^e.
\eeq
Here, $B$ is a disk containing $\ol{D_1\cup D_2}$. Moreover, if $\BH=\BH_e$, {\it i.e.}, $\BH$ satisfies \eqnref{Hesymm}, then $\Jcal_2=0$; If $\BH$ satisfies \eqnref{Hosymm}, then $\Jcal_1=0$.
\end{theorem}

\pf
Since $D_1$ and $D_2$ are symmetric, the number $\tau$ defined by \eqnref{def_tau} is $0$. So it follows from \eqnref{|s_j|} that $s_j=0$ for $j=1,2,3$.

Now it remains to show that $\Kcal_1=\Jcal_1$, $\Kcal_1=\Jcal_2$ and $\Kcal_3=0$, for which it is enough to show that $\Jcal_3=0$ by the definition \eqnref{def_Kcal} of $\Kcal_j$. Recall that
$$
\Jcal_3=\int_{\p D^e} {\p_\nu \Bh_3}\cdot \BH.
$$
Let $\Bh_3=(h_{31},h_{32})^T$. Thanks to the symmetry of the inclusions and the boundary condition of $\Bh_3$, one can see that the following two functions are also solutions of \eqnref{hj_def}:
$$
\begin{bmatrix} -h_{31}(x,-y) \\ h_{32}(x,-y) \end{bmatrix}, \quad
\begin{bmatrix} -h_{31}(-x,y) \\ h_{32}(-x,y) \end{bmatrix}.
$$
So, by the uniqueness of the solution we see that $\Bh_3$ satisfies the following symmetry:
\beq\label{h3symm}
\begin{cases}
h_{31}(x,y)=-h_{31}(x,-y)=-h_{31}(-x,y), \\
h_{32}(x,y)=h_{32}(x,-y)=h_{32}(-x,y).
\end{cases}
\eeq
The outward normal $\Bn=(n_{1}, n_{2})^T$ to $\p D^e$ satisfies
\beq\label{nsymm}
\begin{cases}
n_{1}(x,y)=n_{1}(x,-y)=-n_1(-x,y), \\
n_{2}(x,y)=-n_{2}(x,-y)=n_{2}(-x,y).
\end{cases}
\eeq
So, the conormal derivative $\Bf:=\p_\Gv \Bh_3$ on $\p D^e$ enjoys the following symmetry:
\beq\label{fsymm}
\begin{cases}
f_1(x,y) = -f_1(x,-y) = -f_1(-x,y), \\
f_2(x,y) = f_2(x,-y) = f_2(-x,y).
\end{cases}
\eeq

Let $\BH = \BH_e + \BH_o$ be the decomposition as in \eqnref{H_symm_decomp}. We write $\Jcal_3$ as
$$
\Jcal_3 = \int_{\p D^e}\Bf\cdot \BH_e + \int_{\p D^e}\Bf\cdot \BH_o=: I+II.
$$
Using the symmetry with respect to the $y$-axis in \eqnref{Hesymm} and \eqnref{fsymm}, we have
\begin{align*}
I&=\int_{\p D_1 } (f_1,f_2)\cdot (H_{e1},H_{e2}) + \int_{\p D_2} (f_1,f_2)\cdot (H_{e1},H_{e2})
\\
&=\int_{\p D_1 } (f_1,f_2)\cdot (H_{e1},H_{e2}) + \int_{\p D_1} (-f_1,f_2)\cdot (-H_{e1},H_{e2})
\\
&=2 \int_{\p D_1} f_1 H_{e1} + f_2 H_{e2}.
\end{align*}
Then, the symmetry with respect to the $x$-axis in \eqnref{Hesymm} and \eqnref{fsymm}, we obtain
\begin{align*}
I&= 2\int_{\p D_1 \cap \{y\geq 0\}} \big( f_1 H_{e1} + f_2 H_{e2} \big) + 2 \int_{\p D_1 \cap \{y< 0\}} \big( f_1 H_{e1} + f_2 H_{e2} \big)
\\
&= 2 \int_{\p D_1 \cap \{y\geq 0\}} \big( f_1 H_{e1} + f_2 H_{e2} \big) + 2 \int_{\p D_1 \cap \{y\geq 0\}} \big( (-f_1) H_{e1} + f_2 (-H_{e2}) \big)
=0.
\end{align*}
By the exactly same way, we can show $II=0$.

Suppose that $\BH$ has the symmetry property \eqnref{Hesymm}. Let $\Bg:=\p_\Gv \Bh_2$. Then \eqnref{htwosymm} and \eqnref{nsymm} show that $\Bg$ has the following symmetry properties:
\beq\label{gsymm}
\begin{cases}
g_1(x,y) = -g_1(x,-y)= g_1(-x,y) , \\
g_2(x,y) = g_2(x,-y)= -g_2(-x,y) .
\end{cases}
\eeq
So one can see as before that
$$
\Jcal_2 = \int_{\p D^e}\Bg\cdot \BH_e =0.
$$
Similarly, one can show that $\Jcal_1=0$ if $\BH$ has the symmetry property \eqnref{Hosymm}. This completes the proof. \qed

\begin{cor}\label{10000}
Under the same hypothesis as in Theorem \ref{thm_symm_charact_free}, we have
\beq
(|{\Jcal}_{1}|+|{\Jcal}_{2}|) \frac{1}{\sqrt{\Ge}}\lesssim  \| \nabla (\Bu-\BH)\|_{L^\infty(D^e)} \lesssim \frac{1}{\sqrt{\Ge}}.
\eeq
\end{cor}

We emphasize that $\Jcal_{1}$ and $\Jcal_{2}$ do depend on $\Ge$. In the next subsection we show that there are some cases such that $1 \lesssim |\Jcal_{1}| + |\Jcal_{2}|$ by considering circular inclusions. It implies that $\Ge^{-1/2}$ is the optimal blow-up rate of the gradient.

%%%%%%%%%%%%%%%%%%%%%%%%%%%%%%%%%%%%%%%%%%%%%%%%%%%%%%%%%%%%%%%%%%%%%%%%%%
\subsection{Circular inclusions and optimality of the blow-up rate}
%%%%%%%%%%%%%%%%%%%%%%%%%%%%%%%%%%%%%%%%%%%%%%%%%%%%%%%%%%%%%%%%%%%%%%%%%%

In this subsection we show that $\Ge^{-1/2}$ is a lower bound on the  blow-up rate of the gradient by considering two circular inclusions under a uniform loading.

\begin{prop}\label{prop:J1_J2_lower}
Suppose that $D_1$ and $D_2$ are disks with the same radius $r_0$ and let $\Bu$ be the solution to \eqnref{elas_eqn_free}.
Let
\beq
\Ga^*= \Ga^*(\Gl,\Gm) := \frac{\Gl+\Gm}{\Gm}.
\eeq
\begin{itemize}
\item[(i)] If $\BH(x,y)=(A x,B y)^T$ with $A\neq 0$, then there are $\Ge_0>0$ and $\Ga_0>0$ independent of $\Ge$ such that for any $(\Gl,\Gm)$ satisfying $\Ga^*(\Gl,\Gm) \le \Ga_0$ and $\Ge \le \Ge_0$ the following holds:
\beq\label{Jcaluniform1}
1\lesssim|\Jcal_1| \quad \mbox{and} \quad \Jcal_2=0.
\eeq

\item[(ii)] If $\BH(x,y) = C(y, x)^T$ with $C\neq 0$, then there are $\Ge_0>0$ and $\Ga_0>0$ independent of $\Ge$ such that for any $(\Gl,\Gm)$ satisfying $\Ga^*(\Gl,\Gm) \le \Ga_0$ and $\Ge \le \Ge_0$ the following holds:
\beq\label{Jcaluniform2}
\Jcal_1=0\quad \mbox{and} \quad 1\lesssim|\Jcal_2|.
\eeq
\end{itemize}
\end{prop}

We emphasize that the condition $\Ga^*(\Gl,\Gm) \le \Ga_0$ can be satisfied even if $\Ga_0$ is small. In fact, the strong convexity condition requires $\Gm>0$ and $\Gl+\Gm>0$. So, by taking negative $\Gl$, the condition is satisfied.

\medskip

\pf
We only prove (i) since (ii) can be proved similarly.

Since $\BH(x,y)=(A x,B y)^T$ satisfies \eqnref{Hesymm}, $\Jcal_{2}=0$ by Theorem \ref{thm_symm_charact_free}.
To prove $1 \lesssim|\Jcal_{1}|$, we define $\Br_1$ by
\beq\label{Br1def}
(1+ \frac{m_1}{\sqrt\Ge} t_1)\Bh_1 = \frac{m_1}{\sqrt\Ge} \Bq_1 + \Br_1,
\eeq
where $t_1$ is the constant appearing in Lemma \ref{lem:Bq_on_circle}. Then $\Br_1$ satisfies $\Lcal_{\Gl,\Gm} \Br_1=0$ in $D^e$ and $\Br_1\in\mathcal{A}^*$. It also satisfies, according to Lemma \ref{lem:Bq_on_circle}, the boundary condition
\beq\label{tilde_Br1_bdry}
\Br_1 = - \frac{m_1 \Ga_2 a}{\sqrt\Ge r_0^2} \begin{bmatrix} x \\ y \end{bmatrix} \quad\mbox{on } \p D_1 \cup \p D_2.
\eeq
Since
$$
\Jcal_1=\int_{\p D^e} {\p_\nu \Bh_1}\cdot \BH,
$$
we may write, using \eqnref{Br1def}, $(1+\frac{m_1}{\sqrt\Ge} t_1)\Jcal_{1}$ as
\begin{align}
(1+\frac{m_1}{\sqrt\Ge} t_1)\Jcal_{1} &=
\int_{\p D^e} \p_\Gv ((1+\frac{m_1}{\sqrt\Ge} t_1)\Bh_1 ) \cdot \BH
-
\int_{\p D^e} \p_\Gv \BH \cdot ((1+\frac{m_1}{\sqrt\Ge} t_1)\Bh_1)
\nonumber
\\
&=
\int_{\p D^e} \p_\Gv (\frac{m_1}{\sqrt\Ge} \Bq_1 + \Br_1 ) \cdot \BH
-
\int_{\p D^e} \p_\Gv \BH \cdot (\frac{m_1}{\sqrt\Ge} \Bq_1 + \Br_1)
\nonumber
\\
&=
\frac{m_1}{\sqrt\Ge} \left(\int_{\p D^e} \p_\Gv \Bq_1 \cdot \BH
- \p_\Gv \BH \cdot \Bq_1 \right)
 + \int_{\p D^e} \p_\Gv \Br_1 \cdot \BH - \int_{\p D^e} \p_\Gv \BH \cdot \Br_1
 \nonumber
\\
&=: I_1+I_2+I_3.
\label{J1_decomp_GlGm}
\end{align}

To estimate $I_1$, we first recall that $m_1 := \big[(\Ga_1 - \Ga_2)\sqrt{2(\Gk_1+\Gk_2)}\,\big]^{-1}$. Since $\Gk_1=\Gk_2=1/r_0$,
we have
\beq\label{monesqrt}
\frac{m_1}{\sqrt\Ge}= \frac{\sqrt{r_0}}{2(\Ga_1-\Ga_2)\sqrt{\Ge}}.
\eeq
Then Green's formula for the Lam\'{e} system and \eqnref{qonedirac} yield
\begin{align*}
I_1 &= \frac{m_1}{\sqrt\Ge} \int_{D_1\cup D_2} \Lcal_{\Gl,\Gm}\Bq_1 \cdot \BH \\
&= \frac{\sqrt{r_0}}{2(\Ga_1-\Ga_2)\sqrt{\Ge}} \left[ (\BH(\Bp_1)-\BH(\Bp_2))\cdot\Be_1
- \sum_{j=1}^2 \frac{\Ga_2 a}{\Ga_1-\Ga_2} (\p_1 \BH(\Bp_j)\cdot\Be_1+\p_2 \BH(\Bp_j)\cdot \Be_2) \right].
\end{align*}
Since $\Bp_1=(-a,0)$, $\Bp_2=(a,0)$ and $\BH(x,y)=(Ax, By)^T$, we arrive at
\beq\label{Ione}
I_1 =- \frac{a\sqrt{r_0}}{2(\Ga_1-\Ga_2)\sqrt{\Ge}} \left[ 2A +
\frac{2\Ga_2(A+B)}{\Ga_1-\Ga_2} \right].
\eeq

Since $\Ga_1=\frac{\Gl+3\mu}{4\pi\mu(\Gl+2\mu)}$ and $\Ga_2=\frac{\Gl+\mu}{4\pi\mu(\Gl+2\mu)}$, we have
\begin{align}
\frac{1}{\Ga_1-\Ga_2} &= {2\pi (\Gl+2\mu)} = 2\pi \Gm (1+ \Ga^*),
\label{Ga1Ga2_Gastar1}
\\
\frac{\Ga_2}{\Ga_2-\Ga_2} &= \frac{\Gl+\mu}{2\mu} = \frac{\Ga^*}{2}.
\label{Ga1Ga2_Gastar2}
\end{align}
Since $a=\sqrt{r_0 \Ge}+O(\Ge^{3/2})$, we have
\beq\label{ar0}
\frac{a\sqrt{r_0}}{\sqrt{\Ge}} = r_0 + O(\Ge).
\eeq
Substituting \eqnref{Ga1Ga2_Gastar1}-\eqnref{ar0} into \eqnref{Ione}, we obtain
\begin{align}
I_1 &= - 2\pi \mu(1+\Ga^*)(r_0 + O(\Ge)) \left[ 2A + \Ga^*(A+B) \right]
\nonumber
\\
&= -2\pi r_0 \mu \big(2A  + O(\Ga^* + \Ge+ \Ga^*\Ge) \big).
\label{lim_JI_GlGm}
\end{align}

To estimate $I_2$, let $B_r$ be the disk of radius $r$ centered at the origin containing $\overline{D_1 \cup D_2}$. Choose $r$ and $R$ so that $r<R$.
Let $\chi$ be a smooth radial function such that $\chi(\Bx)=1$ if $|\Bx| \le r$ and $\chi(\Bx)=0$ if $|\Bx| \ge R$. Let $\Bv:= \chi \BH$. Then, using Green's formula, we obtain
\beq\label{eqn_JII_estim_GlGm1}
|I_2| =
\left| \int_{\p D^e} \frac{\p \Br_1}{\p\Gv} \cdot \Bv \right|
= \left| \int_{D^e} \Cbb \hatna \Br_1 : \hatna\Bv \right|
\leq {\Ecal_{D^e}[\Br_1]}^{1/2}{\Ecal_{D^e}[\Bv]}^{1/2}.
\eeq
Let $\Bw(\Bx):= \chi(\Bx)\Bx$ for $\Bx \in D^e$. Then
\beq
- \frac{m_1 \Ga_2 a}{\sqrt\Ge r_0^2} \Bw = \Br_1 \quad\mbox{on } \p D^e.
\eeq
Then, the variational principle \eqnref{variation} yields
\beq
\Ecal_{D^e}[\Br_1] \le \left(\frac{m_1 \Ga_2 a}{\sqrt\Ge r_0^2} \right)^2 \Ecal_{D^e}[\Bw].
\eeq
We emphasize that the variational principle holds since $\Br_1$ is a solution of the Lam\'e system. So, we arrive at
\beq\label{Itwo}
|I_2| \le \frac{m_1 \Ga_2 a}{\sqrt\Ge r_0^2} \Ecal_{D^e}[\Bw]^{1/2} \Ecal_{D^e}[\Bv]^{1/2}.
\eeq

Note that $\| \nabla \Bw\|^2_{L^2(D^e)} \lesssim 1$ and
$$
\Ecal_{D^e}[\Bw]= \int_{D^e}\mathbb{C}\hatna \Bw: \hatna\Bw\lesssim  (\Gl+2\Gm) \| \nabla \Bw\|^2_{L^2(D^e)} \lesssim \Gl+2\Gm.
$$
Here and throughout this proof, $X \lesssim Y$ implies $X \le C Y$ for some constant $C$ independent of $(\Gl, \Gm)$ and $\Ge$. Similarly, one can see that
$$
\Ecal_{D^e}[\Bv] \lesssim \Gl+2\Gm.
$$
So we infer from \eqnref{Itwo} that
\beq\label{Itwo2}
|I_2| \lesssim \frac{m_1 \Ga_2 a}{\sqrt\Ge r_0^2} (\Gl+2\Gm) .
\eeq
It then follows from \eqnref{monesqrt}, \eqnref{Ga1Ga2_Gastar1}--\eqnref{ar0} that
\begin{align}
|I_2|&\lesssim \frac{\Ga_2}{\Ga_1-\Ga_2} (\Gl+2\Gm) \lesssim \Ga^* \mu(1+\Ga^*) \lesssim \mu \Ga^*.
\label{J_II_estim_GlGm}
\end{align}

Since $\BH=(Ax,By)$, it is easy to see that
$$
|\p_\Gv \BH|_{\p D^e}| \lesssim \Gl+2\Gm.
$$
So, by \eqnref{tilde_Br1_bdry} and the fact that $a\approx \sqrt\Ge$, we see that
\beq\label{J_III_estim_GlGm}
|I_3| \lesssim  m_1 \Ga_2 (\Gl+2\mu) \lesssim \frac{\Ga_2}{\Ga_1-\Ga_2}(\Gl+2\Gm)\lesssim \mu\Ga^*.
\eeq

Combining \eqnref{J1_decomp_GlGm}, \eqnref{lim_JI_GlGm}, \eqnref{J_II_estim_GlGm} and \eqnref{J_III_estim_GlGm}, we have
$$
(1+\frac{m_1}{\sqrt\Ge} t_1)\Jcal_{1} = -4\pi r_0 \mu \big(A  + O(\Ga^* + \Ge+ \Ga^*\Ge) \big).
$$
Recall from \eqnref{tj_estim} that $|t_j| \lesssim (\Ga_1+\Ga_2)\Ge^{3/2}$.
So, using \eqnref{Ga1Ga2_Gastar2}, we have
$$
\left| \frac{m_1}{\sqrt\Ge} t_1 \right| \lesssim \frac{\Ga_1+\Ga_2}{\Ga_1-\Ga_2}\Ge \lesssim (1 + \frac{2\Ga_2}{\Ga_1-\Ga_2})\Ge \lesssim (1+\Ga^*)\Ge.$$
Therefore, we finally arrive at
$$
\Jcal_1 = -\frac{{4\pi r_0} \mu (A+ O(\Ga^*+\Ge+ \Ga^*\Ge))} {1 + O({\Ge}+\Ga^* \Ge)}.
$$

Since $A\neq 0$, there are $\Ga_0$ and $\Ge_0$ such that
$$
1 \lesssim |\Jcal_1|
$$
for all $(\Gl,\Gm)$ satisfying $\Ga^* \le \Ga_0$ and $\Ge \le \Ge_0$. This completes the proof.
\qed

\medskip

Corollary \ref{10000} and Proposition \ref{prop:J1_J2_lower} shows that $\Ge^{-1/2}$ is a lower bound on $\nabla \Bu$ as $\Ge \to 0$. In fact, we have
\beq
\|\nabla \Bu \|_{L^\infty(B \setminus (D_1 \cup D_2))} \approx \Ge^{-1/2}
\eeq
as $\Ge \to 0$, provided that $\Ga^* \le \Ga_0$.
In fact, we have a more refined estimate for $\nabla \Bu(0,0)$.

\begin{theorem}
Let $\Bu$ be the solution to \eqnref{elas_eqn_free} when $D_1$ and $D_2$ are disks with the same radius. Suppose that the Lam\'e parameters $(\lambda,\mu)$ satisfy $\Ga^* \le \Ga_0$. Then the following holds as $\Ge \to 0$.
\begin{itemize}
\item[(i)] If $\BH(x,y)=(A x,B y)$ with $A\neq 0$, then
\beq\label{divest}
|\p_1 u_1 (0,0)| \approx \Ge^{-1/2} \quad\mbox{and}\quad  |\p_2 u_1 (0,0)| + |\p_1 u_2 (0,0)| + |\p_2 u_2 (0,0)| \lesssim 1.
\eeq

\item[(ii)] If $\BH(x,y) = C( y, x)$ with $C\neq 0$, then
\beq\label{shearest}
|\p_1 u_2 (0,0)| \approx \Ge^{-1/2} \quad\mbox{and}\quad  |\p_1 u_1 (0,0)| + |\p_2 u_1 (0,0)| + |\p_2 u_2 (0,0)| \lesssim 1.
\eeq
\end{itemize}
\end{theorem}
\pf
Suppose that $\BH(x,y)=(A x,B y)$ with $A\neq 0$. It then follows from \eqnref{BuBHBq} and \eqnref{Jcaluniform1} that
\beq
\nabla(\Bu-\BH)(\Bx) = \nabla \Bb(\Bx)+  {\Jcal}_{1} \nabla \Bq_1(\Bx).
\eeq
Then \eqnref{p1q11}, \eqnref{nablabH} and \eqnref{Jcaluniform1} yield \eqnref{divest}. (ii) can be proved similarly. \qed

\medskip
The estimates \eqnref{divest} and \eqnref{shearest} yield, in particular,
\beq\label{divest2}
|\nabla \cdot \Bu(0,0)| \approx \Ge^{-1/2} \quad\mbox{and}\quad  |\p_2 u_1 (0,0)| + |\p_1 u_2 (0,0)| \lesssim 1
\eeq
if $\BH(x,y)=(A x,B y)$ with $A\neq 0$, and
\beq\label{shearest2}
|\p_1 u_2 (0,0)| + |\p_2 u_1 (0,0)|\approx \Ge^{-1/2} \quad\mbox{and}\quad   |\nabla \cdot \Bu(0,0)| \lesssim 1
\eeq
if $\BH(x,y) = C( y, x)$ with $C\neq 0$. Note that $\nabla \cdot \Bu$ represents the bulk force while $|\p_1 u_2| + |\p_2 u_1|$ the magnitude of the shear force. These estimates are in accordance with results of numerical experiments in \cite{KK-Conm-16}.

%%%%%%%%%%%%%%%%%%%%%%%%%%%%%%%%%%%%%%%%%%%%%%%%%%%%%%%%%%%%
\section*{Conclusion}
%%%%%%%%%%%%%%%%%%%%%%%%%%%%%%%%%%%%%%%%%%%%%%%%%%%%%%%%%%%%

We investigate the problem of characterizing the stress concentration in the narrow region between two hard inclusions and deriving optimal estimates of the magnitude of the stress in the context of the isotropic linear elasticity. We introduce singular functions which are constructed using nuclei of strain, and then show that they capture precisely the singular behavior of the stress as the distance between two inclusions tends to zero. As consequences we are able to derive an upper bound of the blow-up rate of the stress, namely, $\Ge^{-1/2}$ where $\Ge$ is the distance between two inclusions. We then show that $\Ge^{-1/2}$ is an optimal blow-up rate in the sense that it is also a lower bound on the rate of the stress blow-up in some cases. We show that it is a lower bound in the case when inclusions are disks of the same radii.

To show that $\Ge^{-1/2}$ is a lower bound in the case of circular inclusions, we impose a certain condition on the Lam\'e parameters. This condition does not seem natural and may be removed.
In fact, it is likely, as suggested in numerical experiments in \cite{KK-Conm-16}, that $\Ge^{-1/2}$ is a lower bound without any assumption on Lam\'e parameters if the background field is a uniform loading.  It is quite interesting and challenging to clarify this.

\appendix

%%%%%%%%%%%%%%%%%%%%%%%%%%%%%%%%%%%%%%%%%%%%%%%%%%%%%%%%%%%%%%%%%%%%
\section{The Neumann-Poincar\'e operator and the exterior problem}\label{sec:ext}
%%%%%%%%%%%%%%%%%%%%%%%%%%%%%%%%%%%%%%%%%%%%%%%%%%%%%%%%%%%%%%%%%%%%

In this section we prove Propositions \ref{thm_u_layer_bdd_case} and \ref{thm_u_layer}, and Theorem \ref{thm_ext_diri}. The proofs are based on the layer potential technique.

%%%%%%%%%%%%%%%%%%%%%%%%%%%%%%%%%%%%%%%%%%%%%%%%%%%%%%%%%%%%%
\subsection{The NP operator}
%%%%%%%%%%%%%%%%%%%%%%%%%%%%%%%%%%%%%%%%%%%%%%%%%%%%%%%%%%%%%

Let us begin by reviewing well-known results on the layer potentials on simple closed curves. Let $D$ be a simply connected bounded domain in $\Rbb^2$ with the $C^{1,\Ga}$ ($\Ga >0$) smooth boundary. The co-normal derivative of the single layer potential and the double layer potential satisfy the following jump formulas:
\begin{align}
\p_\Gv\Scal_{\p D} [\BGvf]|_{\pm} (\Bx) &= \left(\pm \frac{1}{2}I + \Kcal_{\p D}^* \right) [\BGvf](x), \quad \Bx \in \p D, \label{singlejump} \\
\Dcal_{\p D} [\BGvf]|_{\pm} (\Bx) &= \left(\mp \frac{1}{2}I + \Kcal_{\p D} \right) [\BGvf](x), \quad \Bx \in \p D, \label{doublejump}
\end{align}
where $\Kcal_{\p D}$ is the boundary integral operator defined by
\beq
\Kcal_{\p D} [\BGvf](x) := \mbox{p.v.} \int_{\p D} \p_{\Gv_{\By}} \BGG  (\Bx-\By) \BGvf(\By) d \Gs(\By), \quad \Bx \in \p D,
\eeq
and $\Kcal_{\p D}^*$ is the adjoint operator of $\Kcal_{\p D}$ on $L^2(\p D)^2$. Here, p.v. stands for the Cauchy pricincipal value. The operators $\Kcal_{\p D}$ and $\Kcal_{\p D}^*$ are called the Neumann-Poincar\'e (NP) operators.

It is known that the operator $-1/2I + \Kcal_{\p D}^*$ is a Fredholm operator of index $0$, it is invertible on $H^{-1/2}(\p D)^2$, and its kernel is of three dimensions (see, for example, \cite{DKV-Duke-88}). It is worth mentioning that the NP operator can be realized as a self-adjoint operator on $H^{-1/2}(\p D)^2$ by introducing a new inner product, and it is polynomially compact (see \cite{AJKKY-arXiv}).

We now consider $D^e = \Rbb^2 \setminus \overline{(D_1 \cup D_2)}$, whose boundary $\p D^e$ consists of two disjoint curves $\p D_1$ and $\p D_2$. To define the NP operator in this case, we consider the solution to \eqnref{elas_eqn_free} in the form of \eqnref{singlerep}, namely,
$$
\Bu(\Bx) = \BH(\Bx) + \Scal_{\p D_1}[\BGvf_1](\Bx) + \Scal_{\p D_2}[\BGvf_2](\Bx).
$$
The boundary condition on $\p D^e$ in \eqnref{elas_eqn_free} amounts to
$$
\p_\Gv \left(\BH + \Scal_{\p D_1}[\BGvf_1] + \Scal_{\p D_2}[\BGvf_2] \right) |_- =0 \quad\mbox{on } \p D^e,
$$
which, according to \eqnref{singlejump}, is equivalent to the following system of integral equations:
$$
\left\{
\begin{array}{rl}
\ds \left(- \frac{1}{2}I + \Kcal_{\p D_1}^* \right) [\BGvf_1] + \p_{\Gv} \Scal_{\p D_2}[\BGvf_2]|_{\p D_1} &= -\p_{\Gv} \BH \quad\mbox{on } \p D_1, \\
\ds \p_{\Gv} \Scal_{\p D_1}[\BGvf_1]|_{\p D_2} + \left(- \frac{1}{2}I + \Kcal_{\p D_2}^* \right) [\BGvf_2] &= -\p_{\Gv} \BH \quad\mbox{on } \p D_2.
\end{array}
\right.
$$
This system of integral equations can be rewritten as
\beq\label{inteqnBH}
\left( -\frac{1}{2} \Ibb + \Kbb^* \right)
\begin{bmatrix} \BGvf_1 \\ \BGvf_2 \end{bmatrix} = - \begin{bmatrix} \p_{\Gv} \BH |_{\p D_1} \\ \p_{\Gv} \BH |_{\p D_2} \end{bmatrix},
\eeq
where $\Ibb$ is the identity operator and $\Kbb^*$, which is the NP operator on $\p D^e$, is defined by
\beq
\Kbb^* := \begin{bmatrix} \Kcal_{\p D_1}^* & \p_{\Gv_1} \Scal_{\p D_2} \\ \p_{\Gv_2} \Scal_{\p D_1} & \Kcal_{\p D_2}^* \end{bmatrix}.
\eeq
A special attention is necessary for the off-diagonal entries in the above: For example, $\p_{\Gv_1} \Scal_{\p D_2}$ means that the single layer potential is defined on $\p D_2$ and the co-normal derivative is evaluated on $\p D_1$, so the operator maps $H^{-1/2}(\p D_2)^2$ into $H^{-1/2}(\p D_1)^2$. One can see that the adjoint operator $\Kbb$ of $\Kbb^*$ on $L^2(\p D^e)^2$ is given by
\beq
\Kbb = \begin{bmatrix} \Kcal_{\p D_1} & \Dcal_{\p D_2}|_{\p D_1} \\ \Dcal_{\p D_1}|_{\p D_2} & \Kcal_{\p D_2} \end{bmatrix}.
\eeq
Here $\Dcal_{\p D_2}|_{\p D_1}$ means the double layer potential on $\p D_2$ evaluated on $\p D_1$. We emphasize that
\beq
(\Dcal_{\p D_1}[\BGvf_1] + \Dcal_{\p D_2}[\BGvf_2])|_+
= \left( -\frac{1}{2} \Ibb + \Kbb \right)
\begin{bmatrix} \BGvf_1 \\ \BGvf_2 \end{bmatrix} \quad\mbox{on } \p D^e.
\eeq

\begin{lemma}
The operator $-1/2 \Ibb + \Kbb^*$ is of Fredholm index $0$ on $H^{-1/2}(\p D^e)^2$.
\end{lemma}

\pf
We express $-1/2 \Ibb + \Kbb^*$ as
\beq\label{comppert}
-1/2 \Ibb + \Kbb^* = \begin{bmatrix} -1/2I + \Kcal_{\p D_1}^* & 0 \\ 0 & -1/2I + \Kcal_{\p D_2}^* \end{bmatrix} +
\begin{bmatrix} 0 & \p_{\Gv_1} \Scal_{\p D_2} \\ \p_{\Gv_2} \Scal_{\p D_1} & 0 \end{bmatrix}.
\eeq
Since $-1/2I + \Kcal_{\p D_j}^*$ is of Fredholm index $0$ for $j=1,2$, so is the first operator on the right-hand side above. Since $\p D_1$ and $\p D_2$ are apart, the second operator on the right-hand side is compact. Since the Fredholm index is invariant under a compact perturbation, $-1/2 \Ibb + \Kbb^*$ is of Fredholm index $0$.
\qed

\medskip

In the following we prove Propositions \ref{thm_u_layer_bdd_case} and \ref{thm_u_layer}, and Theorem \ref{thm_ext_diri}. We prove Proposition \ref{thm_u_layer} first since it is simpler.

%%%%%%%%%%%%%%%%%%%%%%%%%%%%%%%%%%%%%%%%%%%%%%%%%%%%%%%%
\subsection{Proof of Proposition \ref{thm_u_layer}}
%%%%%%%%%%%%%%%%%%%%%%%%%%%%%%%%%%%%%%%%%%%%%%%%%%%%%%%%

We first prove the following lemma.
\begin{prop}\label{appinv}
The operator $-1/2 \Ibb + \Kbb^*$ is invertible on $H^{-1/2}_\Psi(\p D_1) \times H^{-1/2}_\Psi(\p D_2)$.
\end{prop}

\pf
As we see from \eqnref{comppert} that $-1/2 \Ibb + \Kbb^*$ is a compact perturbation of an operator which is invertible on $H^{-1/2}_\Psi(\p D_1) \times H^{-1/2}_\Psi(\p D_2)$. So by the Fredholm alternative it suffices to prove the injectivity of $-1/2 \Ibb + \Kbb^*$.

Suppose that
\beq\label{injectassume}
\left( -\frac{1}{2} \Ibb + \Kbb^* \right)
\begin{bmatrix} \BGvf_1 \\ \BGvf_2 \end{bmatrix} = 0
\eeq
for some $(\BGvf_1 , \BGvf_2) \in H^{-1/2}_\Psi(\p D_1) \times H^{-1/2}_\Psi(\p D_2)$ and let
$$
\Bu(\Bx):= \Scal_{\p D_1}[\BGvf_1](\Bx) + \Scal_{\p D_2}[\BGvf_2](\Bx), \quad \Bx \in \Rbb^2.
$$
Then \eqnref{injectassume} implies that $\Lcal_{\Gl, \Gm} \Bu=0$ in $D_i$ and $\p_\Gv \Bu=0$ on $\p D_i$ for $i=1,2$. So
\beq\label{BuDi}
\Bu= \sum_{j=1}^3 a_{ij} \Psi_j \quad \mbox{in } D_i
\eeq
for some constants $a_{ij}$. Since the single layer potential is continuous across $\p D_i$, we have $\Bu|_+= \sum_{j=1}^3 a_{ij} \Psi_j$ on $\p D_i$. Moreover, by the jump formula \eqnref{singlejump} for the single layer potential, we have
$$
\int_{\p D_i} \p_{\Gv} \Bu|_+ \cdot \Psi_j = \int_{\p D_i} (\BGvf_i + \p_{\Gv} \Bu|_-) \cdot \Psi_j = \int_{\p D_i} \BGvf_i \cdot \Psi_j =0
$$
since $\BGvf_i \in H^{-1/2}_\Psi(\p D_i)$. So $\Bu$ is a solution to \eqnref{elas_eqn_free} with $\BH=0$. It is worth mentioning that the decay condition at $\infty$ is satisfies because $\BGvf_j \in H^{-1/2}_\Psi(\p D_j)$. We then have
$$
\int_{D^e} \Cbb \hatna \Bu:\hatna \Bu = \int_{\p D^e} \p_\Gv \Bu|_+ \cdot \Bu = \sum_{j=1}^3 \sum_{i=1}^2 a_{ij} \int_{\p D_i} \p_\Gv \Bu \cdot \Psi_j =0,
$$
where the last equality follows from \eqnref{int_zero}. Hence $\Bu=0$ in $D^e$. By the jump formula \eqnref{singlejump} for the single layer potential, we have
$$
\BGvf_j = \p_\Gv \Bu|_+ - \p_\Gv \Bu|_- = 0 \quad\mbox{on } \p D_i
$$
for $i=1,2$. This completes the proof. \qed

\medskip
\noindent{\sl Proof of Proposition \ref{thm_u_layer}}. Note that since $\Lcal_{\Gl, \Gm} \BH=0$ in $\Rbb^2$, $\p_{\Gv} \BH|_{\p D_i} \in H^{-1/2}_\Psi(\p D_i)$ for $i=1,2$. So
we solve \eqnref{inteqnBH} for $(\BGvf_1 , \BGvf_2)$ on $H^{-1/2}_\Psi(\p D_1) \times H^{-1/2}_\Psi(\p D_2)$. Then $\Bu$ defined by \eqnref{singlerep} is the solution to \eqnref{elas_eqn_free}.

%%%%%%%%%%%%%%%%%%%%%%%%%%%%%%%%%%%%%%%%%%%%%%%%%%%%%%%%
\subsection{Proof of Proposition \ref{thm_u_layer_bdd_case}}
%%%%%%%%%%%%%%%%%%%%%%%%%%%%%%%%%%%%%%%%%%%%%%%%%%%%%%%%

Let $\Bu$ be the solution to \eqnref{elas_eqn_bdd} and let $\Bf:= \p_\Gv\Bu$ on $\p\GO$. Let $\BH_{\GO}$ be the function defined by \eqnref{eqn_def_H_Omega}. We emphasize that $\BH_\GO(\Bx)$ is defined not only for $\Bx \in \GO$, but also for $\Bx \in \Rbb^2 \setminus \overline{\GO}$. Moreover, one can see from \eqnref{singlejump} and \eqnref{doublejump} that the following holds:
\beq\label{HGOjump}
\BH_\GO|_- - \BH_\GO|_+ = \Bg, \quad \p_\Gv\BH_\GO|_- - \p_\Gv\BH_\GO|_+ = \Bf \quad \mbox{on } \p\GO.
\eeq

Let $(\BGvf_1,\BGvf_2)\in H^{-1/2}_\Psi (\p D_1) \times H^{-1/2}_\Psi(\p D_2)$ be the unique solution to \eqnref{inteqnBH} with $\BH$ replaced by $\BH_\GO$, and let
\beq
\Bv_1(\Bx) = \BH_\GO(\Bx) + \Scal_{\p D_1}[\BGvf_1](\Bx) + \Scal_{\p D_2}[\BGvf_2](\Bx), \quad \Bx \in D^e \setminus \p\GO.
\eeq
Then $\Bv_1$ is a solution to
\beq\label{elas_eqn_bdd2}
 \ \left \{
 \begin{array} {ll}
\ds \Lcal_{\Gl,\Gm} \Bv= 0 \quad &\mbox{ in } D^e \setminus \p\GO,\\[2mm]
\ds \Bv=\sum_{j=1}^3 a_{ij} \Psi_j(\Bx), \quad &\mbox{ on } \p D_i, \quad i=1,2 ,\\[2mm]
\ds \Bv|_- - \Bv|_+ = \Bg, \quad \p_\Gv\Bv|_- - \p_\Gv\Bv|_+ = \Bf \quad & \mbox{on } \p\GO, \\[2mm]
\ds \Bv(\Bx)= O(|\Bx|^{-1}) \quad &\mbox{as } |\Bx| \to \infty,
 \end{array}
 \right.
 \eeq
where the constants $a_{ij}$ are determined by the condition \eqnref{int_zero}.

Let
$$
\Bv_2(\Bx) := \left\{
 \begin{array} {ll}
\ds \Bu(\Bx) \quad & \Bx \in \GO \setminus \overline{D_1 \cup D_2}, \\[2mm]
\ds 0 \quad & \Bx \in \Rbb^2 \setminus \overline{\GO}.
 \end{array}
 \right.
$$
Then $\Bv_2$ is also a solution to \eqnref{elas_eqn_bdd2} with the same $\Bg$ and $\Bf$.

Let $\Bv:=\Bv_1-\Bv_2$. Then $\Bv$ is a solution to \eqnref{elas_eqn_bdd2} with $\Bg=0$ and $\Bf=0$. So, we have
$$
\int_{D^e} \Cbb \hatna \Bv:\hatna \Bv = \int_{\p D^e} \p_\Gv \Bv|_+ \cdot \Bv = \sum_{j=1}^3 \sum_{i=1}^2 c_{ij} \int_{\p D_i} \p_\Gv \Bu \cdot \Psi_j =0,
$$
where the last equality follows from \eqnref{int_zero}. So we infer $\Bv=0$ in $D^e$. In particular, $\Bu=\Bv_1$ in $\GO \setminus \overline{D_1 \cup D_2}$ as desired.
\qed

%%%%%%%%%%%%%%%%%%%%%%%%%%%%%%%%%%%%%%%%%%%%%%%%%%%
\subsection{Proof of Theorem \ref{thm_ext_diri}}
%%%%%%%%%%%%%%%%%%%%%%%%%%%%%%%%%%%%%%%%%%%%%%%%%%%

Let
$$
V = \left\{\Bf \in H^{1/2}(\p D^e)^2: \Kbb[\Bf]=\frac{1}{2}\Bf \right\}
$$
and
$$
W = \left\{ \Bf \in H^{-1/2}(\p D^e) :\Kbb^*[\Bf]=\frac{1}{2}\Bf \right\},
$$
which are null spaces of $-1/2\Ibb + \Kbb$ and $-1/2\Ibb + \Kbb^*$, respectively. In particular, we have $\dim V=\dim W$. For $j=1,2,3$, let
$$
\Ga_j^1(\Bx) = \begin{cases}
 \Psi_j(\Bx) &\quad \mbox{if } \Bx \in \p D_1,
 \\
 0 &\quad \mbox{if } \Bx \in \p D_2,
 \end{cases}
$$
and
$$
\Ga_j^2(\Bx) = \begin{cases}
 0 &\quad \mbox{if } \Bx \in \p D_1, \\
 \Psi_j(\Bx) &\quad \mbox{if } \Bx \in \p D_2.
\end{cases}
$$

\begin{lemma}\label{prop_alpha}
The following holds:
\begin{itemize}
\item[(i)] $\dim V=\dim W=6$.
\item[(ii)] $\{ \Ga_j^1, \Ga_j^2 : j=1,2,3 \}$ is a basis of $V$.
\end{itemize}
\end{lemma}

\pf
If $\Bx \in \Rbb^2 \setminus \overline{D_i}$, then
\beq
\Dcal_{\p D_i} [\Psi_j](\Bx) = \int_{D_i} \Cbb \hatna_\By \BGG(\Bx-\By) : \hatna \Psi_j(\By)=0,
\eeq
for $i=1,2$ and $j=1,2,3$. In particular, we have $(-1/2I + \Kcal_{\p D_i})[\Psi_j]=0$ on $\p D_i$. So, we infer that $\{ \Ga_j^1, \Ga_j^2 \} \subset V$.
Since $\Ga_j^1$ and $\Ga_j^2$, $j=1,2,3$, are linearly independent, we infer $\dim V \ge 6$.

On the other hand, since $-1/2 \Ibb + \Kbb^*$ is a Fredholm operator of index 0, we have $H^{-1/2}(\p D^e)^2 = \mathrm{Range}(-1/2\Ibb + \Kbb^*) \oplus W$. According to Proposition \ref{appinv}, $H^{-1/2}_\Psi(\p D_1) \times H^{-1/2}_\Psi(\p D_2) \subset \mathrm{Range}(-1/2\Ibb + \Kbb^*)$. Since $H^{-1/2}_\Psi(\p D_i)^2$ has co-dimension 3 in $H^{-1/2}(\p D_i)$. we infer that the co-dimension of $\mathrm{Range}(-1/2\Ibb + \Kbb^*)$ in $H^{-1/2}(\p D^e)^2$ is larger than or equals to $6$. So, $\dim W \le 6$. This completes the proof. \qed

\begin{lemma}\label{WPsi}
Let
\beq
W_\Psi := \Big\{ \Bf= (\Bf_1, \Bf_2)\in W : \int_{\p D_1} \Bf_1 \cdot \Psi_j + \int_{\p D_2} \Bf_2 \cdot \Psi_j =0 \, \mbox{ for } j=1,2,3 \Big\}.
\eeq
Then, $\dim W_\Psi = 3$.
\end{lemma}

\pf
For $i=1,2$, define $\Gb_j^i$ by
\beq
\Gb_j^i:= \begin{cases}
\Ga_j^i \quad&\mbox{if } j=1,2, \\
\Ga_3^i + c_i \Ga_1^i + d_i \Ga_2^i \quad&\mbox{if } j=3,
\end{cases}
\eeq
where the constants $c_i$ and $d_i$ are chosen so that
\beq
\la \Gb_j^i, \Gb_l^k \ra =0 \quad\mbox{if } (i,j) \neq (k,l).
\eeq
Then there is an eigenfunction $\Bf_j^i \in W$ such that
\beq
\la \Bf_j^i, \Gb_l^k \ra =
\begin{cases}
1 \quad\mbox{if } (i,j) = (k,l), \\
0 \quad\mbox{if } (i,j) \neq (k,l).
\end{cases}
\eeq
In fact, since $( -1/2\Ibb + \Kbb^* )[\Gb_j^i] \in H^{-1/2}_\Psi(\p D_1) \times H^{-1/2}_\Psi(\p D_2)$,
there is a unique $\Bg = (\Bg_1, \Bg_2) \in H^{-1/2}_\Psi(\p D_1) \times H^{-1/2}_\Psi(\p D_2)$ such that
$$
\left( -\frac{1}{2}\Ibb + \Kbb^* \right)[\Bg] = \left( -\frac{1}{2}\Ibb + \Kbb^* \right)[\Gb_j^i].
$$
Let $\Bf := \Gb_j^i - \Bg$. Then, $\Bf \in W$. Moreover, we have
$$
\la \Bf, \Gb_l^k \ra = \la \Gb_j^i, \Gb_l^k \ra.
$$
So $\Bf_j^i:= \la \Gb_j^i, \Gb_j^i \ra^{-1} \Bf$ is the desired function.

Let for $i=1,2$
$$
\Bg_1^i := \Bf_1^i-c^i \Bf_3^i, \quad \Bg_2^i := \Bf_2^i-d^i \Bf_3^i, \quad \Bg_3^i:= \Bf_3^i.
$$
Then, one can see that
\beq
\la \Bg_j^i, \Ga_l^k \ra =
\begin{cases}
1 \quad\mbox{if } (i,j) = (k,l), \\
0 \quad\mbox{if } (i,j) \neq (k,l).
\end{cases}
\eeq
Then $\Bg_j^1-\Bg_j^2$ ($j=1,2,3$) three linearly independent functions belonging to $W_\Psi$, while $\Bg_j^1+ \Bg_j^2$ ($j=1,2,3$) does not belong to $W_\Psi$. So $\dim W_\Psi =3$. \qed

Define the operator $\Sbb: H^{-1/2}(\p D^e)^2 \to H^{1/2}(\p D^e)^2$ as follows: for $\Bf = (\Bf_1, \Bf_2) \in H^{-1/2}(\p D^e)$ let
\beq\label{defBv}
\Bv(\Bx):= \Scal_{\p D_1}[\Bf_1](\Bx) + \Scal_{\p D_2}[\Bf_2](\Bx),
\eeq
and
\beq
\Sbb[\Bf] := \begin{bmatrix} \Bv|_{\p D_1} \\ \Bv|_{\p D_2} \end{bmatrix}.
\eeq

\begin{lemma}\label{Vdecom}
Let
\beq
V^+ := \mbox{span} \, \{ \Psi_1, \Psi_2, \Psi_3 \}.
\eeq
Then, the following holds:
\begin{itemize}
\item[(i)] $\Sbb$ maps $W$ into $V$, and $\Sbb$ is injective on $W_\Psi$.
\item[(ii)] $V= \Sbb(W_\Psi) \oplus V^+$.
\end{itemize}
\end{lemma}

\pf
Let $\Bf \in W$ and define $\Bv$ by \eqnref{defBv}. Then we have
$$
\p_\Gv \Bv|_-=(-1/2\Ibb + \Kbb^*)[\Bf]=0 \quad\mbox{on } \p D^e.
$$
Since $\Lcal_{\Gl, \Gm} \Bv=0$ in $D_i$ ($i=1,2$), we infer that $\Bv= \sum_j a_{ij}\Psi_j$ on $\p D_i$ for some constants $a_{ij}$. So $\Sbb[\Bf] \in V$.

If further $\Bf \in W_\Psi$, then $\Bv(\Bx)=O(|\Bx|^{-1})$ as $|\Bx| \to 0$. So, if $\Sbb[\Bf]=0$, then $\Bv=0$ in $D^e$, and hence $\Bv=0$ in $\Rbb^2$. Thus we have $\Bf=\p_\Gv \Bv|_+ - \p_\Gv \Bv|_-=0$ on $\p D^e$. This proves (i).

We now show that $\Sbb(W_\Psi) \cap V^+ = \{ 0\}$. In fact, if $\Bf = (\Bf_1, \Bf_2) \in W_\Psi$ satisfies $\Sbb[\Bf]= \sum_{j} a_j \Psi_j$ on $\p D^e$, let
$$
\Bv(\Bx):= \Scal_{\p D_1}[\Bf_1](\Bx) + \Scal_{\p D_2}[\Bf_2](\Bx), \quad \Bx \in \Rbb^2.
$$
Then $\Bv \in \Acal$, and
\begin{align*}
\int_{D^e} \Cbb \hatna\Bv: \hatna\Bv &= -\int_{\p D^e} \p_\Gv \Bv|_+ \cdot \Bv \\
&= -\int_{\p D^e} \Bf \cdot (\sum_{j} a_j \Psi_j) - \int_{\p D^e} \p_\Gv \Bv|_- \cdot (\sum_{j} a_j \Psi_j)=0.
\end{align*}
So, $\Bv =0$ in $D^e$, and hence $\sum_{j} a_j \Psi_j=0$ on $\p D^e$.

Since $\Sbb$ is injective on $W_\Psi$, $\dim \Sbb(W_\Psi)=3$. So $\dim \Sbb(W_\Psi) \oplus V^+=6$. This yields (ii). \qed

\medskip

Since $-1/2\Ibb + \Kbb$ is fredholm, we have $H^{1/2}(\p D^e)^2=\mathrm{Range}(-1/2\Ibb + \Kbb) \oplus V^+$. So we obtain the following proposition.

\begin{prop}\label{last}
$H^{1/2}(\p D^e)^2 =\mathrm{Range}(-1/2\Ibb + \Kbb) \oplus \Sbb(W_\Psi) \oplus V^+$.
\end{prop}

\medskip
\noindent{\sl Proof of Theorem \ref{thm_ext_diri}}. Let $\Bg \in H^{1/2}(\p D^e)^2$. According to the previous proposition, there is $\Bf=(\Bf_1, \Bf_2) \in H^{1/2}(\p D^e)^2$, $\BGvf=(\BGvf_1, \BGvf_2) \in W_\Psi$, and constants $a_1,a_2,a_3$ such that
$$
\Bg= \left(-\frac{1}{2}\Ibb + \Kbb \right)[\Bf] + \Sbb[\BGvf] + \sum_{j=1}^3 a_j \Psi_j.
$$
Then the solution $\Bu$ is given by
$$
\Bu = \sum_{i=1}^2 \left( \Dcal_{\p D_i}[\Bf_i] + \Scal_{\p D_i}[\BGvf_i] \right) +  \sum_{j=1}^3 a_j \Psi_j \quad\mbox{in } D^e.
$$
Note that $\sum_{i=1}^2 \left( \Dcal_{\p D_i}[\Bf_i] + \Scal_{\p D_i}[\BGvf_i] \right) \in \Acal$ by Lemma \ref{lem:Acal}. So $\Bu \in \Acal^*$.

For uniqueness, assume that $\Bu$ and $\Bv$ are two solutions in $\Acal^*$, and let $\Bw:=\Bu-\Bv$. Then $\Bw \in \Acal^*$ and $\Bw=0$ on $\p D^e$. Let $\Bw=\Bw_1 + \Bw_2$ be such that $\Bw_1 \in \Acal$ and $\Bw_2= \sum_{j=1}^3 a_j \Psi_j$. Then by Lemma \ref{cor_betti}, we have
$$
0= \int_{D^e} \Cbb \hatna\Bw:\hatna\Bw= \int_{D^e} \Cbb \hatna\Bw_1:\hatna\Bw_1 .
$$
Since $\Bw_1(\Bx) \to 0$ as $|\Bx| \to \infty$, $\Bw_1=0$ in $D^e$. So, $\sum_{j=1}^3 a_j \Psi_j=0$ on $\p D^e$, which implies $a_j=0$, $j=1,2,3$. This completes the proof. \qed

\end{document}